\documentclass[11pt]{amsart}
\usepackage[a4paper,margin=1.15in]{geometry}
\usepackage{amsmath,amssymb,amsthm,mathtools,mathrsfs,booktabs,array,enumitem,xcolor,needspace}
\usepackage[hidelinks]{hyperref}
\hypersetup{pdftitle={Levi Dynamics and Largest Domains of Discontinuity on Flag Manifolds},pdfauthor={Angel Cano, Marlon M. Lopez-Flores, Jose Seade, Andrew Shaw},pdfsubject={Singular Cartan dynamics, Levi incidence, largest domains of discontinuity}}
\numberwithin{equation}{section}
\newtheorem{theorem}{Theorem}[section]

\newtheorem*{destformula}{Destination formula}
\newtheorem{thmx}{Theorem}
\newtheorem{proposition}[theorem]{Proposition}
\newtheorem{lemma}[theorem]{Lemma}
\newtheorem{corollary}[theorem]{Corollary}
\theoremstyle{definition}
\newtheorem{definition}[theorem]{Definition}
\newtheorem{example}[theorem]{Example}

\theoremstyle{remark}
\newtheorem{remark}[theorem]{Remark}
\newcommand{\C}{\mathbb C}\newcommand{\K}{\mathbb K}\newcommand{\R}{\mathbb R}\newcommand{\PP}{\mathbb P}
\newcommand{\g}{\mathfrak g}\newcommand{\ak}{\mathfrak a}
\newcommand{\kk}{\mathfrak k}\newcommand{\pk}{\mathfrak p}
\DeclareMathOperator{\PSL}{PSL}\DeclareMathOperator{\PGL}{PGL}\DeclareMathOperator{\Ad}{Ad}
\DeclareMathOperator{\diag}{diag}\DeclareMathOperator{\Hom}{Hom}\DeclareMathOperator{\PU}{PU}
\newcommand{\G}{\mathcal G}\newcommand{\Lims}{\mathscr L}\newcommand{\I}{\mathscr I}
\newcommand{\D}{\mathscr D}\DeclareMathOperator{\Eq}{Eq}\newcommand{\Rec}{\mathscr R}
\newcommand{\Kul}{\mathrm{Kul}}\newcommand{\Myr}{\mathrm{Myr}}\newcommand{\CoG}{\mathrm{CoG}}
\title[Levi dynamics and largest domains of discontinuity]{Levi Dynamics and Largest Domains of Discontinuity on Flag Manifolds}
\author{\'Angel Cano}
\address[\'Angel Cano]{Instituto de Matem\'aticas, Unidad Cuernavaca, Universidad Nacional Aut\'onoma de M\'exico, Cuernavaca, Morelos, M\'exico}
\email{angelcano@im.unam.mx}
\author{Marlon M. L\'opez-Flores}
\address[Marlon M. L\'opez-Flores]{Laborat\'orio de Intelig\^encia Artificial, Rob\'otica e Cibern\'etica (LIARC), Instituto Militar de Engenharia, Rio de Janeiro, Brazil}
\email{marlon.michael@ime.eb.br}
\author{Jos\'e Seade}
\address[Jos\'e Seade]{Instituto de Matem\'aticas, Unidad Cuernavaca, Universidad Nacional Aut\'onoma de M\'exico, Cuernavaca, Morelos, M\'exico}
\email{jseade@im.unam.mx}
\author{Andrew Shaw}
\address[Andrew Shaw]{Instituto de Matem\'aticas, Unidad Cuernavaca, Universidad Nacional Aut\'onoma de M\'exico, Cuernavaca, Morelos, M\'exico}
\email{andreshaw@ciencias.unam.mx}

\date{}
\providecommand{\Z}{\mathbb Z}\DeclareMathOperator{\Sym}{Sym}

\begin{document}
\begin{abstract}
Divergent sequences in a semisimple Lie group have a bounded Levi part that is invisible to Cartan projections and limit flags. We show that this Levi part is a genuine dynamical invariant and determines much of the large-scale behaviour of discrete groups on flag manifolds. For every real flag manifold, we compute the complete moving-source destination relation of a divergent sequence, including the Levi contribution. This yields a stratified convergence principle, a boundary compactification encoding the relation, and an exact properness criterion. For discrete groups the relation organises the classical limit sets into a single flag, from Quint's limit set through the recurrent and Kulkarni sets to the Myrberg set, the complement of equicontinuity. Sequences with coalescing limit flags satisfy a diagonal formula, which gives a recurrence criterion for the largest domain of proper discontinuity. For compact twists of Kleinian groups, the criterion becomes an exact incidence condition on the twist closure on every Grassmannian and Lagrangian Grassmannian. Two strongly irreducible subgroups of $\PGL(8,\C)$ with identical marked Cartan projections, limit flags and destination relations on $\C\PP^7$ are distinguished on $\operatorname{Gr}(4,\C^8)$ alone. We also show that marked boundary data recover the compact twist while unmarked data recover only its orbit partition. Further applications treat critical projective dynamics, cusped Hitchin groups and symmetric powers.
\end{abstract}

\maketitle
\section*{Introduction}
\label{sec:introduction}

What do the divergent sequences of a discrete group remember? For a discrete subgroup of a semisimple Lie group acting on a flag manifold, the usual answer is recorded by Cartan projections and limit flags: which root gaps diverge, and along which attracting and repelling flags. A divergent sequence carries one more piece of data, a bounded Levi factor acting between its divergent directions. This article describes the limiting dynamics with this factor included, and asks what it governs: how moving sources are distributed, how the classical limit sets fit together, which points return, and what can be recovered from the limiting relations themselves.

A sharp test is a classical question about properness. When do the non-recurrent points of the action belong to a single domain of proper discontinuity, and how much of the answer is determined by Cartan projections and limit flags? Here the Levi factor turns out to be decisive.

Two discrete subgroups of $\PGL(8,\C)$ can be indistinguishable by their marked Cartan projections, their limit flags and their unmarked destination relations on $\C\PP^7$, and still differ on one Grassmannian. Let $\Gamma_0<\mathrm{SL}_2(\C)$ be a non-elementary free Schottky group, and twist it by dense representations into $\mathrm{Sp}(2)$ and into $\mathrm{SU}(4)$, acting on $\C^2\otimes\C^4=\C^8$. Both twisted groups are strongly irreducible. Every divergent sequence in either group has the same singular Cartan type, and the attracting and repelling $4$-planes of the two groups coincide. Their unmarked destination relations on $\C\PP^7$ are equal, and on every Grassmannian $\operatorname{Gr}(k,\C^8)$ they have the same canonical wall. Yet the existence of a largest domain of proper discontinuity separates them, exactly in the middle dimension:
\begin{center}\small
\begin{tabular}{@{}lccccccc@{}}
\toprule
$k$ & $1$ & $2$ & $3$ & $4$ & $5$ & $6$ & $7$\\\midrule
twist by $\mathrm{SU}(4)$ & yes & yes & yes & yes & yes & yes & yes\\
twist by $\mathrm{Sp}(2)$ & yes & yes & yes & no & yes & yes & yes\\\bottomrule
\end{tabular}
\end{center}
The table records whether the action on $\operatorname{Gr}(k,\C^8)$ has a largest domain; when it does, the two domains coincide (Theorem~\ref{thmC}). The difference comes from the compact Levi part of the limiting dynamics, which moves destinations inside their Schubert cells without changing types or limit flags. The group $\mathrm{SU}(4)$ exchanges any two $2$-planes of $\C^4$, whereas $\mathrm{Sp}(2)$ cannot exchange an isotropic plane with a symplectic one. This article develops the theory that detects such phenomena. It gives an exact description of the limiting dynamics of singular Cartan sequences on flag manifolds, Levi component included, and a recurrence criterion through coalescing sequences.

\medskip
On a flag manifold, a singular Cartan sequence can have several divergent gaps and bounded Levi blocks between them. The limiting flags record the expanding directions; the bounded blocks retain residual Levi motion. To capture the full limiting dynamics, we use the dynamical relation of Frances~\cite{Frances} and Kapovich--Leeb--Porti~\cite[\S5.1]{KLP}, considering every possible limit of nearby source points. For an extracted sequence $(g_m)$, define
\[
\D_{(g_m)}(x)=\{y:\text{there are }x_m\to x\text{ with }g_mx_m\to y\}.
\]
This relation records every destination available to a moving source. The first step is to compute its exact fibres, which yields a stratified convergence principle (Theorems~\ref{thm:dest} and~\ref{thm:S}).

Let $G$ be a connected centre-free real semisimple group and $X=G/Q$. After extraction, write $g_m=k_m\exp(\mu_m)k'_m$ with $k_m\to k$ and $k'_m\to k'$. Let $\Theta$ be the set of simple roots whose Cartan gaps diverge, let $b$ be the limiting bounded Levi factor, and write $P_\Theta=L_\Theta N_\Theta$. On each source stratum, the divergent central part retracts to a Levi flag; denote this retraction by $\rho$. Then
\begin{destformula}
For every $x\in X$,
\[
\D_{(g_m)}(x)=k\,\overline{N_\Theta b\rho(k'x)}.
\]
Every point on the right is realised along the sequence, and the graphs converge to this closed relation. The inverse sequence gives the transposed relation. The fibres form a locally trivial family of translated Schubert varieties, with dimension complementary to the source stratum.
The type, limiting parabolics and relative Levi class are encoded by one boundary point in a regular projective compactification of $G$; the destination relation depends only on this point.
\end{destformula}

For a faithful action on a complex flag variety, Brion identifies the
orbit closure of the diagonal in the Hilbert scheme with the wonderful
compactification of $G$ and describes each limiting graph by translates
of Schubert varieties under the diagonal Levi action
\cite[Theorem~3 and Proposition~10]{r2grBrion2003}. Its support is the
destination relation computed by the formula above; a direct proof is given
in Proposition~\ref{prop:hilbert-support} below. The pointwise formula here also
applies to real flag manifolds with arbitrary restricted root
multiplicities. It constructs the moving sources along the whole
extracted sequence and records how the bounded Levi factor places the
destination over each source. This last dependence is what distinguishes
the two actions in Theorem~\ref{thmC}.

For a discrete group, we take the union of the sequence relations. The resulting relation gives an exact properness criterion and organizes higher-rank limit sets into a flag. Theorem~\ref{thm:L} places Quint's limit set at one end and the Myrberg set at the other, with the recurrent and Kulkarni sets between them. In projective space the Myrberg set is the union of the kernels of the pseudo-projective limits of Cano--Seade~\cite{CanoSeade2010}, and its complement is the equicontinuity region.

For regular groups with minimal limit set, the Bruhat-pair formula identifies the recurrent wall whenever the canonical thickening is fat; its complement is the largest domain of proper discontinuity, without hyperbolicity. Stecker proves that for $\Delta$-Anosov representations of non-elementary word-hyperbolic groups every maximal domain comes from a minimal fat ideal~\cite[Corollary~3.16]{SteckerBalanced}. The balanced-ideal counts are one on $\operatorname{Fl}(\R^3)$, ten on $\operatorname{Fl}(\R^4)$, and two on $\operatorname{Gr}(3,6)$~\cite[\S\S4.1--4.2 and Appendix~A.1]{SteckerBalanced}. Thus maximality on $\operatorname{Fl}(\R^3)$ is already known in the Anosov setting. Our contribution there is the recurrent wall, the full relation $\D_\Gamma$, and the extension beyond hyperbolic groups.

We next state the singular recurrence theorem. Fix a self-opposed type $\Theta$ and write $\mathcal C_\Theta$ for the relative Bruhat positions on $X=G/Q$, with opposition $\tau\mapsto w_0\tau$. The canonical thickening is
\[
\mathcal T_\star=\{\tau\in\mathcal C_\Theta:\tau\le w_0\tau\}.
\]
Call a position set fat if its union with its opposite is all of $\mathcal C_\Theta$. A recurrence seed is a limit flag at which every Levi point in this thickening returns under some escaping sequence whose attracting and repelling $\Theta$-flags coalesce there; the precise condition appears in Definition~\ref{def:recurrence-seed}. Call $\Gamma$ $\Theta$-convergent if every sequence of distinct elements has a subsequence of type $\Theta'\supseteq\Theta$, and let $\Lambda_\Theta\subset G/P_\Theta$ be its set of attracting $\Theta$-flags. For a limit flag $\xi$, let $\mathcal T_\star(\xi)$ be the union of the Bruhat strata indexed by $\mathcal T_\star$, and set $\mathcal T_\star(\Lambda_\Theta)=\bigcup_{\xi\in\Lambda_\Theta}\mathcal T_\star(\xi)$.

\begin{thmx}\label{thmB}
Let $\Gamma<G$ be discrete and $\Theta$-convergent, acting on $X=G/Q$, with $\Theta$ self-opposed. Assume that $\Lambda_\Theta$ is minimal, that $\mathcal T_\star$ is fat, and that $\Lambda_\Theta$ contains a recurrence seed. Then
\[
\mathscr R_\Gamma=\mathcal T_\star(\Lambda_\Theta),
\]
and the complement of this recurrent wall is the largest domain of proper discontinuity.
\end{thmx}

Theorem~\ref{thmB} supplies a sufficient criterion for general singular actions. For compact twists of Kleinian groups, the complete relation gives a converse as well, in every degree. Let $\Gamma_0<\mathrm{SL}_2(\C)$ be discrete and non-elementary, with $-I\notin\Gamma_0$ and $\Lambda\subsetneq\PP^1$, and let $\vartheta:\Gamma_0\to U(r)$ have compact closure $H$. Put $W=\C^2\otimes\C^r$, $E_p=p\otimes\C^r$ for $p\in\PP^1$, and $\Gamma_\vartheta=\{[\gamma\otimes\vartheta(\gamma)]\}$. For $1\le k\le2r-1$, set
\[
 \ell=\min(k,2r-k),\qquad a=\lceil\ell/2\rceil,\qquad b=\lfloor\ell/2\rfloor,
\]
\[
\Omega^\star_k=\{S\in\operatorname{Gr}(k,W):\ \dim(S\cap E_p)<k/2\ \text{for all }p\in\Lambda\}.
\]

\begin{thmx}\label{thmGr}
The action of $\Gamma_\vartheta$ on $\operatorname{Gr}(k,W)$ has a largest domain of proper discontinuity if and only if
\[
 \forall K\in\operatorname{Gr}(a,\C^r),\ P\in\operatorname{Gr}(b,\C^r)
 \quad\exists h\in H:\quad hP+h^{-1}P\subseteq K.
\]
In that case the largest domain is the nonempty set $\Omega^\star_k$. Otherwise two non-recurrent $k$-planes are related by the destination relation of $\Gamma_\vartheta$.
\end{thmx}

This is Theorem~\ref{thm:grk-converse}. When $\ell=2b$, the condition says that $H$ exchanges every pair of $b$-planes. The odd-dimensional condition is weaker than pair reversal: twists with compact closure $\mathrm{SO}(3)$ have a largest domain on $\operatorname{Gr}(3,\C^6)$ but none on $\operatorname{Gr}(2,\C^6)$ (Corollary~\ref{cor:so3-odd}). For $k=2$, pair reversal is equivalent to projective transitivity and separates the strongly irreducible twists with closures $\mathrm{SU}(r)$ and $\operatorname{Sym}^{r-1}\mathrm{SU}(2)$, $r\ge3$ (Corollary~\ref{grlevi:cor:two-actions}). The case $r=4$, $k=4$ gives the opening example.

Let $\Gamma_0$ be a non-elementary free Schottky group, so that $\Lambda\subsetneq\PP^1$, and take $r=4$, so that $W=\C^2\otimes\C^4=\C^8$.

\begin{thmx}\label{thmC}
There are two faithful, discrete, strongly irreducible twists $\Gamma_{\mathrm{Sp}}$ and $\Gamma_{\mathrm{SU}}$ of $\Gamma_0$, with compact closures $\mathrm{Sp}(2)$ and $\mathrm{SU}(4)$, such that:
\begin{enumerate}[label=\textup{(\alph*)},nosep]
\item the two twists have identical Cartan projections under the natural marking; every divergent sequence has Cartan type $(4,4)$, and their attracting and repelling $4$-planes are the same planes $E_p$, $p\in\Lambda$;
\item their complete unmarked destination relations on $\PP(W)=\C\PP^7$ coincide;
\item for every $k\ne4$, both have the same largest domain $\Omega^\star_k$ on $\operatorname{Gr}(k,W)$;
\item on $\operatorname{Gr}(4,W)$, $\Omega^\star_4$ is the largest domain for $\Gamma_{\mathrm{SU}}$, whereas $\Gamma_{\mathrm{Sp}}$ has two non-recurrent $4$-planes related by its destination relation, and hence no largest domain.
\end{enumerate}
\end{thmx}

Neither the marked Cartan projections, nor the limit flags, nor the unmarked destination relation on $\C\PP^7$ distinguishes the two groups, and neither do their largest domains on the other Grassmannians. The difference lies in the compact Levi geometry: $\mathrm{SU}(4)$ exchanges any two $2$-planes of $\C^4$, while $\mathrm{Sp}(2)$ cannot exchange an isotropic plane with a symplectic one. This becomes visible only on the middle Grassmannian.

Strong irreducibility here refers to the defining module $W$. The natural decomposition of $\bigwedge^2W$ has two proper invariant summands, so these examples do not satisfy strong irreducibility in all exterior powers; see Remark~\ref{grlevi:rem:exterior-sip}.

The same mechanism works for Lagrangian Grassmannians, with real orthogonal twists $\vartheta:\Gamma_0\to O(n)$ and $\Gamma_\vartheta<\mathrm{Sp}(2n,\C)$. Under the same proper-limit-set hypothesis, put $j=\lceil n/2\rceil$ and let $\beta$ be the standard symmetric bilinear form on $\C^n$. The action on the Lagrangian Grassmannian of $\C^2\otimes\C^n$ has a largest domain if and only if
\[
 \forall A\in\operatorname{Gr}(j,\C^n)\quad
 \exists h\in H:\quad h(A^{\perp\beta})\subseteq A
\]
(Theorem~\ref{thm:lag-converse}). The criterion is necessary and sufficient in every dimension; for $n$ even the inclusion is an equality. For $n=2$ it becomes a statement about one rotation.

\begin{thmx}\label{thmD}
Let $\Gamma_0<\mathrm{SL}_2(\C)$ be discrete and non-elementary, with $-I\notin\Gamma_0$ and $\Lambda\subsetneq\PP^1$, and let $\vartheta:\Gamma_0\to O(2)$ be a real orthogonal twist with closure $H$. The Lagrangian Grassmannian of $\C^2\otimes\C^2$ is the three-dimensional quadric $Q^3$. The action of $\Gamma_\vartheta$ on $Q^3$ has a largest domain of proper discontinuity if and only if $H$ contains the rotation $R_{\pi/2}$. In that case it is $\{S:S\cap E_p=0\text{ for all }p\in\Lambda\}$; otherwise there are two non-recurrent related Lagrangians.
\end{thmx}

A twist onto $\langle R_{\pi/2}\rangle\cong\Z/4$ already gives a largest domain. The trivial twist does not. The algebraic map $A\mapsto A\otimes I_2$ is the complexification of the real embedding used by Stecker, who obtained infinitely many maximal domains on the real Lagrangian Grassmannian in the Fuchsian case~\cite[\S3.6]{SteckerBalanced}. Theorem~\ref{thmD} identifies the obstruction as the missing rotation. All these groups have the same Cartan type, the same limit planes and the same wall.

Section~\ref{sec:applications} gives three further applications of the destination and recurrence framework. In projective space it gives a critical recurrence theorem: for a discrete group of self-opposed uniform critical type $(a,p)$ with full and critically strongly irreducible critical module, the critical wall is the recurrent set whenever $p=a+1$ or the critical flags contain a regular critical seed, and its complement is the largest domain of proper discontinuity (Theorem~\ref{thm:crit-recurrence}). When $p=a+1$ no return estimate is needed: the diagonal formula behind Theorem~\ref{thmB} shows that coalescence alone produces the recurrence. For cusped Hitchin groups the destination relation can be computed completely, with no Zariski-density condition and beyond the word-hyperbolic Borel--Anosov setting; the largest domain is the complement of the middle osculating wall (Theorem~\ref{thm:cusped}). For the irreducible symmetric powers of a non-elementary discrete subgroup of $\PSL(2,\C)$, acting on degree-$N$ divisors of $\PP^1$, the largest domain is given by the multiplicity threshold $\operatorname{mult}_\xi D<\lceil N/2\rceil$ on the limit set; when the classical limit set is proper, it coincides with the equicontinuity region for $N\le2$ and is strictly larger for $N\ge3$ (Theorem~\ref{thm:C}). No finite generation, convex cocompactness or absence of parabolics is assumed.

The general recurrence criterion does not give a necessary combinatorial condition. Proposition~\ref{prop:diagonal-triples} exhibits a diagonal uniform lattice in $\PSL(2,\C)^3$ whose action on $(\C\PP^1)^3$ has a nonempty largest domain, although the canonical thickening is not fat. Here the largest domain consists of the ordered triples of distinct points. This separates fatness of a position ideal from the geometry of its thickening along a particular limit set.

The four questions of the first paragraph organise the article. How moving sources are distributed is answered by the destination formula and the stratified convergence principle (Sections~\ref{sec:cartan} and~\ref{sec:boundary}). How the classical limit sets fit together is the subject of Sections~\ref{sec:groups} and~\ref{sec:chain}, which assemble the sequence relations into the group relation and the flag from Quint's set to the Myrberg set. What can be recovered from the limiting relations is treated in Section~\ref{sec:levi-recovery}, through the marked and unmarked correspondences of twisted Kleinian groups; this section also proves the joint realisation of boundary data used later. Which points return is the subject of Sections~\ref{sec:domains} and~\ref{sec:recurrence}, on properness domains and on regular and singular recurrence. Section~\ref{grlevi:sec:grassmann-levi} applies this answer to Grassmannians and Lagrangian Grassmannians, where it decides the existence of a largest domain in every degree. Section~\ref{sec:applications} treats critical projective dynamics, cusped Hitchin groups and symmetric powers, and Appendix~\ref{sec:proofcells} gives the root-coordinate proof of the cell lemma.

\section{Cartan dynamics and destination formula}\label{sec:cartan}
\subsection{The classical convergence principle}\label{sec:classical}
Rank-one dynamics explains why the sources must be allowed to move.

\begin{proposition}\label{prop:classical}
Let $g_m\in\PSL(2,\C)$ leave every compact set. After extraction, $g_m=k_m\diag(e^{t_m},e^{-t_m})k'_m$ with $t_m\to\infty$, $k_m\to k$, $k'_m\to k'$ in $\mathrm{PSU}(2)$, and for every $x\in\C\PP^1$
\[
\D_{(g_m)}(x)=\begin{cases}\{k[e_1]\},&x\ne k'^{-1}[e_2],\\ \C\PP^1,&x=k'^{-1}[e_2].\end{cases}
\]
\end{proposition}
Here $\D_{(g_m)}(x)$ consists of limits $g_mx_m$ with $x_m\to x$; see Section~\ref{sec:thmA}.
\begin{proof}
Let $a_m=\diag(e^{t_m},e^{-t_m})$ and write points of $\C\PP^1$ as $[u:1]$ or $[e_1]=[1:0]$. For $v_m\to[u:1]$ with $u\ne0$, the images are $a_mv_m=[e^{2t_m}u_m:1]\to[e_1]$, and the same holds for $v_m\to[e_1]$.

If $v_m\to[e_2]=[0:1]$, write $v_m=[w_m:1]$ with $w_m\to0$, so $a_mv_m=[e^{2t_m}w_m:1]$. Given $s\in\C$, choose $w_m=e^{-2t_m}s$: then $a_mv_m=[s:1]$. To reach $[e_1]$, choose $w_m=e^{-t_m}$, so $w_m\to0$ and $e^{2t_m}w_m=e^{t_m}\to\infty$. Hence $\D_{(a_m)}([e_2])=\C\PP^1$.

For general $g_m$, compose with $k_m\to k$ and $k'_m\to k'$, which converge uniformly on the compact space $\C\PP^1$.
\end{proof}
At the fixed repelling source $[e_2]$, transverse perturbations are
magnified by $e^{2t_m}$. Their destinations fill the closure of its
upper-unipotent orbit, the description that extends to higher rank.

%--------------------------------------------------------------
\subsection{An example with nontrivial Levi data}\label{sec:sl3}
For $a_m=\diag(e^{2t_m},ce^{-t_m},c^{-1}e^{-t_m})$ on $\R\PP^2$,
where $t_m\to+\infty$ and $c>0$, the bounded factor
$\diag(1,c,c^{-1})$ moves the exceptional destination line.

\begin{proposition}\label{prop:sl3}
For this sequence, $\D_{(a_m)}([x_1:x_2:x_3])=\{[e_1]\}$ if $x_1\ne0$, and
\[
\D_{(a_m)}([0:x_2:x_3])=\PP\bigl(\R e_1+\R(0,cx_2,c^{-1}x_3)\bigr).
\]
\end{proposition}
\begin{proof}
If $x_1\neq0$ and $v_m\to x$, then $a_mv_m=[e^{3t_m}v_{m,1}:cv_{m,2}:c^{-1}v_{m,3}]\to[e_1]$.

If $x_1=0$ and $v_m=[\varepsilon_m:x_2^{(m)}:x_3^{(m)}]\to x$, then $a_mv_m=[e^{3t_m}\varepsilon_m:cx_2^{(m)}:c^{-1}x_3^{(m)}]$. Every limit lies on the projective line $\PP(\R e_1+\R(0,cx_2,c^{-1}x_3))$.

Its finite points $[s:cx_2:c^{-1}x_3]$, $s\in\R$, occur for $\varepsilon_m=se^{-3t_m}$ and $x_j^{(m)}=x_j$; its remaining point $[e_1]$ occurs for $\varepsilon_m=e^{-t_m}$, since $e^{3t_m}\varepsilon_m=e^{2t_m}\to\infty$.
\end{proof}
These fibres are orbit closures of the upper block-unipotent group
of type $(1,2)$. When $x_2x_3\ne0$, varying $c$ moves the destination
line while preserving the source stratum and the unipotent group.

%--------------------------------------------------------------
\subsection{The setting}\label{sec:setting-sum}
\paragraph{The group.} $G$ is the identity component of $\mathrm{Aut}(\g)$ for a real semisimple Lie algebra $\g$ (so $G$ has trivial centre). Let $\g=\kk\oplus\pk$ be a Cartan decomposition, $K$ the maximal compact subgroup, and $\ak\subset\pk$ a maximal abelian subspace. Write $\Phi\subset\ak^*$ for the \emph{restricted} root system, with root spaces $\g_\lambda$ of any dimension. The system may be non-reduced, as for $\mathrm{SU}(p,q)$ with $p\ne q$. Let $\Delta$ be a set of simple roots and $\ak^+$ the corresponding closed chamber. Then $G=K\exp(\ak^+)K$.

\paragraph{Types.} For $\Theta\subset\Delta$ put $\ak_\Theta=\bigcap_{\alpha\notin\Theta}\ker\alpha$, $A_\Theta=\exp\ak_\Theta$, and $L_\Theta=Z_G(A_\Theta)$. Let $N_\Theta$ be the group whose Lie algebra is spanned by the root spaces of the positive roots \emph{not} in the span of $\Delta\setminus\Theta$. The associated parabolics are $P_\Theta=L_\Theta N_\Theta$ and its opposite $P_\Theta^-=L_\Theta N_\Theta^-$.

Let $W=N_K(\ak)/Z_K(\ak)$ be the Weyl group, with longest element $w_0$, and let $W_\Theta\subset W$ be the Weyl group of the Levi factor $L_\Theta$, generated by the reflections in the roots of $\Delta\setminus\Theta$. The flag manifold is $X=G/Q$ with $Q=P_\theta$ for some $\theta\subset\Delta$, and $W_\theta$ is defined likewise. \emph{We use $\theta$ for the flag manifold and $\Theta$ for the type of a sequence; they are unrelated.} The extreme cases illustrate the convention: $\Theta=\Delta$ gives $P_\Theta=B$, the minimal parabolic; $\theta=\Delta$ gives the full flag manifold $G/B$; $\theta=\{\alpha\}$ for a single $\alpha$ gives a ``Grassmannian-like'' $G/P_\alpha$.

\paragraph{The projective model.} $\PSL(n+1,\C)$, viewed as a real group, has $\ak$ the real diagonal matrices, restricted roots of type $A_n$ with multiplicity $2$, and $X_\theta$ the complex flag manifolds.
\begin{center}\small
\begin{tabular}{@{}p{.31\linewidth}p{.31\linewidth}p{.28\linewidth}@{}}
\toprule
$\PSL(n+1,\C)$ & General $G$ & $\PSL(2,\C)$\\\midrule
Gaps $\sigma_j/\sigma_{j+1}\to\infty$; block decomposition of $n+1$ & Type $\Theta=\{\alpha:\alpha(\mu(g_m))\to\infty\}$ & single gap, $\Theta=\Delta$\\
Source and target flags & Limit parabolics in the frames $k,k'\in K$ & points $x^-,x^+$\\
Graded isomorphisms of the blocks & Relative class in the Levi fibre $L_\Theta/Z_\Theta$ & trivial flag data\\
$\C\PP^n$ & $X_\theta$ for all $\theta$ & $\C\PP^1$\\\bottomrule
\end{tabular}
\end{center}
The projective case in the first column gives a concrete realisation
of the general Cartan and Levi data.

%--------------------------------------------------------------
\subsection{The type, the Levi class and the compact factors}
The Cartan projection separates into diverging directions and a
bounded component in the Levi factor.

Let $(g_m)$ leave every compact set and write $g_m=k_me^{\mu_m}k'_m$ with $\mu_m\in\ak^+$.
After extraction, $k_m\to k$, $k'_m\to k'$, and each $\alpha(\mu_m)$ converges in $[0,\infty]$. Put $\Theta=\{\alpha\in\Delta:\alpha(\mu_m)\to\infty\}$. Since $\Delta$ is a basis of $\ak^*$ and $(g_m)$ diverges, $\Theta\ne\varnothing$.

\begin{lemma}\label{lem:split}\
Write $\mu_m=\nu_m+\eta_m$ with $\nu_m\in\ak_\Theta$ and $\eta_m$ in the span of the coroots of $\Delta\setminus\Theta$. Then $\eta_m\to\eta$ and $\beta(\nu_m)\to\infty$ for $\beta\in\Theta$. Put $z_m=\exp\nu_m$ (central in $L_\Theta$), $b_m=\exp\eta_m\to b\in L_\Theta$ and $a_m=e^{\mu_m}=z_mb_m$. For a restricted root $\lambda$:
\[
\lambda(\nu_m)\to+\infty\ (\lambda\in\Phi(N_\Theta)),\qquad
\lambda(\nu_m)\to-\infty\ (\lambda\in\Phi(N_\Theta^-)),\qquad
\lambda(\nu_m)=0\ (\lambda\in\Phi(L_\Theta)).
\]
\end{lemma}
\begin{proof}
For $\alpha\notin\Theta$ we have $\alpha(\nu_m)=0$, so $\alpha(\eta_m)=\alpha(\mu_m)$ is bounded; the restrictions of these $\alpha$ to the span of their coroots are linearly independent (the Cartan matrix of the Levi factor is invertible), so $\eta_m$ converges because all these coordinates $\alpha(\mu_m)$ already converge. For $\beta\in\Theta$, $\beta(\nu_m)=\beta(\mu_m)-\beta(\eta_m)\to\infty$. A positive root $\lambda=\sum n_\alpha\alpha$ lies in $\Phi(N_\Theta)$ iff $n_\beta>0$ for some $\beta\in\Theta$, and then $\lambda(\nu_m)=\sum_{\beta\in\Theta}n_\beta\beta(\nu_m)\to\infty$. The other two cases follow in the same way.
\end{proof}
Thus $g_m=k_mz_mb_mk'_m$.

%--------------------------------------------------------------
\subsection{Cells and the stratified retraction}
The opposite parabolic orbits give a finite stratification of $X$.
Within each stratum, contraction by the central factor of $L_\Theta$
has a limit in a Levi flag variety. This limit determines the
unipotent orbit whose closure is the destination.

On a complex flag variety the attracting and repelling strata for a
one-parameter subgroup, together with their retractions to its fixed
components, belong to the decomposition of Bia\l ynicki-Birula
\cite{BB1973}. We use the following real Bruhat description, which
also allows the simple gaps in $\Theta$ to diverge at different rates.
Let $X=G/Q$ and fix representatives $\sigma\in N_K(\ak)$ of the double cosets $W_\Theta\backslash W/W_\theta$, each of minimal length in its entire double coset $W_\Theta\sigma W_\theta$. Since $A\subset\sigma Q\sigma^{-1}$, $A$ fixes each $\sigma Q$.

\begin{lemma}\label{lem:cells}
\begin{enumerate}[label=(\roman*),nosep]
\item $X=\bigsqcup_\sigma O_\sigma$, with $O_\sigma=P_\Theta^-\,\sigma Q/Q$.
\item The fixed point set $X^{A_\Theta}=\bigsqcup_\sigma L_\Theta\sigma Q/Q$ is a disjoint union of compact flag manifolds of $L_\Theta$.
\item For $x=n\,l\,\sigma Q\in O_\sigma$ ($n\in N_\Theta^-$, $l\in L_\Theta$) the limit $\rho(x):=\lim zx=l\,\sigma Q$ exists as $z\to\infty$ in $A_\Theta$ (meaning $\beta(\log z)\to+\infty$ for all $\beta\in\Theta$) and does not depend on how $z\to\infty$.
\item Let $\Sigma$ be the set of roots of the unipotent group $N_\sigma=\sigma N_\theta^-\sigma^{-1}$. The map $n\mapsto n\sigma Q$ is an $A$-equivariant chart of a neighbourhood of $\sigma Q$. Split $\Sigma=\Sigma_+\sqcup\Sigma_0\sqcup\Sigma_-$ according to $\Phi(N_\Theta)$, $\Phi(L_\Theta)$, $\Phi(N_\Theta^-)$, and put $M_+=\exp\bigoplus_{\Sigma_+}\g_\lambda$, $M_-=\exp\bigoplus_{\Sigma_0\cup\Sigma_-}\g_\lambda$. Both are subgroups, $N_\sigma=M_-M_+$ with unique decomposition, and $N_\Theta\,\sigma Q=M_+\,\sigma Q$.
\end{enumerate}
\end{lemma}
On each cell $O_\sigma$, the restriction of $\rho$ is a smooth retraction onto $L_\Theta\sigma Q/Q$, holomorphic when $G$ is complex. The map $\rho$ need not be continuous across cells; throughout, the term ``retraction'' refers to these restrictions.
Contraction of $N_\Theta^-$ gives (iii). The root-coordinate proof,
including double restricted roots and the ordered factorisation in (iv),
is given in Appendix~\ref{sec:proofcells}.

\paragraph{The projective plane.} In Proposition~\ref{prop:sl3},
$\Theta=\{\alpha_1\}$ and there are two double cosets, represented by
$1$ and $s$ with $sQ=[e_2]$. The cell $O_1=\{x_1\ne0\}$ retracts
to $[e_1]$, whereas $O_s=\{x_1=0\}$ is fixed pointwise by $\rho$.
At $sQ$, the two roots of $\Sigma$ lie in $\Phi(N_\Theta)$ and
$\Phi(L_\Theta)$ respectively, explaining the one-dimensional
destination over the exceptional line.

%--------------------------------------------------------------
\subsection{The destination formula}\label{sec:thmA}
For a sequence $f_m$ of homeomorphisms of a compact space $X$, let
\[
\D_{(f_m)}(x)=\{y:\ \exists\,x_m\to x\text{ with }f_mx_m\to y\},
\]
where $m$ runs over the \emph{whole} given sequence. Extraction is made once and for all at the level of the Cartan data.

Allowing the sources to move brings us to the theory of dynamically related
points of Frances~\cite{Frances}. Cano--Liu--L\'opez define a destination set using accumulation points
and compute its line fibres for complex hyperbolic groups
\cite[arXiv:1506.08113v5, Definition~3.2 and Proposition~3.6]{CLL}.
M\'endez studies the sequence and group relations, shows how they behave under
converging compact factors, and analyses them on a complex quadric
\cite[arXiv:1809.06953v1, \S4, Proposition~4.1 and Propositions~4.3--4.6]{Mendez2020}.
We now compute these fibres on every flag manifold, with no further extraction
depending on the source or its destination.

\Needspace{9\baselineskip}
\begin{theorem}\label{thm:dest}\
For a sequence $g_m=k_ma_mk'_m$ of type $\Theta$ with data $(k,k',b)$ as in Lemma~\ref{lem:split}, and every $x\in X$:
\begin{enumerate}[label=(\alph*),nosep]
\item $\D_{(a_m)}(x)=\overline{N_\Theta\cdot b\,\rho(x)}$; every point of the closure is reached along the whole sequence $(a_m)$.
\item $\D_{(g_m)}(x)=k\,\D_{(a_m)}(k'x)=k\,\overline{N_\Theta\,b\rho(k'x)}$.
\item If $k'x\in O_\sigma$ the fibre is the closure of a cell of real dimension $\sum_{\lambda\in\Sigma_+}\dim\g_\lambda$.
\item The graphs of $g_m$ converge in the Hausdorff topology to the closed relation $R=\{(x,y):y\in\D_{(g_m)}(x)\}$, and $\D_{(g_m^{-1})}$ is the transpose of $\D_{(g_m)}$.
\end{enumerate}
\end{theorem}

\begin{proof}
At $\sigma Q$, small source perturbations in the bounded and
contracting root directions have vanishing effect on the image.
The expanding directions produce the unipotent orbit, and a diagonal
choice reaches its closure. We then transport this calculation to an
arbitrary source, where the bounded Levi factor reappears.

\emph{Step 1: the base point.} Take $x=\sigma Q$. A source $x_m\to\sigma Q$ is, by Lemma~\ref{lem:cells}(iv), $x_m=q_mp_m\sigma Q$ with $q_m\in M_-$, $p_m\in M_+$, and $q_m,p_m\to1$. On each root space of $\mathrm{Lie}(M_-)$, $\Ad(a_m)$ is multiplication by $e^{\lambda(\mu_m)}$, which is bounded (for $\Sigma_0$) or tends to $0$ (for $\Sigma_-$). There are only finitely many root spaces, so these restrictions have a common operator-norm bound; since $\log q_m\to0$, we have $\Ad(a_m)\log q_m\to0$. Hence $c_m:=a_mq_ma_m^{-1}\to1$, and since $a_m$ fixes $\sigma Q$,
\[
a_mx_m=c_m\,(a_mp_ma_m^{-1})\,\sigma Q.
\]
The second factor lies in $M_+\sigma Q=N_\Theta\sigma Q$, and $c_m\to1$ acts on the compact space $X$ uniformly close to the identity. So every limit point of $a_mx_m$ lies in $\overline{N_\Theta\sigma Q}$. This is true for any subsequence.

Conversely, let $y\in\overline{N_\Theta\sigma Q}$ and choose $p^{(j)}\in M_+$ with $p^{(j)}\sigma Q\to y$. For fixed $j$, $u^{(j)}_m:=a_m^{-1}p^{(j)}a_m\to1$ as $m\to\infty$, because $a_m^{-1}$ multiplies every root space of $\mathrm{Lie}(M_+)$ by $e^{-\lambda(\mu_m)}\to0$. Choose integers $M_1<M_2<\cdots$ with $d(u_m^{(j)},1)<1/j$ for $m\ge M_j$, and put $x_m=u^{(j)}_m\sigma Q$ for $M_j\le m<M_{j+1}$. Then $x_m\to\sigma Q$ and $a_mx_m=p^{(j)}\sigma Q\to y$. This is a construction along the whole sequence. Hence $\D_{(a_m)}(\sigma Q)=\overline{N_\Theta\sigma Q}$.

\emph{Step 2: an arbitrary point.} Let $x=h\sigma Q$ with $h=nl$, $n\in N_\Theta^-$, $l\in L_\Theta$. Since $z_m$ is central in $L_\Theta$ and contracts $N_\Theta^-$,
\[
a_mha_m^{-1}=(a_mna_m^{-1})(b_mlb_m^{-1})\longrightarrow c:=blb^{-1}.
\]
Sources $x_m\to x$ correspond to $h^{-1}x_m\to\sigma Q$, and $a_mx_m=(a_mha_m^{-1})\,a_m(h^{-1}x_m)$ with the first factor converging uniformly to $c$. So $\D_{(a_m)}(x)=c\,\D_{(a_m)}(\sigma Q)=c\,\overline{N_\Theta\sigma Q}$. Now $c\in L_\Theta$ normalises $N_\Theta$, and $b^{-1}$ fixes $\sigma Q$ (as $b\in A$), so
\[
c\,\overline{N_\Theta\sigma Q}=\overline{N_\Theta\,c\,\sigma Q}=\overline{N_\Theta\,b\,l\,\sigma Q}=\overline{N_\Theta\,b\rho(x)},
\]
which proves~(a).

\emph{Step 3: compact factors.} With $v_m=k'_mx_m$ we have $x_m\to x$ iff $v_m\to k'x$, and $g_mx_m=k_ma_mv_m$; since $k_m\to k$ uniformly, this gives (b).

\emph{Step 4: dimension and graphs.} The orbit $N_\Theta b\rho(x)$ is the translate by $bl\in L_\Theta$ of $N_\Theta\sigma Q=M_+\sigma Q$, which is a cell of dimension $\sum_{\Sigma_+}\dim\g_\lambda$ by the chart. This orbit is semialgebraic by the polynomial nilpotent exponential and the algebraic flag action, so its closure has the same dimension.

For (d), let $R^*$ be the set of cluster points of the graphs along arbitrary subsequences. In the compact metric space $X\times X$ it is closed, since $R^*=\bigcap_{M\ge1}\overline{\bigcup_{m\ge M}\operatorname{graph}(g_m)}$. The upper inclusion in Steps~1--3 applies to every subsequence, so $R^*\subset R$; the whole-sequence construction gives $R\subset R^*$. Thus $R$ is closed.

If the graphs failed to lie in an $\varepsilon$-neighbourhood of $R$, compactness would give a cluster point outside $R$. Conversely, take a finite $\varepsilon$-net in $R$ and use the whole-sequence construction at its finitely many points: all are approached for every sufficiently large $m$. These two inclusions give Hausdorff convergence.

Finally, $x_m\to x$ and $g_mx_m\to y$ iff $y_m=g_mx_m\to y$ and $g_m^{-1}y_m\to x$, which gives the transpose.
\end{proof}

\begin{remark}\label{rem:destA}
The relation depends only on $(\Theta,k,k',b)$, not on the relative rates of the divergent gaps. With $g_m=k_ma_mk'_m$, the source is carried by $k'$, so the relevant point is $k'x$. The destination is the closure of a strong unstable leaf, which need not be a manifold. No discreteness is assumed here; it enters when sequence relations are assembled into a group relation in Section~\ref{sec:groups}.
\end{remark}

%--------------------------------------------------------------
\subsection{Dimensions and the rank-one case}\label{sec:consequences}
\begin{corollary}\label{cor:conservation}
If $k'x\in O_\sigma$ then
\[
\dim_{\R}\D_{(g_m)}(x)+\dim_{\R}O_\sigma=\dim_{\R}X.
\]
In particular $\D_{(g_m)}(x)$ is a single point if and only if $k'x$ lies in the open cell, and the set $\mathcal E$ of \emph{exceptional} sources (those with non-trivial destination) is $k'^{-1}$ of the complement of the open cell, a proper algebraic subvariety.
\end{corollary}
\begin{proof}
The chart in Lemma~\ref{lem:cells}(iv) gives $\dim X=\sum_\Sigma\dim\g_\lambda$. The tangent-space calculation in Appendix~\ref{sec:proofcells} gives $\dim O_\sigma=\sum_{\Sigma_0\cup\Sigma_-}\dim\g_\lambda$, and Theorem~\ref{thm:dest}(c) gives the fibre dimension $\sum_{\Sigma_+}\dim\g_\lambda$. The complement of the open source cell is a finite union of opposite Schubert varieties, hence is algebraic.
\end{proof}
The equality concerns \emph{dimension}; the orbit closure can still be singular.

\begin{corollary}\label{cor:rank-one}
If $G$ has real rank one and $P$ is a proper parabolic subgroup, then $X=G/P$ is the sphere at infinity of the symmetric space, $\Theta=\Delta$, there are two cells (a point $\{p_-\}$ and the open cell), and for every divergent sequence with Cartan data as in Lemma~\ref{lem:split}
\[
\D_{(g_m)}(x)=\begin{cases}\{k\,p_+\},&k'x\ne p_-,\\ X,&k'x=p_-.\end{cases}
\]
\end{corollary}
\begin{proof}
Here $Q=P_\Theta$ is the minimal parabolic, $p_+=eQ$ is fixed by $N$, and $p_-=wQ$ is the fixed point of $P^-$. The cells are $O_1=P^-eQ=X\setminus\{p_-\}$ (open) and $O_w=\{p_-\}$. The group $N$ acts simply transitively on $X\setminus\{p_+\}$, so $\rho(x)=p_+$ for $x\ne p_-$ and $\overline{N\,b\,p_-}=X$. In the normalisation of Lemma~\ref{lem:split}, $\Theta=\Delta$ gives $\eta_m=0$, so $b=1$.
\end{proof}
This is the boundary convergence property in rank one. On the perfect boundary of a non-elementary hyperbolic group, the same formula follows along a convergence subsequence: use $x_m=g_m^{-1}y$ for $y\ne x^+$, and targets $y_{j(m)}\to x^+$ with $x_m=g_m^{-1}y_{j(m)}$ for the exceptional destination. In higher rank, divergent root directions and a residual Levi class replace this two-point picture.

%--------------------------------------------------------------
\subsection{Projective space}\label{sec:proj}
For $\PSL(n+1,\mathbb K)$ acting on $\mathbb K\PP^n$ ($\mathbb K=\R$ or $\C$), vector coordinates give a direct proof of Theorem~\ref{thm:dest} and a model for later applications.

Let $V=\mathbb K^{n+1}=V_1\oplus\cdots\oplus V_r$ be a block decomposition into coordinate subspaces of dimensions $d_1,\dots,d_r$. Let $a_m$ be the diagonal matrix with $a_m|_{V_i}=z_{i,m}\,b_{i,m}$, where $z_{i,m}>0$ are scalars with $z_{i+1,m}/z_{i,m}\to0$ and $b_{i,m}$ are positive diagonal matrices converging to invertible positive diagonal matrices $b_i$. For a Cartan sequence one can take $z_{i,m}$ to be the largest singular value of its $i$th block. Bounded ratios within each block then keep the normalised entries uniformly away from zero, and after extraction the matrices $b_{i,m}$ converge. Thus this form covers every extracted Cartan sequence of type $\Theta=\{\alpha_{d_1},\alpha_{d_1+d_2},\dots\}$, with the simple roots at the block boundaries.

\begin{proposition}\label{prop:proj}
Let $x=[y]\in\mathbb K\PP^n$ with $y=(y_1,\dots,y_r)$, $y_i\in V_i$, and let $i_0=\min\{i:y_i\ne0\}$. Then
\[
\D_{(a_m)}(x)=\PP\bigl(V_1\oplus\cdots\oplus V_{i_0-1}\oplus\mathbb K\,b_{i_0}y_{i_0}\bigr).
\]
In particular the fibre is a projective subspace of $\mathbb K$-dimension $d_1+\cdots+d_{i_0-1}$; it is a point exactly when $i_0=1$.
\end{proposition}
\begin{proof}
For representatives $y^{(m)}\to y$, write
\[
 z_{i_0,m}^{-1}a_my^{(m)}=w_m+u_m,
 \qquad w_m\in V_{<i_0},\quad u_m\longrightarrow b_{i_0}y_{i_0}\ne0.
\]
The convergence follows from $z_{j,m}/z_{i_0,m}\to0$ for $j>i_0$
and boundedness of the Levi blocks. Along a convergent image subsequence,
extract further so that $(w_m)$ is bounded or $\|w_m\|\to\infty$.
In the first case the limit has the form $[w+b_{i_0}y_{i_0}]$;
in the second it lies in $\PP(V_{<i_0})$.
Both lie in the stated subspace.

Conversely, let $w\in V_{<i_0}$. Put $y^{(m)}_j=(z_{i_0,m}/z_{j,m})\,b_{j,m}^{-1}w_j\to0$ for $j<i_0$ and $y^{(m)}_j=y_j$ for $j\ge i_0$. Then $a_my^{(m)}/z_{i_0,m}\to w+b_{i_0}y_{i_0}$ along the whole sequence. Points at infinity of this affine chart are obtained by letting $|w|\to\infty$ slowly (diagonal choice), exactly as in Theorem~\ref{thm:dest}, Step~1.
\end{proof}
\begin{remark}
Here $\rho(x)=[0:\dots:y_{i_0}:\dots:0]$, and the block-upper-unipotent orbit is the affine space $V_{<i_0}+b_{i_0}y_{i_0}$; its projective closure adds $\PP(V_{<i_0})$. Thus the source and fibre dimensions sum to $n$ over $\mathbb K$, or to $(\dim_{\R}\mathbb K)n=\dim_{\R}X$ as in Corollary~\ref{cor:conservation}. The blocks $(1,2)$ in $\PSL(3,\R)$ give Proposition~\ref{prop:sl3}.
\end{remark}

%--------------------------------------------------------------
\section{Boundary compactification and stratified convergence}\label{sec:boundary}

\subsection{The exterior boundary of a general group}\label{sec:extG}
A regular highest-weight representation provides a sufficient
graded Cartan datum for describing destinations. It need not be minimal
for a given flag manifold: compact normal factors may be visible in
the projective boundary even when they act trivially on that flag
manifold. We use a projectively faithful representation of $G$; its
highest restricted weight detects every simple root, while its
boundary operators retain the relative Levi class.

\paragraph{The representation.} Let $\varrho:\tilde G\to\mathrm{GL}(V)$ be a finite-dimensional complex irreducible representation of the universal cover whose highest \emph{restricted} weight $\lambda$ is regular, $\langle\lambda,\alpha^\vee\rangle>0$ for all $\alpha\in\Delta$. Choose a complex Cartan subalgebra containing $\ak_\C$ and a positive complex root system adapted to the restricted positive roots: a positive complex root restricts either to a positive restricted root or to zero. Let $\rho_{\C}$ be the sum of the complex fundamental weights, equivalently half the sum of all positive complex roots, and let $\varrho$ have complex highest weight $\rho_{\C}$. Restriction to $\ak$ gives
\[
\lambda=\rho_{\C}|_{\ak}=\frac12\sum_{\beta\in\Phi^+}m_\beta\,\beta,
\qquad m_\beta=\dim_{\R}\g_\beta.
\]
For a simple restricted root $\alpha$, its simple reflection permutes the positive restricted roots other than $\alpha$ and, if present, $2\alpha$, and reverses those exceptional roots. Since multiplicities are Weyl-invariant, the reflection formula yields
\[
\langle\lambda,\alpha^\vee\rangle=m_\alpha+2m_{2\alpha}>0
\qquad(m_{2\alpha}=0\text{ if }2\alpha\notin\Phi).
\]
This proves restricted regularity, including non-reduced root systems. The chosen complex highest weight is nonzero on every complex simple ideal, including the complexifications of compact real factors.

The centre of $\tilde G$ acts by scalars, so $\varrho$ induces $\bar\varrho:G\to\PGL(V)$. We require it to act non-trivially on every simple factor, including the compact ones. Its differential then has zero projective kernel: on a simple ideal a scalar action would be zero, while a nonzero representation is faithful at the Lie-algebra level. The projective kernel in the connected centre-free group $G$ is consequently discrete and central, hence trivial. The construction using all complex fundamental weights has this property, so no factor is hidden.

Choose a Hermitian form compatible with a Cartan involution $\varsigma$, so $d\varrho(\varsigma X)=-d\varrho(X)^*$; it is obtained by making the compact form $\kk\oplus i\pk$ unitary. Thus $K$ acts unitarily and the restricted weight spaces are orthogonal.

The projective image is closed in $\PGL(V)$. Indeed, $s(H)=\max_\nu\nu(H)-\min_\nu\nu(H)$ is positive for every $0\ne H\in\ak$: otherwise $d\varrho(H)$ would be scalar, contrary to the zero projective differential kernel. By homogeneity and compactness of the unit sphere, $s(H)\ge c\|H\|$ for some $c>0$. If $\|\mu_m\|\to\infty$, the condition number $e^{s(\mu_m)}$ tends to infinity, so the projective operators cannot converge to an invertible one. Hence $\bar\varrho$ is an embedding, and its closure in $\PP(\mathrm{End}\,V)$ compactifies $G$. In all formulas $\varrho(g)$ for $g\in G$ means any lift to $\tilde G$; the resulting projective operator is independent of that choice.

Write $V=\bigoplus V_\nu$ for the weight decomposition under $A$ (restricted weights); every weight is $\nu=\lambda-\sum_{\alpha\in\Delta}n_\alpha\alpha$ with $n_\alpha\in\mathbb Z_{\ge0}$. For $\Theta\subseteq\Delta$ put
\[
V^\Theta=\bigoplus\{V_\nu:\ n_\alpha(\nu)=0\ \ \forall\alpha\in\Theta\},\qquad V=V^\Theta\oplus V^{\Theta,\perp},
\]
$V^{\Theta,\perp}$ the sum of the other weight spaces, and $\Pi_\Theta$ the projection onto $V^\Theta$ along $V^{\Theta,\perp}$. $V^\Theta$ is the eigenspace of $\varrho(A_\Theta)$ of largest weight, hence $L_\Theta$-stable, and $\Pi_\Theta$ commutes with $\varrho(L_\Theta)$.

\begin{lemma}\label{lem:Pi}\
(i) Use the canonical exponential lifts of $N_\Theta$ and $N_\Theta^-$. Then $\varrho(n)\Pi_\Theta=\Pi_\Theta$ for $n\in N_\Theta$ (a raising operator applied to $V^\Theta$ leaves the set of weights) and $\Pi_\Theta\varrho(n^-)=\Pi_\Theta$ for $n^-\in N_\Theta^-$ (lowering operators have no component in $V^\Theta$). (ii) The stabiliser of the subspace $V^\Theta$ in $G$ is $P_\Theta$, and that of $V^{\Theta,\perp}$ is $P_\Theta^-$.
\end{lemma}
\begin{proof}
(i) If $\gamma$ is a root of $N_\Theta$, it has a positive coefficient at some $\alpha\in\Theta$. For a weight in $V^\Theta$, raising by $\gamma$ would make that coefficient of $\lambda-\nu$ negative, so the raising operator is zero. Lowering cannot send any weight outside $V^\Theta$ into $V^\Theta$. Exponentiating gives the two identities.

(ii) The stabiliser $H$ of $V^\Theta$ contains $P_\Theta$, since $L_\Theta$ preserves it and $N_\Theta$ fixes it pointwise. In particular $H$ contains the minimal real parabolic $B$. Closed overgroups of $B$ are the standard real parabolics, obtained by adjoining the negative simple restricted root spaces \cite[Propositions~7.76, 7.82(d) and~7.83(b)]{KnappLG}. To prove $H=P_\Theta$ it is therefore enough to exclude $\g_{-\alpha}$ for each simple $\alpha\in\Theta$; no assertion about an arbitrary negative root is needed.

To exclude these negative simple root spaces also for nonsplit real groups, choose $0\ne X\in\g_\alpha$ and put $Y=-\varsigma X\in\g_{-\alpha}$. Then $[X,Y]=cH_\alpha$ with $c>0$ and $H_\alpha$ the positive restricted coroot direction. Indeed $[X,Y]$ centralises $\ak$ and is sent to its negative by $\varsigma$; because $\g_0=\mathfrak m\oplus\ak$ with $\mathfrak m\subset\kk$, the bracket lies in $\ak$. Its sign follows by pairing with $H\in\ak$ and using the positive form $-B(\cdot,\varsigma\cdot)$. For $0\ne v\in V_\lambda$, $d\varrho(X)v=0$ and $d\varrho(Y)=d\varrho(X)^*$. Hence
\[
 \|d\varrho(Y)v\|^2
 =\langle d\varrho([X,Y])v,v\rangle
 =c\lambda(H_\alpha)\|v\|^2>0.
\]
The resulting vector has restricted weight $\lambda-\alpha$, outside $V^\Theta$. Thus $\g_{-\alpha}$ does not preserve $V^\Theta$. This uses neither one-dimensional restricted root spaces nor a reduced restricted root system.

The classification now gives $\mathrm{Lie}(H)=\mathfrak p_\Theta$. Every element of $H$ normalises its Lie algebra, and $N_G(\mathfrak p_\Theta)=P_\Theta$ \cite[Proposition~7.83(b)]{KnappLG}. Hence $P_\Theta\subseteq H\subseteq N_G(\mathfrak p_\Theta)=P_\Theta$, which proves the equality for the full stabiliser, including all components.

Finally $V^{\Theta,\perp}$ is the orthogonal complement for the chosen form, and $\varrho(g)^*=\varrho(\varsigma(g)^{-1})$ up to scalar, so its stabiliser is $\varsigma(H)=P_\Theta^-$.
\end{proof}

\begin{lemma}\label{lem:levi-scalar}
For the representation in this section, let
\[
Z_\Theta=\{z\in L_\Theta:\ \varrho(z)|_{V^\Theta}
\text{ is scalar}\}.
\]
Every element of $Z_\Theta$ fixes every point on every flag manifold
of the reductive Levi group $L_\Theta$, with its split centre and
compact normal factors acting trivially on those flag quotients.
Consequently it fixes every value of the stratified retraction
$\rho:X\to X^{A_\Theta}$, and the destination formula is unchanged
when its relative Levi representative $b$ is replaced by $bz$ or $zb$.
\end{lemma}
\begin{proof}
Let $J\subseteq\Delta\setminus\Theta$. The nested top restricted
weight spaces satisfy $V^{\Theta\cup J}\subseteq V^\Theta$.
By Lemma~\ref{lem:Pi}(ii), applied to $\Theta\cup J$, the stabiliser
of $V^{\Theta\cup J}$ in $G$ is $P_{\Theta\cup J}$.
Thus the Levi flag variety
\[
L_\Theta/(L_\Theta\cap P_{\Theta\cup J})\longrightarrow
\mathrm{Gr}_{\dim V^{\Theta\cup J}}(V^\Theta),\qquad
l\longmapsto\varrho(l)V^{\Theta\cup J}
\]
is equivariantly embedded: equality of stabilisers gives injectivity,
and compactness of the flag quotient makes the map an embedding.
Since $z\in Z_\Theta$ acts by a scalar on $V^\Theta$, it fixes every
point of the displayed Grassmannian and hence every point of each
embedded flag quotient. Every parabolic of the reductive Levi is
conjugate to a standard one of this form, and its split centre and
compact normal factors are contained in all these parabolics.
The components of $X^{A_\Theta}$ are precisely these Levi flag
quotients (Lemma~\ref{lem:cells}), so $z$ fixes $\rho(x)$ for every
$x$. Finally $Z_\Theta\subset L_\Theta$ normalises $N_\Theta$;
normalising the orbit inside its closure gives
$\overline{N_\Theta bz\rho(x)}=
 \overline{N_\Theta b\rho(x)}$, and similarly for $zb$.
\end{proof}
\begin{theorem}\label{thm:extG}\
Let $g_m=k_me^{\mu_m}k'_m$ have Cartan data $(\Theta,k,k',b)$ in the sense of Lemma~\ref{lem:split}. Then in $\PP(\mathrm{End}\,V)$
\[
\lim_m[\varrho(g_m)]=[\varrho(k)\,\Pi_\Theta\varrho(b)\,\varrho(k')].
\]
Let $X_\varrho$ be the closure of $\bar\varrho(G)$ in $\PP(\mathrm{End}\,V)$ and $\partial X_\varrho=X_\varrho\setminus\bar\varrho(G)$.
\begin{enumerate}[label=\textup{(\alph*)},nosep]
\item $\partial X_\varrho=\bigsqcup_{\varnothing\ne\Theta\subseteq\Delta}\mathcal S_\Theta$, with $\mathcal S_\Theta=\{[\varrho(k)\Pi_\Theta\varrho(\ell)\varrho(k')]:k,k'\in K,\ \ell\in L_\Theta\}$, and $\mathcal S_\Theta$ is the $G\times G$-orbit of $[\Pi_\Theta]$ for the action $(g,h)\cdot\Psi=\varrho(g)\Psi\varrho(h)^{-1}$.
\item The map $\mathcal S_\Theta\to G/P_\Theta\times G/P^-_\Theta$, $[\Psi]\mapsto(\text{stabiliser of }\mathrm{Im}\Psi,\ \text{stabiliser of }\ker\Psi)$, is a $G\times G$-equivariant fibre bundle with fibre $L_\Theta/Z_\Theta$, where $Z_\Theta$ is the subgroup of $L_\Theta$ acting by scalars on $V^\Theta$ (it contains $A_\Theta$). For the limit of $(g_m)$, the stabiliser of the image is the attracting parabolic $kP_\Theta k^{-1}$ and the stabiliser of the kernel is the repelling parabolic $k'^{-1}P_\Theta^-k'$ of Section~\ref{sec:groups}.
\item \textup{(Realisation)} Every point of $\mathcal S_\Theta$ is the limit of $g_m=k\exp(m\nu_\Theta)\,\ell\,k'$, where $\nu_\Theta\in\ak_\Theta$ has $\alpha(\nu_\Theta)=1$ for $\alpha\in\Theta$ and $0$ for $\alpha\notin\Theta$; this sequence has type exactly $\Theta$. Hence $\partial X_\varrho$ is exactly the set of limits of divergent sequences, and \emph{every datum $(\Theta,\xi^+,\xi^-,[\ell])$ occurs}.
\item The destination relation $\D_{(g_m)}$ depends on $(g_m)$ only through the limit point in $\partial X_\varrho$.
\end{enumerate}
\end{theorem}
\begin{proof}
\emph{The limit.} On $V_\nu$, $\varrho(e^{\mu_m})$ is the scalar $e^{\nu(\mu_m)}=e^{\lambda(\mu_m)}e^{-\sum n_\alpha\alpha(\mu_m)}$. If $n_\alpha>0$ for some $\alpha\in\Theta$ this is $o(e^{\lambda(\mu_m)})$, because $\alpha(\mu_m)\to\infty$ and all $\alpha'(\mu_m)\ge0$. If $n_\alpha=0$ for $\alpha\in\Theta$, then only $\alpha\notin\Theta$ occur, $\alpha(\mu_m)=\alpha(\eta_m)\to\alpha(\eta)$, and the scalar divided by $e^{\lambda(\mu_m)}$ tends to $e^{-\sum n_\alpha\alpha(\eta)}$, which is $\varrho(b)|_{V_\nu}$ divided by $e^{\lambda(\eta)}$ ($b=\exp\eta$). Since $k_m\to k$, $k'_m\to k'$, the claim follows.

\emph{(a)} By compactness every sequence in $G$ leaving compacts has a subsequence with a Cartan datum, so limits lie in the right-hand side; conversely (c) shows that every point is attained.

For the orbit statement, Lemma~\ref{lem:Pi}(i) shows that if $g=kp$ with $p=ln\in P_\Theta$ and $h^{-1}=p'^-k'$ with $p'^-=n^-l'\in P_\Theta^-$ (Iwasawa-type decompositions $G=KP_\Theta=P_\Theta^-K$), then $[\varrho(g)\Pi_\Theta\varrho(h)^{-1}]=[\varrho(k)\varrho(l)\Pi_\Theta\varrho(l')\varrho(k')]$.

\emph{(b)} By Lemma~\ref{lem:Pi}(ii), $\mathrm{Im}=\varrho(k)V^\Theta$ has stabiliser $kP_\Theta k^{-1}$, and $\ker=\varrho(k')^{-1}V^{\Theta,\perp}$ has stabiliser $k'^{-1}P^-_\Theta k'$. Standard parabolics of distinct types are not conjugate, so the image stabiliser recovers $\Theta$ and the strata in (a) are disjoint. The $G\times G$-stabiliser of $[\Pi_\Theta]$ is $\{(ln,\,l'n^-)\}$ with $\varrho(l)\Pi_\Theta=c\,\Pi_\Theta\varrho(l')$, i.e.\ $l'\equiv l$ modulo $Z_\Theta$; this gives the fibre $L_\Theta/Z_\Theta$. Here $Z_\Theta$ is a closed normal subgroup, namely the projective kernel of the action on $V^\Theta$. The displayed fibre is identified with this quotient after choosing frames for the two parabolics; the identification is not canonical over the base.

To verify the topology, suppose $[\varrho(\ell_j)|_{V^\Theta}]$ converges in $\PGL(V^\Theta)$. Using the Levi polar decomposition described in (c), write
\[
 \ell_j=u_j\exp(\nu_j+\eta_j)v_j,\qquad
 u_j,v_j\in K\cap L_\Theta,\quad \nu_j\in\ak_\Theta,\quad
 \eta_j\in\operatorname{span}\{\alpha^\vee:\alpha\notin\Theta\}.
\]
The factor $\exp\nu_j$ lies in $Z_\Theta$. For every $\alpha\notin\Theta$, the root vector calculation in Lemma~\ref{lem:Pi} gives the weights $\lambda$ and $\lambda-\alpha$ in $V^\Theta$. Thus the difference $s_\Theta(\eta_j)$ between the largest and smallest values on $\eta_j$ of weights occurring in $V^\Theta$ satisfies $s_\Theta(\eta_j)\ge|\alpha(\eta_j)|$. The condition number of $\varrho(\ell_j)|_{V^\Theta}$ is $e^{s_\Theta(\eta_j)}$, so convergence to an invertible projective operator bounds all these root values. The Levi Cartan matrix is invertible, hence $\eta_j$ stays bounded; compactness of $K\cap L_\Theta$ then makes $\ell_jZ_\Theta$ relatively compact. This proves that $L_\Theta/Z_\Theta\to\PGL(V^\Theta)$ is a closed embedding.

Image and kernel vary continuously on this stratum of fixed rank, and their flag orbits are compact. Local sections of $G\to G/P_\Theta$ and $G\to G/P_\Theta^-$ put them in the standard frames, where the operator is determined by its invertible restriction to $V^\Theta$. The closed embedding just proved gives product charts with the subspace topology of $\PP(\mathrm{End}\,V)$, proving the bundle assertion. No connectedness assumption on $L_\Theta$ is used.

\emph{(c)} $\varrho(\exp(m\nu_\Theta))$ acts on $V^\Theta$ (top eigenspace of $\varrho(A_\Theta)$) by the scalar $e^{m\lambda(\nu_\Theta)}$ and on the other eigenspaces by strictly smaller powers; $\ell\in L_\Theta$ commutes with $A_\Theta$ and preserves the eigenspaces. Hence $e^{-m\lambda(\nu_\Theta)}\varrho(\exp(m\nu_\Theta)\ell)\to\Pi_\Theta\varrho(\ell)$.

For the type, write $\ell=u\exp(\eta)v$ with $u,v\in K\cap L_\Theta$ and $\eta$ dominant for the Levi roots. This polar decomposition stays in the centraliser of $A_\Theta$: the positive polar factor and its logarithm commute with $A_\Theta$, and the compact factor does as well. The logarithm can be conjugated into $\ak$ by $K\cap L_\Theta$ and then into the closed Levi chamber. Since $u,v$ commute with $\exp(m\nu_\Theta)$ and $m\nu_\Theta+\eta$ belongs to $\ak^+$ for all sufficiently large $m$, this is exactly the Cartan projection of $k\exp(m\nu_\Theta)\ell k'$. Thus the complementary root values are constant and the root values on $\Theta$ are $m+\alpha(\eta)$.

\emph{(d)} Replace $(k,k',b)$ by $(kq_1,q_2k',\,q_1^{-1}bq_2^{-1})$ with $q_1,q_2\in K\cap L_\Theta$; the limit point is unchanged since $\Pi_\Theta$ commutes with $\varrho(L_\Theta)$. By Theorem~\ref{thm:dest} and the $L_\Theta$-equivariance of $\rho$ and the normalisation of $N_\Theta$ by $L_\Theta$, $k q_1\,\overline{N_\Theta\,q_1^{-1}b q_2^{-1}\rho(q_2k'x)}=k\,\overline{N_\Theta\,b\rho(k'x)}$.

Two sets of data with the same limit point have the same $\Theta$, parabolics and (by (b)) Levi classes modulo $Z_\Theta$, hence differ by these frame changes and by an element of $Z_\Theta$. Lemma~\ref{lem:levi-scalar} shows that this remaining ambiguity fixes every point of each Levi flag component and leaves $\overline{N_\Theta\,b\rho(y)}$ unchanged. This proves (d).
\end{proof}
\begin{remark}
The normalised limits are obtained from restricted-weight estimates,
so no reducedness or root-multiplicity condition is imposed.
The stabilisers of their image and kernel are identified by the
negative simple root spaces (Lemma~\ref{lem:Pi}), and the scalar
kernel of the Levi action is accounted for by
Lemma~\ref{lem:levi-scalar}. For a real form of $G$, the boundary
considered here consists of limits of sequences in $G$ itself.
It need not coincide with the full real locus of the complex
wonderful compactification.
\end{remark}

%--------------------------------------------------------------
\subsection{The stratified convergence principle}\label{sec:strat}
Over each source stratum, the boundary datum places a family of
translated Schubert varieties, parametrised by a Levi flag manifold.
Inversion gives the corresponding target description. Proposition~\ref{prop:proj}
is the projective model.

\paragraph{Notation.} For $\Theta\subseteq\Delta$ non-empty and $\sigma\in\mathcal C_\Theta=W_\Theta\backslash W/W_\theta$ let $O^-_\sigma=P^-_\Theta\sigma Q$ and $O^+_\sigma=P_\Theta\sigma Q$ (source and target cells), and with $\Sigma_+,\Sigma_0,\Sigma_-$ as in Lemma~\ref{lem:cells}(iv) put
\[
d_\pm(\sigma)=\sum_{\lambda\in\Sigma_\pm(\sigma)}\dim\g_\lambda,\qquad d_0(\sigma)=\sum_{\lambda\in\Sigma_0(\sigma)}\dim\g_\lambda,\qquad d_++d_0+d_-=\dim X.
\]
Then $\dim O^+_\sigma=d_++d_0$, $\dim O^-_\sigma=d_0+d_-$, and $d_0(\sigma)$ is the dimension of the flag manifold $L_\Theta\sigma Q$ of $L_\Theta$. If $\sigma$ is the representative of minimal length in its double coset, put $X_\sigma=\overline{B\sigma Q}$, the \emph{Schubert variety} of $\sigma$. Finally $\rho^+:X\to X^{A_\Theta}$ is the stratified retraction of Lemma~\ref{lem:cells}(iii) for the opposite direction, $\rho^+(x)=\lim_{z\to\infty}z^{-1}x$, with smooth restrictions to the cells $O^+_\tau$.

\begin{lemma}\label{lem:schubert-cell}
Let $\sigma$ be the minimal-length representative of a double coset in $W_\Theta\backslash W/W_\theta$. Then
\[
N_\Theta\sigma Q/Q=B\sigma Q/Q,
\qquad
\overline{N_\Theta\sigma Q/Q}=X_\sigma.
\]
In particular, the closure is a Schubert variety of dimension $d_+(\sigma)$.
\end{lemma}
\begin{proof}
Write $N=N_\Theta N_{L_\Theta}$ for the unipotent radical of the chosen minimal parabolic. Minimality of $\sigma$ implies $N_{L_\Theta}\sigma Q=\sigma Q$, and hence $N\sigma Q=N_\Theta\sigma Q$. Since $B=MAN$, where $M=Z_K(A)$, and $\sigma^{-1}M\sigma=M\subset Q$, both $M$ and $A$ fix $\sigma Q$. As $MA$ normalises $N$, this gives $B\sigma Q=N\sigma Q=N_\Theta\sigma Q$. Taking closures gives the Schubert variety $X_\sigma$. Its dimension is the dimension of the expanding cell, namely $d_+(\sigma)$.
\end{proof}

\Needspace{8\baselineskip}
\begin{theorem}\label{thm:S}
Let $G$ be as in Section~\ref{sec:setting-sum}, $X=G/Q$ a flag manifold, $\varrho$ the representation of Theorem~\ref{thm:extG}, and $(g_m)$ a sequence leaving every compact set of $G$.
\begin{enumerate}[label=\textup{(\roman*)},nosep]
\item \textup{(The asymptotic datum.)} After extraction $[\varrho(g_m)]$ converges to a point $\xi$ of the exterior boundary $\partial X_\varrho$, which lies in a stratum $\mathcal S_\Theta$, $\Theta\ne\varnothing$. The point $\xi$ is the same thing as the \emph{graded datum} $\mathcal L(\xi)=(\Theta,\ \xi^+,\ \xi^-,\ [b])$: the type $\Theta$, the attracting parabolic $\xi^+=kP_\Theta k^{-1}$, the repelling parabolic $\xi^-=k'^{-1}P^-_\Theta k'$, and a point $[b]$ of the fibre $L_\Theta/Z_\Theta$ over $(\xi^+,\xi^-)$ \textup{(}the Levi class\textup{)}.
\item \textup{(Strata.)} The sets
\[
\mathcal O^-_\sigma(\xi)=k'^{-1}O^-_\sigma\quad(\text{sources}),\qquad\mathcal O^+_\tau(\xi)=k\,O^+_\tau\quad(\text{targets}),\qquad\sigma,\tau\in\mathcal C_\Theta,
\]
decompose $X$ into finitely many strata, which depend only on $\xi^-$ and $\xi^+$ respectively.
\item \textup{(Destination formula, source side.)} Let $x\in\mathcal O^-_\sigma(\xi)$, with $\sigma$ of minimal length, and $\rho(k'x)=l\,\sigma Q$, $l\in L_\Theta$. Then
\[
\D_{(g_m)}(x)=(k\,b\,l)\cdot X_\sigma,
\]
a translate of one Schubert variety of dimension $d_+(\sigma)$, which depends on $x$ only through the point $\rho(k'x)$ of the flag manifold $L_\Theta\sigma Q$. It lies in $\bigcup_{\tau\le\sigma}\mathcal O^+_\tau(\xi)$, and its intersection with $\mathcal O^+_\sigma(\xi)$ is open and dense in the destination fibre. The fibre is a single point iff $\Sigma_+(\sigma)=\varnothing$.
\item \textup{(Inverse formula, target side.)} Let $y\in\mathcal O^+_\tau(\xi)$ and $\rho^+(k^{-1}y)=l'\tau Q$. Then
\[
\D_{(g_m^{-1})}(y)=k'^{-1}\,\overline{N_\Theta^-\,b^{-1}l'\,\tau Q}=(k'^{-1}b^{-1}l')\cdot\overline{N^-_\Theta\tau Q},
\]
the closure of a cell of dimension $d_-(\tau)$; and $\D_{(g_m^{-1})}$ is the transpose of $\D_{(g_m)}$.
\item \textup{(Bundles.)} The map $x\mapsto\rho(k'x)$ makes $\mathcal O^-_\sigma(\xi)$ a fibre bundle over $L_\Theta\sigma Q$ whose fibres are cells of dimension $d_-(\sigma)$, and $x\mapsto\D_{(g_m)}(x)$ is constant on its fibres. The translates of $X_\sigma$ form a locally trivial family over $L_\Theta\sigma Q$, and the destination incidence over $\mathcal O^-_\sigma(\xi)$ is its pullback. For complex $G$ the bundles and local trivialisations are holomorphic, with complex dimensions equal to half the displayed real dimensions.
\item \textup{(Uniqueness and realisation.)} The relation $\D_{(g_m)}$ depends on $(g_m)$ only through $\xi$: every sequence with the same limit point has the same relation. Every $\Theta\ne\varnothing$ and every datum $(\Theta,\xi^+,\xi^-,[b])$ occurs, for a sequence of type exactly $\Theta$.
\end{enumerate}
\end{theorem}
\begin{proof}
(i) and the realisation in (vi) follow from Theorem~\ref{thm:extG}(b),(c). For uniqueness in (vi), consider any sequence whose projective operators converge to $\xi$. Every subsequence has a further Cartan extraction with the same exterior limit. By Theorem~\ref{thm:extG}(d) all these extracted relations agree, and Theorem~\ref{thm:dest}(d) gives Hausdorff convergence of their graphs to this common relation. The graphs of the full sequence therefore converge to the same relation. Each pair in it is realised along the full sequence by choosing points in the graphs approaching the pair, and every full-sequence limit belongs to it.

\noindent(ii) follows from Lemma~\ref{lem:cells}(i) and its mirror image for $P_\Theta$, translated by $k'^{-1}$ and $k$.

\noindent(iii) By Theorem~\ref{thm:dest}(b), $\D_{(g_m)}(x)=k\overline{N_\Theta\,b\rho(k'x)}=k\overline{N_\Theta\,bl\sigma Q}$. As $bl\in L_\Theta$ normalises $N_\Theta$, this is $kbl\,\overline{N_\Theta\sigma Q}$. For $\sigma$ minimal, the proof of Lemma~\ref{lem:schubert-cell} shows $N_\Theta\sigma Q=N\sigma Q=B\sigma Q$, so $\overline{N_\Theta\sigma Q}=X_\sigma$, of dimension $d_+(\sigma)$ by Theorem~\ref{thm:dest}(c). The cell $B\sigma Q$ is open and dense in its closure. It lies in $O^+_\sigma$, while the rest of $X_\sigma$ is a union of $B\tau'Q$ with $\tau'<\sigma$, which lie in the strata $O^+_{[\tau']}$, $[\tau']<\sigma$ by the minimality of $\sigma$ in its double coset and Lemma~\ref{lem:schubert-cell}. The elements $b,l\in L_\Theta\subset P_\Theta$ preserve every $O^+_\tau$. The point criterion follows from the cell dimension $d_+(\sigma)$.

\noindent(iv) Write $g_m^{-1}=k_m'^{-1}a_m^{-1}k_m^{-1}$. Now $a_m^{-1}\in\exp(-\ak^+)$ is a Cartan sequence for the opposite positive system: the roles of $N_\Theta,N_\Theta^-$ and of $\rho,\rho^+$ are exchanged, the Levi class is $b^{-1}$, and the proof of Theorem~\ref{thm:dest}(a),(b) is symmetric. This gives the formula; the last sentence is Theorem~\ref{thm:dest}(d), and the dimension is as in (iii) with $N_\Theta^-$ in place of $N_\Theta$.

\noindent(v) Put $S_0=L_\Theta\cap\sigma Q\sigma^{-1}$ and $U_0=N_\Theta^-\cap\sigma Q\sigma^{-1}$. The stabiliser of $\sigma Q$ in $P_\Theta^-$ is $U_0S_0$. Indeed, if $n l$ fixes $\sigma Q$ with $n\in N_\Theta^-$ and $l\in L_\Theta$, conjugation by $\exp(t\nu_\Theta)$ contracts $n$ and fixes $l$. It preserves the stabiliser because $\sigma Q$ is fixed by $A_\Theta$. Taking $t\to+\infty$ shows $l\in S_0$ and then $n\in U_0$.

Consequently the source stratum, in the source frame, is the associated bundle $L_\Theta\times^{S_0}(N_\Theta^-/U_0)$ over $L_\Theta/S_0$. Its fibre is a cell of dimension $d_-(\sigma)$ by Lemma~\ref{lem:cells}(iv) for the opposite system.

Since $S_0$ fixes $\sigma Q$ and normalises $N_\Theta$, it preserves $X_\sigma=\overline{N_\Theta\sigma Q}$. The bundle $L_\Theta\times^{S_0}X_\sigma$ therefore parametrises the destination fibres by $[l,y]\mapsto (l\sigma Q,kbl\,y)$, by (iii). Pulling it back along $x\mapsto\rho(k'x)$ gives the claimed incidence.

Local smooth sections of $L_\Theta\to L_\Theta/S_0$, holomorphic for complex $G$, give the local trivialisations.
\end{proof}

\begin{remark}
In rank one, $\mathcal C_\Theta=\{1,w\}$, with $X_1$ a point
and $X_w=X$, recovering Corollary~\ref{cor:rank-one}. For
$\PSL(2,\C)$ the fibre $L_\Theta/Z_\Theta$ is a point, so the
boundary datum is just $(x^+,x^-)$.
\end{remark}

\subsection{Comparison with complex graph limits}

\begin{proposition}\label{prop:hilbert-support}
Suppose that $G$ is a connected complex semisimple group of adjoint
type acting faithfully on $X=G/Q$.
For a sequence with extracted Cartan data $(\Theta,k,k',b)$,
the graphs $\mathcal G_{g_m}=\{(x,g_mx):x\in X\}$ converge in
Brion's Hilbert compactification. The limiting subscheme is reduced
and has support
\[
 \{(x,y)\in X\times X:
       y\in k\overline{N_\Theta b\rho(k'x)}\}.
\]
Every pair in this support is approached along the whole sequence
of graphs.
\end{proposition}
\begin{proof}
Write $g_m=k_m z_m b_m k'_m$ as in Lemma~\ref{lem:split},
where $z_m=\exp\nu_m$ and $b_m\to b$, and put
$J=\Delta\setminus\Theta$. We use Brion's identification of the
Hilbert orbit closure with the wonderful compactification and his
description of its reduced fibres
\cite[Theorem~3 and Proposition~10]{r2grBrion2003}.
The group parameter $h$ in that compactification corresponds to
$\mathcal G_{h^{-1}}$: indeed, $(a,c)$ sends the diagonal to
$\mathcal G_{ca^{-1}}$, whereas its group parameter is $ac^{-1}$.
Thus $\mathcal G_{z_m}$ has torus-chart coordinates
$q_\alpha=e^{-\alpha(\nu_m)}$, tending to $1$ for $\alpha\in J$
and to $0$ otherwise. Its Hilbert limit is the canonical fibre of
type $J$. In the notation of the present article, Brion's component
description gives its support as
\[
 C_\Theta=\bigcup_{\sigma\in\mathcal C_\Theta}
 \overline{(N^-_\Theta\times N_\Theta)
                  \operatorname{diag}(L_\Theta)
                  \cdot(\sigma Q,\sigma Q)}.
\]
Here the representatives are of minimal length. Equivalently,
each component is
$\operatorname{diag}(L_\Theta)
 (\overline{N^-_\Theta\sigma Q}\times X_\sigma)$:
the first factor in Brion's notation is
$w_J\overline{B^-w_J\sigma Q}$, where $w_J$ is longest in
$W_\Theta$, and
$w_JB^-w_J^{-1}=N^-_\Theta(B\cap L_\Theta)$;
the group $B\cap L_\Theta$ fixes $\sigma Q$.

We identify this support directly using the destination formula.
For the central sequence $(z_m)$, Theorem~\ref{thm:dest} gives
the relation
\[
 R_\Theta=\bigcup_{\sigma\in\mathcal C_\Theta}
 \operatorname{diag}(L_\Theta)
       (N^-_\Theta\sigma Q\times X_\sigma).
\]
Indeed, a source $n^-l\sigma Q$ retracts to $l\sigma Q$ and
has destination $lX_\sigma$; conjugation by $l$ preserves
$N^-_\Theta$. Since $N_\Theta\sigma Q$ is dense in $X_\sigma$,
each piece displayed here has the same closure as the corresponding
orbit defining $C_\Theta$. The relation $R_\Theta$ is closed by
Theorem~\ref{thm:dest}(d), and the union is finite. Hence
$R_\Theta=C_\Theta$.

Finally,
$\mathcal G_{g_m}=(k_m'^{-1},k_m b_m)\mathcal G_{z_m}$.
The two factors converge, so the Hilbert limit has support
$(k'^{-1},kb)C_\Theta$. This is the stated relation by the
destination formula, since $b$ normalises $N_\Theta$.
The whole-sequence realisation is already part of
Theorem~\ref{thm:dest}(b),(d).
\end{proof}

\section{Destination relations of discrete groups}\label{sec:groups}

Properness for a discrete group depends on the union of its sequence relations. We first treat this union for actions on compact spaces.

Let $\Gamma$ be a countable group acting by homeomorphisms on a compact metrizable space $X$. Define
\[
\D_\Gamma=\{(x,y):\exists\text{ distinct }\gamma_m\in\Gamma,
\ x_m\to x,\ \gamma_mx_m\to y\},
\qquad
\mathscr R_\Gamma=\{x:(x,x)\in\D_\Gamma\}.
\]
The set $\mathscr R_\Gamma$ is the \emph{recurrent set}, and $\mathscr W_\Gamma=X\setminus\mathscr R_\Gamma$. A \emph{domain of proper discontinuity} is a $\Gamma$-invariant open set $U$ such that $\{\gamma:\gamma K\cap K\neq\varnothing\}$ is finite for every compact $K\subset U$.
For $A\subseteq X$, $\D_\Gamma(A)$ denotes the union of the destination fibres over points of $A$. A largest domain contains every domain of proper discontinuity; a maximal domain admits no strictly larger domain of proper discontinuity.

\begin{lemma}\label{lem:Dbasic}
$\D_\Gamma$ is closed, symmetric, and invariant under $\Gamma\times\Gamma$. In particular $\mathscr R_\Gamma$ is closed and $\Gamma$-invariant, and it contains every point with infinite stabiliser.
\end{lemma}
\begin{proof}
If $(x_r,y_r)\in\D_\Gamma$ tend to $(x,y)$, with witnesses $\gamma_{r,m}$ and $x_{r,m}\to x_r$, choose for each $r$ an index $m(r)$ such that $\gamma_{r,m(r)}$ differs from the earlier choices and $x_{r,m(r)}$, $\gamma_{r,m(r)}x_{r,m(r)}$ are $1/r$-close to $x_r$, $y_r$. Symmetry follows from $y_m=\gamma_mx_m\to y$, $\gamma_m^{-1}y_m\to x$, and invariance from the witnesses $h\gamma_mg^{-1}$, $gx_m$.
\end{proof}

\begin{proposition}\label{prop:recur-equiv}
For $x\in X$ the following are equivalent: \textup{(i)} $x\notin\mathscr R_\Gamma$; \textup{(ii)} $x$ has a neighbourhood $N$ with $\{\gamma:\gamma N\cap N\neq\varnothing\}$ finite; \textup{(iii)} $x$ lies in a domain of proper discontinuity.
\end{proposition}
\begin{proof}
If (ii) fails, choose pairwise distinct $\gamma_r$ and $x_r\in B(x,1/r)$ with $\gamma_rx_r\in B(x,1/r)$; then $x\in\mathscr R_\Gamma$. Given (ii), $U=\Gamma N$ is a domain: a compact $K\subset U$ is covered by finitely many $g_iN$, and $\gamma K\cap K\neq\varnothing$ forces $g_j^{-1}\gamma g_i$ into a finite set. Finally a recurrent point of a domain $U$ would give, for a compact neighbourhood $K\subset U$, infinitely many $\gamma$ with $\gamma K\cap K\neq\varnothing$.
\end{proof}

\begin{theorem}\label{thm:relational}
\begin{enumerate}[label=\textup{(\roman*)},nosep]
\item A $\Gamma$-invariant open set $U$ is a domain of proper discontinuity if and only if $(U\times U)\cap\D_\Gamma=\varnothing$.
\item $\mathscr W_\Gamma$ is the union of all domains of proper discontinuity, and every domain is contained in a maximal one.
\item A largest domain exists if and only if no two non-recurrent points are related by $\D_\Gamma$; it is then $\mathscr W_\Gamma$.
\item If there is no largest domain, there are at least two incomparable maximal domains.
\end{enumerate}
\end{theorem}
\begin{proof}
(i) A pair of $\D_\Gamma$ in $U\times U$ makes the action on the union of two compact neighbourhoods non-proper. Conversely, if $K\subset U$ is compact and $\gamma_mx_m\in K$ for infinitely many distinct $\gamma_m$ and $x_m\in K$, a limit of $(x_m,\gamma_mx_m)$ is a pair of $\D_\Gamma$ in $U\times U$.

(ii) Proposition~\ref{prop:recur-equiv} gives the union assertion. A chain of domains has a domain as its union by (i): any two points of the union lie in one member of the chain. Zorn's lemma therefore extends every domain to a maximal one.

(iii) A largest domain contains every domain, hence equals $\mathscr W_\Gamma$, and by (i) $\mathscr W_\Gamma$ is a domain if and only if $(\mathscr W_\Gamma\times\mathscr W_\Gamma)\cap\D_\Gamma=\varnothing$.

(iv) If the maximal domain given by (ii) were unique, it would contain every domain.
\end{proof}

Call a closed invariant set $L\subset X$ a \emph{vertex cover} if every pair of $\D_\Gamma$ has an endpoint in $L$. By Theorem~\ref{thm:relational}(i), $U\mapsto X\setminus U$ is an inclusion-reversing bijection between domains of proper discontinuity and vertex covers.

\begin{corollary}\label{cor:covers}
$\mathscr R_\Gamma$ is the intersection of all vertex covers. It is itself a vertex cover if and only if a largest domain exists, and then it is the smallest vertex cover.
\end{corollary}
\begin{proof}
Every cover contains the diagonal points of $\D_\Gamma$, that is, $\mathscr R_\Gamma$. If $x\notin\mathscr R_\Gamma$, Proposition~\ref{prop:recur-equiv} puts $x$ in a domain, whose complement is a cover missing $x$. The rest is (iii).
\end{proof}

\begin{remark}
Frances~\cite[Definition~1 and Proposition~1]{Frances} uses moving sources and subsequential accumulation for sequences escaping compact subsets of $\mathrm{Homeo}(X)$, with $X$ a manifold. Here distinctness refers to the abstract group $\Gamma$. For a closed discrete subgroup of $\mathrm{Homeo}(X)$ the unions agree, since subsequences of escaping sequences still escape. Kapovich--Leeb--Porti~\cite[\S5.1, Definition~5.6]{KLP} identify the diagonal with the nonwandering set and characterise properness by the absence of related pairs; their Example~5.5 exhibits distinct maximal domains. See also \cite[Theorem~11]{KapovichProper}. Theorem~\ref{thm:dest} computes these pairs on $G/Q$, including the non-recurrent pairs obstructing a largest domain.
\end{remark}

\begin{proposition}\label{prop:sequence-union}
If $\Gamma<G$ is discrete and $X=G/Q$, then $\D_\Gamma$ is the union of the relations $\D_{(\gamma_m)}$ over all sequences of pairwise distinct elements of $\Gamma$, after extraction of Cartan data. Each relation in this union is given by the destination formula.
\end{proposition}
\begin{proof}
A sequence of pairwise distinct elements of a discrete subgroup leaves every compact subset of $G$. Given a pair in $\D_\Gamma$, extract Cartan data from a witnessing sequence; both limits in the definition persist under extraction. The converse follows because every point in a sequence destination is reached along that extracted sequence (Theorem~\ref{thm:dest}).
\end{proof}

\begin{proposition}\label{prop:D-rank1}
Let $G$ have real rank one, $X=G/P$ its boundary sphere, and $\Gamma<G$ an infinite discrete subgroup with no global fixed point. Let $\Lambda$ be the set of attracting points of sequences of distinct elements. Then
\[
\D_\Gamma=(X\times\Lambda)\cup(\Lambda\times X),\qquad\mathscr R_\Gamma=\Lambda,
\]
and $X\setminus\Lambda$ is the largest domain of proper discontinuity.
\end{proposition}
\begin{proof}
Corollary~\ref{cor:rank-one} gives the fibres: a point $x^+\in\Lambda$ off the repeller and all of $X$ over the repeller $x^-$, and $x^-\in\Lambda$ by inversion. This gives the inclusion $\subseteq$. Conversely, for $y\in\Lambda$ with repeller $r$ and any $x\in X$, right multiplication by some $h\in\Gamma$ moves the repeller to $h^{-1}r\neq x$, so $(x,y)\in\D_\Gamma$; symmetry gives $\Lambda\times X$. Theorem~\ref{thm:relational}(iii) concludes.
\end{proof}

In higher rank, the attracting and recurrent sets and the complement of equicontinuity need not coincide. The recurrent set need not be a vertex cover, so a smallest vertex cover may fail to exist. Sections~\ref{sec:chain} and~\ref{sec:critical} describe these distinctions.

\section{The flag of limit sets: from Quint to Myrberg}\label{sec:chain}

Here $G$ is simple with trivial centre, $X=G/Q$ is a non-trivial flag manifold, and $\Gamma<G$ is discrete. The standard limit sets, equal in rank one, form a flag in higher rank. Moving-source, recurrence and Kulkarni definitions also apply to countable actions on compact metrizable spaces; concentration additionally requires a probability measure, whereas Cartan levels use the flag structure.

\medskip
\emph{Kulkarni's set.} Let $L_0$ be the closure of the set of points with infinite stabiliser and $L_1$ the closure of the set of limits $\lim\gamma_mx$, with $x\notin L_0$ and $(\gamma_m)$ pairwise distinct. Put $B_\Gamma=L_0\cup L_1$, the \emph{fixed-source set}. Let $L_2$ be the closure of the set of cluster points of the sets $\gamma_mK$, for compact $K\subset X\setminus B_\Gamma$ and distinct $\gamma_m$. Kulkarni~\cite[Definition~1.1, p.~255]{Kulkarni} puts $\Lambda_{\Kul}=B_\Gamma\cup L_2$ and $\Omega_{\Kul}=X\setminus\Lambda_{\Kul}$.

\begin{proposition}\label{prop:L2}
$\Lambda_{\Kul}=B_\Gamma\cup\overline{\D_\Gamma(X\setminus B_\Gamma)}$, the action on $\Omega_{\Kul}$ is properly discontinuous, and $\mathscr R_\Gamma\subseteq\Lambda_{\Kul}$.
\end{proposition}
\begin{proof}
A generating cluster point $y=\lim\gamma_mx_m$ with $x_m\in K$ gives, after extraction, $x_m\to x\in K$ and $(x,y)\in\D_\Gamma$. Conversely, if $(x,y)\in\D_\Gamma$ and $x\notin B_\Gamma$, the moving sources eventually lie in a compact neighbourhood of $x$ disjoint from the closed set $B_\Gamma$. Hence every destination of a point of $\Omega_{\Kul}$ lies in $\Lambda_{\Kul}$, and Theorem~\ref{thm:relational}(i) gives properness. A recurrent point outside $B_\Gamma$ is its own destination.
\end{proof}

\medskip
\emph{The exceptional set and the Myrberg set.} For a sequence $(\gamma_m)$ of distinct elements put $E_{(\gamma_m)}=\{x:\D_{(\gamma_m)}(x)\text{ has at least two points}\}$, let $\Lambda_{\mathrm{exc}}$ be the union of these sets and put
\[
\Lambda_{\Myr}=\overline{\Lambda_{\mathrm{exc}}},\qquad\Omega_{\Myr}=X\setminus\Lambda_{\Myr}.
\]
The set $\Lambda_{\mathrm{exc}}$ is invariant: if $x$ has two destinations along $(\gamma_m)$, then $\delta x$ has the same destinations along $(\gamma_m\delta^{-1})$. Since both witnesses survive extraction, Theorem~\ref{thm:dest} gives, for extracted Cartan data,
\[
E_{(\gamma_m)}=\{x:\ N_\Theta\text{ does not fix }b\rho(k'x)\}=k'^{-1}\textstyle\bigcup\{O^-_\sigma:\Sigma_+(\sigma)\neq\varnothing\}.
\]
In projective space, Cano--Seade's pseudo-projective boundary consists of the non-invertible limits in $\PP(\mathrm{End}\,V)$ of distinct group elements~\cite{CanoSeade2010}; its kernels form $\Lambda_{\Myr}$ and its complement is the equicontinuity region (Proposition~\ref{prop:kernel}). The name follows Myrberg's classical region of normality~\cite[Chapter~II, \S3, pp.~243--245]{Myrberg}. On a general flag manifold, the exterior boundary $\partial X_\varrho$ in Theorem~\ref{thm:extG} plays the same role.

\begin{lemma}\label{lem:cellfix}
For a cell $\sigma$, $\Sigma_+(\sigma)$ and $\Sigma_-(\sigma)$ are not both empty.
\end{lemma}
\begin{proof}
Otherwise the chart of Lemma~\ref{lem:cells}(iv) at $\sigma Q$ would consist of points of $L_\Theta\sigma Q\subset X^{A_\Theta}$, and the nontrivial group $A_\Theta$ would fix an open subset of the connected analytic manifold $X$, hence act trivially, against effectiveness.
\end{proof}

\begin{theorem}\label{thm:exc}
$\Omega_{\Myr}$ is a domain of proper discontinuity \textup{(}possibly empty\textup{)}.
\end{theorem}
\begin{proof}
Suppose $x,y\in\Omega_{\Myr}$ and $(x,y)\in\D_\Gamma$, witnessed by $(\gamma_m)$ with Cartan data $(\Theta,k,k',b)$. Since $x\notin\Lambda_{\mathrm{exc}}$, $N_\Theta$ fixes $p=b\rho(k'x)=bl\sigma Q$ and $y=kp$. As $bl$ normalises $N_\Theta$, $N_\Theta$ fixes $\sigma Q$, so $\Sigma_+(\sigma)=\varnothing$ and $\Sigma_-(\sigma)\neq\varnothing$ by Lemma~\ref{lem:cellfix}.

Every source $x'\in k'^{-1}N_\Theta^-l\sigma Q$ has $\rho(k'x')=l\sigma Q$ and therefore the single destination $y$. By Theorem~\ref{thm:dest}(d), all these sources are destinations of $y$ for $(\gamma_m^{-1})$. They form a set of dimension $d_-(\sigma)>0$, so $y\in E_{(\gamma_m^{-1})}$, a contradiction. Theorem~\ref{thm:relational}(i) concludes.
\end{proof}

\begin{proposition}\label{prop:kulkarni-myrberg}
$B_\Gamma\subseteq\mathscr R_\Gamma\subseteq\Lambda_{\Kul}\subseteq\Lambda_{\Myr}$. In particular $\Omega_{\Myr}\subseteq\Omega_{\Kul}$.
\end{proposition}
\begin{proof}
Let $U$ be any domain. It contains no point with infinite stabiliser, so $L_0\cap U=\varnothing$. If $y=\lim\gamma_mx\in U$, invariance gives $x\in U$, and for a compact neighbourhood $C\subset U$ of $y$ the compact set $C\cup\{x\}$ meets infinitely many of its translates; so $L_1\cap U=\varnothing$. By Proposition~\ref{prop:recur-equiv}, $B_\Gamma\subseteq\mathscr R_\Gamma$. The middle inclusion is Proposition~\ref{prop:L2}.

For the last one put $\Omega=\Omega_{\Myr}$, a domain by Theorem~\ref{thm:exc}, so $B_\Gamma\cap\Omega=\varnothing$. Suppose a generating cluster point $y$ of $L_2$ lies in $\Omega$: $y=\lim\gamma_mx_m$ with $x_m$ in a compact $K\subset X\setminus B_\Gamma$. Extract so that $x_m\to x\in K$ and $\gamma_m^{-1}y\to z$, and then extract Cartan data for $(\gamma_m^{-1})$. Along this one sequence, both $x$ (with moving sources $\gamma_mx_m\to y$) and $z$ (with the fixed source $y$) are destinations of $y$. Since $y\notin\Lambda_{\mathrm{exc}}$, $z=x$. As $y\notin L_0$, the fixed-source limit $\gamma_m^{-1}y\to x$ puts $x$ in $L_1$, contradicting $x\in K$.
\end{proof}

For complex projective space, the Kulkarni criterion and the pseudo-projective methods of Cano--Navarrete--Seade~\cite[Proposition~3.3.6 and Section~7.4]{CNS} provide the foundations for this comparison. The argument above proves the inclusion $\Lambda_{\Kul}\subseteq\Lambda_{\Myr}$ on every flag manifold.

\medskip
\emph{Quint's limit set.} Let $\varnothing\neq\theta\subseteq\Delta$, $\mathcal F_\theta=G/P_\theta$ and $\nu_\theta$ the $K$-invariant probability measure. For Zariski-dense subgroups, Quint~\cite[author version, \S7.1]{qfQuint2002} defines the limit set on partial flag manifolds by concentration of this measure. We use the same definition for arbitrary discrete $\Gamma$:
\[
\Lambda^Q_\theta(\Gamma)=\{a\in\mathcal F_\theta:\ (\gamma_m)_*\nu_\theta\to\delta_a\text{ weakly, for some distinct }\gamma_m\in\Gamma\}.
\]
In fact concentration already gives a fixed-source limit. If
$(\gamma_m)_*\nu\to\delta_a$, choose increasing $m_j$ such that
$\nu\{x:d(\gamma_{m_j}x,a)\ge1/j\}<2^{-j}$. By countable
subadditivity, for almost every $x$ one has $\gamma_{m_j}x\to a$.
Fix such an $x$. If $x\in L_0$, closedness and invariance imply
$a\in L_0$. Otherwise $a\in L_1$. Thus for every discrete $\Gamma$,
\begin{equation}\label{eq:quint-recurrent}
\Lambda^Q_\theta(\Gamma)\subseteq B_\Gamma\subseteq
\mathscr R_\Gamma\quad\text{on }\mathcal F_\theta.
\end{equation}

\begin{theorem}\label{thm:Quint}
For an arbitrary sequence $(g_m)$ in $G$, $(g_m)_*\nu_\theta\to\delta_a$ if and only if $\min_{\alpha\in\theta}\alpha(\kappa(g_m))\to\infty$ and $k_mP_\theta\to a$, where $\kappa(g_m)$ is the Cartan projection and $g_m=k_m\exp(\kappa(g_m))k'_m$. Consequently $\Lambda^Q_\theta(\Gamma)$ is the compact set of limits $kP_\theta$ of attracting flags of extracted sequences of distinct elements with $\theta\subseteq\Theta$ (for those extracted sequences).
\end{theorem}
\begin{proof}
For each $\alpha\in\theta$, take a Tits representation with one-dimensional highest weight space~\cite[Theorem~7.2]{Tits1971}. In a Cartan-compatible inner product, $K$ is isometric and the restricted weight spaces are orthogonal; the next restricted weight gives
$\sigma_2/\sigma_1=e^{-\alpha(\kappa(g))}$, as in \cite[\S2.5]{Benoist97}. Benoist's concentration lemma~\cite[\S3.5, general case of the proof]{Benoist97}, applied to the highest-weight embeddings of $\mathcal F_\theta$, gives precisely divergence of these gaps and convergence of $k_mP_\theta$. His roots are multiplicative, so their values are $e^{\alpha(\kappa(g_m))}$ in our notation. This sequence criterion requires no Zariski-density hypothesis.

For compactness, if $a_j\to a$ are concentration limits, choose from a sequence realising $a_j$ an unused element $\gamma_j$ with all $\theta$-gaps greater than $j$ and attracting $\theta$-flag within $1/j$ of $a_j$. The equivalence just proved yields concentration at $a$; Cartan extraction gives the stated description.
\end{proof}

For $P_\theta$-divergent groups with opposition-invariant $\theta$, this is the Canary--Zhang--Zimmer limit set~\cite[author version, \S2.3.2]{mlCZZ}. For Zariski-dense groups, minimality follows separately from \cite[\S3.6(ii)]{Benoist97}; this is used in Lemma~\ref{lem:pairs} for the attracting full flags $\Lambda^Q_\Delta(\Gamma)$ of regular groups.

Separate concentration in each Grassmannian does not suffice. In $\mathrm{SL}(4,\R)$ take the discrete diagonal lattice
$\Gamma=\{\exp(mH_1+nH_2):m,n\in\Z\}$, where
\[
 H_1=\tfrac14\diag(5,1,-7,1),\qquad
 H_2=\tfrac14\diag(-1,7,-1,-5).
\]
Positive powers of $\exp H_1$, $\exp(H_1+H_2)$ and $\exp H_2$ attract to the standard $1$-, $2$- and $3$-planes, respectively. But the consecutive diagonal differences of $mH_1+nH_2$ are
\[
 (\alpha_1,\alpha_2,\alpha_3)=(m-2n,\,2m+2n,\,-2m+n),
 \qquad 2\alpha_1+\alpha_2+2\alpha_3=0.
\]
They cannot all tend to $+\infty$, as concentration at the standard full flag requires. Thus the condition $\theta\subseteq\Theta$ must hold jointly.

\medskip
\emph{Cartan levels.} For a divergent sequence with Cartan data $(\Theta,k,k',b)$ and an integer $d\ge0$ put
\[
T_{\le d}(g_m)=k\textstyle\bigcup\{O^+_\tau:\tau\in\mathcal C_\Theta,\ \dim_{\R}O^+_\tau\le d\},\qquad
\Lambda^{(d)}_{\mathrm Z}(\Gamma)=\overline{\textstyle\bigcup_{(\gamma_m)}T_{\le d}(\gamma_m)}.
\]
These sets are closed, increasing in $d$, and $\Gamma$-invariant: left multiplication by $\gamma$ transports the attracting parabolic, by Theorem~\ref{thm:extG}(b), and preserves cell dimensions.

\begin{theorem}\label{thm:L}
For every discrete $\Gamma<G$,
\[
\Lambda^Q_\theta=\Lambda^{(0)}_{\mathrm Z}\subseteq\Lambda^{(1)}_{\mathrm Z}\subseteq\cdots\subseteq\Lambda^{(\dim X-1)}_{\mathrm Z}=\Lambda_{\Myr},
\qquad
\Lambda^Q_\theta\subseteq B_\Gamma\subseteq\mathscr R_\Gamma\subseteq\Lambda_{\Kul}\subseteq\Lambda_{\Myr}.
\]
For every dimension system $\delta$ of Section~\ref{sec:domains}, $\Omega_\delta$ is a domain,
$\mathscr R_\Gamma\subseteq\Lambda_\delta$, and $\Lambda^{(d)}_{\mathrm Z}\subseteq\Lambda_\delta$ for $d\le\min_\Theta\delta_\Theta$.
Corollary~\ref{cor:trichotomy} applies; in particular, $\Lambda_\delta=\mathscr R_\Gamma$ when $\Omega_\delta$ is largest.
\end{theorem}
\begin{proof}
A cell $O^+_\tau$ is a point exactly when $P_\Theta\subseteq\tau P_\theta\tau^{-1}$. The latter parabolic then contains $B$, so is standard and conjugate to $P_\theta$, hence equals it; thus $\tau Q=Q$ and $\theta\subseteq\Theta$. The point cell is $kP_\theta$, and Theorem~\ref{thm:Quint} gives the first endpoint.

For the other endpoint, $T_{\le\dim X-1}(g_m)=k(X\setminus O^+_{\mathrm{open}})$. Lemma~\ref{lem:cellfix} and $d_++d_0+d_-=\dim X$ give
$E_{(g_m)}=k'^{-1}(X\setminus O^-_{\mathrm{open}})$.
The inverse Cartan decomposition
$g_m^{-1}=(k_m'^{-1}\dot w_0)e^{\iota\mu_m}(\dot w_0^{-1}k_m^{-1})$, with $\iota\mu=-w_0\mu$, has type $\iota\Theta$ and
$\dot w_0P_{\iota\Theta}\dot w_0^{-1}=P^-_\Theta$.
Hence $T_{\le\dim X-1}(g_m^{-1})=E_{(g_m)}$, and $\Gamma=\Gamma^{-1}$ identifies the unions. Monotonicity gives the remaining Cartan inclusions.

The second chain is \eqref{eq:quint-recurrent} and Proposition~\ref{prop:kulkarni-myrberg}. The dimension-system assertions follow from their definition, Theorem~\ref{thm:balanced} and Corollary~\ref{cor:trichotomy}.
\end{proof}

\section{Marked and unmarked Levi correspondences}\label{sec:levi-recovery}

\subsection{Marked graphs and the unmarked destination relation}

Let $X$ be a nonempty compact metrizable space and $\Gamma$ an infinite countable group acting on it by homeomorphisms. Write $\G_\gamma=\{(x,\gamma x):x\in X\}\subset X\times X$ and let $\mathcal K(X\times X)$ be the Hausdorff hyperspace of nonempty compact subsets. Define the \emph{graph-limit boundary}
\[
 \Lims_\Gamma(X)=\bigcap_{F\subset\Gamma,\ F\text{ finite}}
     \overline{\{\G_\gamma:\gamma\in\Gamma\setminus F\}}
     \ \subset\mathcal K(X\times X).
\]
Its \emph{marked incidence space} is
\[
 \I_\Gamma(X)=\{(R,x,y):R\in\Lims_\Gamma(X),\ (x,y)\in R\}.
\]
Forgetting the individual limiting graph recovers the relation $\D_\Gamma$ of Section~\ref{sec:groups}.

\begin{proposition}\label{prop:marked}
The spaces $\Lims_\Gamma(X)$ and $\I_\Gamma(X)$ are nonempty compact metrizable spaces. Projection along the marking satisfies
\begin{equation}\label{eq:mark-projection}
 \D_\Gamma=\operatorname{pr}_{X\times X}\I_\Gamma(X)
                 =\bigcup_{R\in\Lims_\Gamma(X)}R.
\end{equation}
The marking records which pairs arise \emph{along the same degeneration}, information generally lost under projection.
\end{proposition}
\begin{proof}
The compact metrizable hyperspace has nonempty closed tails whose finite intersections are nonempty. Their intersection is therefore nonempty and compact; the incidence condition $(x,y)\in R$ is closed under Hausdorff convergence, so $\I_\Gamma(X)$ is compact as well.

For $(x,y)\in\D_\Gamma$, extract a Hausdorff limit $R$ of the witnessing graphs. Distinctness puts $R$ in $\Lims_\Gamma(X)$ and $(x,y)\in R$. Conversely, the closed-tail definition lets us choose distinct $\gamma_m$ with $\G_{\gamma_m}\to R$; Hausdorff convergence then supplies $(x_m,\gamma_mx_m)\to(x,y)$ for every $(x,y)\in R$. This proves \eqref{eq:mark-projection}; closedness of $\D_\Gamma$ was established in Lemma~\ref{lem:Dbasic}.
\end{proof}

\subsection{The effective Levi action on a full-flag transversal}

Let $X=G/B$ be the full real flag manifold and let $(g_m)$ have
Cartan data $(\Theta,k,k',b)$ as in Lemma~\ref{lem:split}. Put
\[
 B_\Theta=L_\Theta\cap B,\qquad
 Y_\Theta=L_\Theta/B_\Theta,\qquad
 C_\Theta=\ker\bigl(L_\Theta\curvearrowright Y_\Theta\bigr).
\]
We identify $Y_\Theta$ with the fixed component $L_\Theta B/B$.
On the open cell $P^-_\Theta B/B$, the stratified retraction is
$\rho(n^-lB)=lB$, for $n^-\in N^-_\Theta$ and $l\in L_\Theta$.

\begin{theorem}\label{thm:levi-transversal}
On $k'^{-1}P^-_\Theta B/B$, the whole-sequence destination is the
singleton given by
\begin{equation}\label{eq:generic-map}
 f_b(x)=kb\rho(k'x).
\end{equation}
In particular,
\begin{equation}\label{eq:levi-graph}
 \D_{(g_m)}\cap\bigl((k'^{-1}Y_\Theta)\times X\bigr)
 =\{(k'^{-1}y,kby):y\in Y_\Theta\}.
\end{equation}
For fixed type and fixed source and target frames, the transverse
relations associated to $b_1,b_2\in L_\Theta$ coincide if and only if
$b_2^{-1}b_1\in C_\Theta$. Thus this transverse relation recovers
exactly the class in the effective Levi group $L_\Theta/C_\Theta$.
\end{theorem}
\begin{proof}
Theorem~\ref{thm:dest} gives
$\D_{(g_m)}(x)=k\overline{N_\Theta b\rho(k'x)}$.
On the stated cell, $\rho(k'x)=lB$. Since $bl\in L_\Theta$
normalises $N_\Theta\subset B$, this unipotent group fixes $blB$.
The orbit is therefore a singleton. On $Y_\Theta$ the retraction
is the identity, giving~\eqref{eq:levi-graph}. Finally the two graphs
coincide precisely when $b_1y=b_2y$ for every $y\in Y_\Theta$,
which is the asserted kernel condition.
\end{proof}

\begin{remark}\label{rem:effective}
Lemma~\ref{lem:levi-scalar} gives $Z_\Theta\subseteq C_\Theta$.
The fibre $L_\Theta/Z_\Theta$ in Theorem~\ref{thm:extG} belongs to
the chosen projective compactification; the effective action recovered
here is its quotient $L_\Theta/C_\Theta$. For real groups the latter
kernel can contain compact normal factors invisible on real flags.
For example, the full real flag manifold of
$\PSL(2,\R)\times\mathrm{SO}(3)$ is $\R\PP^1$, on which the second
factor acts trivially, whereas the projectively faithful tensor
representation used in Section~\ref{sec:extG} retains that factor in
its rank-one degenerations of the first factor. Thus the chosen
boundary is sufficient to determine the destination relation, but
need not be its minimal effective parameter space.
\end{remark}

\subsection{Joint realisation of boundary data}

Let $\Gamma_0<\mathrm{SL}_2(\C)$ be discrete, non-elementary and satisfy $-I\notin\Gamma_0$. Its classical limit set $\Lambda\subset\PP^1$ is compact, infinite and minimal for the action. Let $r\ge2$, $U=\C^2$, $V=\C^r$ and let $\chi:\Gamma_0\to U(r)$ be a homomorphism. Put $H=\overline{\chi(\Gamma_0)}$, a compact Lie subgroup. The induced projective action on
\[
 X_r=\PP(U\otimes V^*)=\PP(\Hom(V,U))
\]
is $\gamma[A]=[\gamma A\chi(\gamma)^{-1}]$. We denote its image by $\Gamma_\chi$. For $p\in\PP(U)$ set $E_p=p\otimes V^*$ and $W_\Lambda=\bigcup_{p\in\Lambda}\PP(E_p)$.
Our convention is the contragredient action $\chi^\vee=(\chi^{-1})^t$ on $V^*$. Thus a twist $\gamma\otimes\vartheta(\gamma)$ on $U\otimes W$ is obtained by taking $V=W^*$ and $\chi=\vartheta^\vee$; the same results apply. The compact group $H$ need not be connected.

For $q\in\PP(U)$ choose $\ell_q\in U^*\setminus0$ with kernel $q$. For $p\in\PP(U)$ and $[\eta]\in\PP(V^*)$, define
\begin{equation}\label{eq:F}
 F(p,[\eta])=\PP\bigl(E_p+\C(v\otimes\eta)\bigr),\qquad v\notin p.
\end{equation}
The definition is independent of $v$ (and of all projective representatives).

A sequence $\gamma_m\in\Gamma_0$ is \emph{escaping} if it is pairwise distinct. Its normalised matrices have, after extraction, rank-one limits $[T_{p,q}]$ with image $p$ and kernel $q$, and the compact factors $\chi(\gamma_m)$ have cluster points in $H$.
The following point is stronger than independent density of these
three parameters: their values can be prescribed along a single
escaping sequence.

\begin{lemma}\label{lem:saturate}
For every $p,q\in\Lambda$ and every $h\in H$ there exists an escaping sequence $(\gamma_m)$ with
\begin{equation}\label{eq:joint}
 [\gamma_m]\longrightarrow[T_{p,q}]
 \quad\text{in }\PP(\operatorname{End}U),
 \qquad \chi(\gamma_m)\longrightarrow h.
\end{equation}
No parabolic element, finite-generation hypothesis or geometric-finiteness hypothesis is needed.
\end{lemma}
\begin{proof}
Choose a loxodromic $a$, diagonal with eigenvalues $\lambda,\lambda^{-1}$, $|\lambda|>1$, and eigenlines $a^\pm$. Put
$B_0=\{w:w(a^-)=a^-\}$ and $S=\Gamma_0\setminus B_0=\{w:w_{12}\ne0\}$.
A loxodromic $c$ not fixing $a^-$ has return powers $j_m\to\infty$ with $\chi(c)^{j_m}\to I$.
For $w\in B_0$, $wc^{j_m}\in S$ and $\chi(wc^{j_m})\to\chi(w)$; hence $\chi(S)$ is dense in $H$.

For fixed $w\in S$,
\[
 a^nwa^{-n}=\begin{pmatrix}
 w_{11}&\lambda^{2n}w_{12}\\
 \lambda^{-2n}w_{21}&w_{22}
 \end{pmatrix}
 \quad\Longrightarrow\quad [a^nwa^{-n}]\to[E_{12}],
\]
whose image and kernel are both $a^+$. Along return powers $n_j\to\infty$ with $\chi(a)^{n_j}\to I$, the compact factors tend to $\chi(w)$. Left and right multiplication by $g_1,g_2$ changes the endpoints to $g_1a^+,g_2^{-1}a^+$ and the compact limit to $\chi(g_1)\chi(w)\chi(g_2)$.

Use the density of $\Gamma_0a^+$ to choose $g_{1,m},g_{2,m}$ with endpoints within $1/m$ of $p,q$, and choose $w_m\in S$ with
$\|\chi(w_m)-\chi(g_{1,m})^{-1}h\chi(g_{2,m})^{-1}\|<1/(2m)$.
After fixing these choices, take a sufficiently large return power $n=n_{j(m)}$ and put
$\gamma_m=g_{1,m}a^nw_ma^{-n}g_{2,m}$ so that simultaneously
\[
 \begin{gathered}
 d\bigl([\gamma_m],[g_{1,m}E_{12}g_{2,m}]\bigr)<1/m,\qquad
 \|\chi(a^nw_ma^{-n})-\chi(w_m)\|<1/(2m),\\
 \|\gamma_m\|>m+\max_{i<m}\|\gamma_i\|,
 \end{gathered}
\]
with empty maximum zero. The outer compact factors are unitary, so the compact errors total less than $1/m$. The endpoint and projective-matrix approximations give $[\gamma_m]\to[T_{p,q}]$, including when $p=q$, while the norm bound ensures distinctness.
\end{proof}

\subsection{The boundary as a space of correspondences}

For $p,q\in\Lambda$ and a projective unitary $[h]\in\PU(r)$, define the closed relation $R_{p,q,[h]}\subset X_r\times X_r$ by its fibres:
\begin{equation}\label{eq:R-fibres}
 R_{p,q,[h]}([A])=
 \begin{cases}
 \{[v_p\otimes(\ell_q A)h^{-1}]\},&\ell_q A\ne0,\\[2pt]
 F(p,[\beta h^{-1}]),& [A]=[v_q\otimes\beta]\in\PP(E_q).
 \end{cases}
\end{equation}
Here $v_p\in p\setminus0$; the relations do not depend on the chosen representatives of $p,q,h$. Let $PH$ denote the image of $H$ in $\PU(r)$.

\begin{theorem}\label{thm:marked-product}
The group $\Gamma_\chi$ is faithful and discrete. Every escaping sequence in $\Gamma_\chi$, after extraction, has singular Cartan type $(r,r)$ and a graph limit of the form $R_{p,q,[h]}$. Conversely, every triple in $\Lambda\times\Lambda\times PH$ occurs. The map
\begin{equation}\label{eq:boundary-product}
 \Lambda\times\Lambda\times PH\longrightarrow\Lims_{\Gamma_\chi}(X_r),
 \qquad(p,q,[h])\longmapsto R_{p,q,[h]}
\end{equation}
is a homeomorphism.
In particular, for fixed $(p,q)$, the \emph{marked} boundary recovers the projective compact Levi subgroup $PH$ as a space of transformations. The \emph{unmarked} relation $\D_{\Gamma_\chi}=\bigcup_{p,q,[h]}R_{p,q,[h]}$ can identify different such subgroups.
\end{theorem}
\begin{proof}
For the Hilbert--Schmidt metric, the singular values are $s,s^{-1}$, those of $\gamma$, each repeated $r$ times. Discreteness of $\Gamma_0$ gives $s\to\infty$ along escaping sequences, proving type $(r,r)$ and discreteness of the projective image. A scalar tensor operator has both factors scalar, and $-I\notin\Gamma_0$ then gives faithfulness.

For a sequence with $[\gamma_m]\to[T_{p,q}]$, $\chi(\gamma_m)\to h$, normalisation gives a singleton limit on $X_r\setminus\PP(E_q)$. On $\PP(E_q)$, apply a singular-value decomposition of $\gamma_m$ with scales $s_m,s_m^{-1}$, $s_m\to\infty$. In the source chart, an ordinary matrix has rows $(\alpha_m,\beta_m)$, with $\alpha_m\to0$ and $\beta_m\to\beta\ne0$. After division by the slow scale its image has rows $(s_m^2\alpha_m,\beta_m)h_m^{-1}$, up to convergent unitary frames. By choosing $\alpha_m$ of order $s_m^{-2}$ one obtains every finite point of the projective span in~\eqref{eq:R-fibres}; allowing $s_m^2\alpha_m$ to tend slowly to infinity gives all points of $\PP(E_p)$. No additional limit is possible, since the slow row converges projectively to $\beta h^{-1}$. These choices work along the whole extracted sequence, and the graphs converge in the Hausdorff topology to $R_{p,q,[h]}$.

Lemma~\ref{lem:saturate} realises every triple, proving surjectivity in~\eqref{eq:boundary-product}. For injectivity, the generic singleton destinations identify $p$, and the exceptional source plane with positive-dimensional destination identifies $q$. With these fixed, the values of singleton destinations as $\ell_q A$ ranges over $V^*$ recover the projective transformation $[h]$: if two such maps agree on all $[\ell_q A]$, they differ by a scalar.

For continuity, let $(p_j,q_j,[h_j])\to(p,q,[h])$. By Lemma~\ref{lem:saturate} choose distinct $\gamma_j$, with $\|\gamma_j\|>j$, whose projective matrix, projective compact factor and graph are within $1/j$ of $[T_{p_j,q_j}]$, $[h_j]$ and $R_{p_j,q_j,[h_j]}$, respectively. Every subsequence admits convergent unitary frames, and the calculation above gives the same graph limit $R_{p,q,[h]}$. Compactness of the hyperspace gives convergence of the full graph sequence, hence of $R_{p_j,q_j,[h_j]}$. The continuous bijection from the compact parameter space to the Hausdorff hyperspace is a homeomorphism.
\end{proof}

\subsection{Unmarked incidence, orbit quotients and the largest domain}

We now forget which individual graph limit produced a given pair.
The resulting union retains the compact orbit partition even though
it need not retain the compact subgroup itself.

\begin{theorem}\label{thm:unmarked}
For any $[A]\in X_r\setminus W_\Lambda$,
\begin{equation}\label{eq:outside}
 \D_{\Gamma_\chi}([A])=
 \{[v_p\otimes(\ell_q A)h^{-1}]:p,q\in\Lambda,\ h\in H\}.
\end{equation}
For $[A]=[v_q\otimes\beta]\in W_\Lambda$,
\begin{equation}\label{eq:inside}
 \D_{\Gamma_\chi}([A])=
 \bigcup_{p\in\Lambda,\ h\in H} F(p,[\beta h^{-1}]).
\end{equation}
Consequently, the recurrent set and the largest domain are
\begin{equation}\label{eq:maximal}
 \mathscr R_{\Gamma_\chi}=W_\Lambda,
 \qquad\Omega_{\max}(\Gamma_\chi)=X_r\setminus W_\Lambda.
\end{equation}
The latter is nonempty for $r\ge2$ and independent of $H$.

If $\Lambda\subsetneq\PP^1$, fix $s\in\PP^1\setminus\Lambda$ and let
\[
 \Sigma_s=\{[v_s\otimes\beta]:[\beta]\in Y\}\simeq Y,
 \qquad Y=\PP(V^*),\quad Q_H=Y/H.
\]
With $W_\Lambda\simeq\Lambda\times Y$ by the Segre embedding, one has the \emph{fibre-product correspondence}
\begin{equation}\label{eq:fibre-product}
 \D_{\Gamma_\chi}\cap(\Sigma_s\times W_\Lambda)
 \ \simeq\ Y\mathop{\times}_{Q_H}(\Lambda\times Y)
 =\{(\beta,(p,\eta)):[\beta]_H=[\eta]_H\}.
\end{equation}
Thus, under the hypothesis $\Lambda\subsetneq\PP^1$, $\D_{\Gamma_\chi}$ determines the compact Hausdorff orbit space $Q_H$ and the orbit partition of $Y$. For two twists of the same $\Gamma_0$ with this proper limit set, the unmarked relations coincide if and only if their projective compact closures have the same orbit partition on $Y$. Under the same proper-limit-set hypothesis, the relation is the full two-sided star
\[
 (X_r\times W_\Lambda)\cup(W_\Lambda\times X_r)
\]
if and only if $H$ acts transitively on $Y$.
\end{theorem}
\begin{proof}
Taking the union of \eqref{eq:R-fibres} over the marked boundary gives \eqref{eq:outside}--\eqref{eq:inside}. Destinations of sources off $W_\Lambda$ lie in $W_\Lambda$, whereas $W_\Lambda\times W_\Lambda\subset\D_{\Gamma_\chi}$ because every fibre on the wall contains every $\PP(E_p)$. Thus Theorem~\ref{thm:relational} gives \eqref{eq:maximal}; rank-two matrices prove nonemptiness.

For $[A]=[v_s\otimes\beta]$ with $s\notin\Lambda$, the covector $\ell_q A$ is a nonzero scalar multiple of $\beta$ for every $q\in\Lambda$. Thus~\eqref{eq:outside} reduces to $\Lambda\times H[\beta]$. This is exactly~\eqref{eq:fibre-product}. Since $H$ is compact, its orbit equivalence relation is closed, and $Q_H$ is compact Hausdorff. Two such relations agree precisely when all their fibres $H[\beta]$ agree. If $H$ is transitive,~\eqref{eq:outside} yields the full wall for sources off the wall. For sources on the wall, any target not in $\PP(E_p)$ has a nonzero quotient covector modulo $p$, which transitivity makes proportional to $\beta h^{-1}$ for some $h\in H$; formula~\eqref{eq:inside} then gives all targets. Conversely, the star relation and~\eqref{eq:fibre-product} force every $H$-orbit in $Y$ to be all of $Y$.
\end{proof}

\begin{remark}
For $\Lambda=\PP^1$, the explicit marked-boundary and complete destination formulas, and the largest-domain assertion, still hold. The transverse slice $\Sigma_s$ used to recover $Q_H$ no longer exists. The theorem makes no claim that the unmarked relation recovers the orbit partition in this case.
\end{remark}

\subsection{Rigidity under projective equivalence}

The orbit partition can also be recovered without fixing tensor
coordinates. Retain $U=\C^2$, $V=\C^r$, $r\ge2$, and the
non-elementary group $\Gamma_0$ above, and assume
$\Lambda\subsetneq\PP(U)$.

\begin{theorem}\label{thm:projective-rigidity}
Let $\chi_i:\Gamma_0\to U(r)$ have compact closures $H_i$, acting
contragrediently on $V^*$. If $T\in\PGL(U\otimes V^*)$ satisfies
\[
 (T\times T)\D_{\Gamma_{\chi_1}}=\D_{\Gamma_{\chi_2}},
\]
then $T=a\otimes c$ with $a\in\PGL(U)$, $a\Lambda=\Lambda$, and
$c\in\PGL(V^*)$, and
\begin{equation}\label{eq:orbit-conjugacy}
 c(H_1[\beta])=H_2(c[\beta])
 \quad\text{for every }[\beta]\in\PP(V^*).
\end{equation}
Conversely, every such pair $(a,c)$ induces a projective equivalence
of the complete destination relations.
\end{theorem}
\begin{proof}
By~\eqref{eq:maximal}, the diagonal of either relation is
$W_\Lambda$, so $T(W_\Lambda)=W_\Lambda$. The set $\Lambda$ is
infinite, hence Zariski dense in the complex projective line.
Consequently the complex Zariski closure of $W_\Lambda$ is the
Segre variety $S=\PP(U)\times\PP(V^*)$, which $T$ preserves.
Its two families of maximal linear subspaces have different
dimensions if $r>2$, and therefore are preserved. If $r=2$, an
exchange is also impossible: $W_\Lambda=\Lambda\times\PP^1$
contains complete fibres of one family and none of the other,
because $\Lambda$ is proper. Thus $T=a\otimes c$ and
$a\Lambda=\Lambda$.

For $s\notin\Lambda$, formula~\eqref{eq:fibre-product} identifies
the destination fibre over $[v_s\otimes\beta]$ with
$\Lambda\times H_1[\beta]$. Transporting it by $T$ and comparing
with the second relation gives~\eqref{eq:orbit-conjugacy}.
Conversely, $\ell_{aq}\circ a$ is proportional to $\ell_q$.
Hence~\eqref{eq:outside} is equivariant under $(a,c)$ whenever
\eqref{eq:orbit-conjugacy} holds. Formula~\eqref{eq:inside} is
equivariant as well, since
$(a\otimes c)F(p,[\eta])=F(ap,[c\eta])$.
These two formulae exhaust the relation and prove the converse.
\end{proof}

\subsection{Different compact groups with the same unmarked relation}

\begin{proposition}\label{prop:sharpness}
There exist two unitary twists of the same non-elementary free Kleinian group into $U(4)$ whose image closures are $H_1=\mathrm{Sp}(2)$ and $H_2=\mathrm{SU}(4)$ (in their standard inclusions into $U(4)$). Their complete unmarked destination relations coincide, whereas their marked graph-limit boundaries are different subsets of $\mathcal K(X_4\times X_4)$.
\end{proposition}
\begin{proof}
Take a classical non-elementary Schottky group $\Gamma_0\simeq F_2<\mathrm{SL}_2(\C)$ with $-I\notin\Gamma_0$. Dense two-generator homomorphisms $F_2\to\mathrm{Sp}(2)$ and $F_2\to\mathrm{SU}(4)$ exist because every connected compact semisimple Lie group is topologically generated by two elements~\cite{Kuranishi}; see also~\cite{BG,BGErratum}. Both compact groups act transitively on the unit sphere in $\C^4$ (the first in its quaternionic model), hence on $\PP^3$. The star criterion of Theorem~\ref{thm:unmarked} gives identical unmarked relations. However their projective images in $\PU(4)$ are different. By the injective marked-boundary parametrisation~\eqref{eq:boundary-product}, their graph-limit sets differ: for any $p,q\in\Lambda$, a projective unitary transformation outside $P\mathrm{Sp}(2)$ gives a graph limit for the second but not for the first.
\end{proof}

\begin{remark}
The compact twist leaves the largest projective domain unchanged. On $\operatorname{Gr}(2,U\otimes V)$, the same residual orbit geometry can change whether a largest domain exists; see Section~\ref{grlevi:sec:grassmann-levi}, using the contragredient convention above.
\end{remark}

\section{Balanced cells and properness domains}\label{sec:domains}

The destination formula lets us construct properness domains by covering the possible source--target cell incidences.

Fix a nonempty type $\Theta$ and write $\mathcal C_\Theta=W_\Theta\backslash W/W_\theta$. For $\sigma\in\mathcal C_\Theta$, let $E(\sigma)$ be the set of labels $\tau$ for which $O^+_\tau$ meets $\overline{N_\Theta\sigma Q/Q}$. A pair of subsets $\mathfrak S_\Theta,\mathfrak T_\Theta\subseteq\mathcal C_\Theta$ is a vertex cover if
\[
\sigma\notin\mathfrak S_\Theta,\quad \tau\in E(\sigma)
\quad\Longrightarrow\quad \tau\in\mathfrak T_\Theta.
\]
A \emph{balanced system} is a choice of such a pair for every nonempty Cartan type.

For a divergent sequence with limiting compact factors $k,k'$, declare the source cells $k'^{-1}O^-_\sigma$ with $\sigma\in\mathfrak S_\Theta$ bad, and the target cells $kO^+_\tau$ with $\tau\in\mathfrak T_\Theta$ bad. For a discrete subgroup $\Gamma<G$, let $\Lambda_\mathfrak B(\Gamma)$ be the closure of the union of these bad source and target sets over all divergent sequences in $\Gamma$.

\begin{theorem}\label{thm:balanced}
For every balanced system $\mathfrak B$, the set $\Lambda_\mathfrak B(\Gamma)$ is closed and $\Gamma$-invariant, and every pair in $\D_\Gamma$ has at least one endpoint in $\Lambda_\mathfrak B(\Gamma)$. Consequently,
\[
\Omega_\mathfrak B=X\setminus\Lambda_\mathfrak B(\Gamma)
\]
is an open invariant set on which the action is properly discontinuous; if it is nonempty, it is a domain of proper discontinuity.
\end{theorem}
\begin{proof}
Take $(x,y)\in\D_\Gamma$ and extract Cartan data $(\Theta,k,k',b)$ for a witnessing sequence. If $k'x$ lies in a bad source cell, then $x\in\Lambda_\mathfrak B(\Gamma)$.

Otherwise its cell label $\sigma$ is not in $\mathfrak S_\Theta$. Write $\rho(k'x)=l\sigma Q$ with $l\in L_\Theta$. The destination formula gives $\D_{(\gamma_m)}(x)=kbl\,\overline{N_\Theta\sigma Q}$. Since $bl\in L_\Theta$ preserves every $P_\Theta$-orbit, the target labels are exactly $E(\sigma)$. The cover condition makes all these target cells bad, so $y\in\Lambda_\mathfrak B(\Gamma)$.

The set is invariant because left and right multiplication of the witnessing sequence transports its attracting and repelling parabolics equivariantly; it is closed by definition. The relational criterion now gives proper discontinuity on the complement.
\end{proof}

When $\Theta=\Delta$, the cell closure order is the Bruhat order, and the construction specialises to the thickening domains of Kapovich--Leeb--Porti. The formulation here keeps track of the destination fibre and allows the cell choice to vary with the Cartan type. Cocompactness is not addressed here.

\begin{corollary}\label{cor:trichotomy}
Let $\Lambda=\Lambda_{\mathfrak B}(\Gamma)$ be a balanced wall and $\mathscr R=\mathscr R_\Gamma$. Then $\mathscr R\subseteq\Lambda$, with three alternatives:
\textup{(i)} $\Lambda=\mathscr R$, so $\Omega_{\mathfrak B}$ is largest and $\Lambda$ is the smallest vertex cover;
\textup{(ii)} $\mathscr R\subsetneq\Lambda$ and $\mathscr R$ is a vertex cover, so the largest domain is $X\setminus\mathscr R\supsetneq\Omega_{\mathfrak B}$;
\textup{(iii)} $\mathscr R$ is not a vertex cover, so two non-recurrent points are related and no largest domain exists.
\end{corollary}
\begin{proof}
Apply Corollary~\ref{cor:covers} and Theorem~\ref{thm:relational} to the closed invariant vertex cover $\Lambda$.
\end{proof}

Since $\sigma Q\in N_\Theta\sigma Q$, diagonal loops force
$\mathfrak S_\Theta\cup\mathfrak T_\Theta=\mathcal C_\Theta$.
A vertex cover with complementary source and target sets is therefore minimal.
Every down-set $\mathcal J_\Theta$ defines such a cover
$(\mathcal C_\Theta\setminus\mathcal J_\Theta,\mathcal J_\Theta)$, since
$E(\sigma)\subseteq\{\tau:\tau\le\sigma\}$ by Lemma~\ref{lem:schubert-cell}.
We call these \emph{sublevel systems}.

\medskip
\emph{The exceptional system} takes $\mathcal J_\Theta$ to be the set of closed cells, those with $\Sigma_+(\sigma)=\varnothing$. Its bad sources are the sources with more than one destination. Under the hypotheses of Section~\ref{sec:chain}, its wall is the Myrberg set $\Lambda_{\Myr}$. \emph{Dimension systems} take $\mathcal J_\Theta=\{\tau:\dim O^+_\tau\le\delta_\Theta\}$ for integers $\delta_\Theta$; a cell below $\sigma$ has smaller dimension, so this is a down-set. We write $\Lambda_\delta$ and $\Omega_\delta=X\setminus\Lambda_\delta$ for the resulting wall and domain. The critical wall of Section~\ref{sec:critical} is a variant in projective space which pairs a middle target level with the complementary level of the inverse sequence.

\medskip
\emph{The cell graph is the Bruhat order.} For all types the edges of the cell graph are exactly the order relations.

\begin{proposition}\label{prop:Eorder}
For every nonempty $\Theta\subseteq\Delta$ and every $\sigma\in\mathcal C_\Theta$, $E(\sigma)=\{\tau\in\mathcal C_\Theta:\tau\le\sigma\}$.
\end{proposition}
\begin{proof}
For the minimal representative $\sigma$, Lemma~\ref{lem:schubert-cell} gives
$\overline{N_\Theta\sigma Q}=\overline{B\sigma Q}=\bigcup_{\tau'\le\sigma}B\tau'Q$.
Each $B\tau'Q$ lies in the $P_\Theta$-orbit of its double coset. Projection is order preserving; conversely, a lower double coset has a representative $\tau'\le\sigma$ with this upper representative fixed~\cite[Lemma~4.1.15(i),(iii)]{SteckerThesis}. Its minimal representative satisfies $\tau\le\tau'\le\sigma$ by the parabolic representative properties~\cite[Propositions~2.4.4 and~2.5.1; cf.\ Chapter~2, Exercise~15]{BjornerBrenti}. These two implications identify the orbit closure order with the asserted edges.
\end{proof}

For the diagonal symmetric-pair interpretation of this order, see \cite[Lemma~10.1.1 and Remark~10.1.3]{RichardsonSpringer1990}.

So for each type the minimal balanced systems are exactly the sublevel systems. Indeed, if $v\in\mathfrak S_\Theta\cap\mathfrak T_\Theta$ in a minimal cover, removing its source copy would leave an edge $v\to w$ with $w\notin\mathfrak T_\Theta$ uncovered; removing its target copy would leave an edge $u\to v$ with $u\notin\mathfrak S_\Theta$ uncovered. Thus $w\le v\le u$, and Proposition~\ref{prop:Eorder} gives an uncovered edge $u\to w$, a contradiction. The two sets are therefore complementary, and the cover condition makes $\mathfrak T_\Theta$ a down-set. The diagonal matching gives the lower bound $|\mathcal C_\Theta|$, attained by every such complementary cover.

\medskip
\emph{The dictionary with thickenings.} Let $\Theta_0\subseteq\Delta$ be invariant under the opposition involution $\iota=-w_0$. A discrete $\Gamma$ is \emph{$\Theta_0$-convergent} if every divergent sequence has a subsequence of type $\Theta\supseteq\Theta_0$; let $\Lambda_{\Theta_0}\subset G/P_{\Theta_0}$ be the closed set of its attracting flags, which contains the repelling ones. Put $\mathcal C_{\Theta_0}=W_{\Theta_0}\backslash W/W_\theta$. For $\Theta\supseteq\Theta_0$ the projection $\pi:\mathcal C_\Theta\to\mathcal C_{\Theta_0}$ is surjective and order preserving, and $w_0$ acts on $\mathcal C_{\Theta_0}$ as an order-reversing involution. For a down-set $\mathrm{Th}\subseteq\mathcal C_{\Theta_0}$ and $\xi=gP_{\Theta_0}$ put $\mathrm{Th}(\xi)=\bigcup_{[\tau']\in\mathrm{Th}}gP_{\Theta_0}\tau'Q$. Following Kapovich--Leeb--Porti~\cite[\S3.4]{KLP}, $\mathrm{Th}$ is \emph{fat} if $\mathrm{Th}\cup w_0\mathrm{Th}=\mathcal C_{\Theta_0}$ and \emph{balanced} if moreover $\mathrm{Th}\cap w_0\mathrm{Th}=\varnothing$.

\begin{theorem}\label{thm:KLPfull}
Let $\Gamma$ be $\Theta_0$-convergent and $\mathrm{Th}\subseteq\mathcal C_{\Theta_0}$ a down-set. For $\Theta\supseteq\Theta_0$ put $\mathfrak T_\Theta=\pi^{-1}(\mathrm{Th})$ and $\mathfrak S_\Theta=\pi^{-1}(w_0\mathrm{Th})$.
\begin{enumerate}[label=\textup{(\roman*)},nosep]
\item These pairs are vertex covers for all $\Theta\supseteq\Theta_0$ if and only if $\mathrm{Th}$ is fat; they are complementary if and only if $\mathrm{Th}$ is balanced.
\item In that case $\Lambda_{\mathfrak B}(\Gamma)=\bigcup_{\xi\in\Lambda_{\Theta_0}}\mathrm{Th}(\xi)$, and $\Gamma$ acts properly discontinuously on its complement.
\end{enumerate}
More generally, for arbitrary $\Theta_0$, a down-set $I\subseteq W$ with $W_{\Theta_0} I=I=IW_\theta$ and $I^\vee=W\setminus w_0I$, every sequence of type $\Theta\supseteq\Theta_0$ with flags $\xi^+\in G/P_{\Theta_0}$, $\xi^-\in G/P_{\iota\Theta_0}$ satisfies
\begin{equation}\label{eq:absorption}
x\notin\mathrm{Th}_{I^\vee}(\xi^-)\ \Longrightarrow\ \D_{(g_m)}(x)\subseteq\mathrm{Th}_I(\xi^+),
\end{equation}
so $\bigcup_{\xi^-}\mathrm{Th}_{I^\vee}(\xi^-)\cup\bigcup_{\xi^+}\mathrm{Th}_I(\xi^+)$ is a vertex cover, with no fatness hypothesis.
\end{theorem}
\begin{proof}
(i) For $\Theta=\Delta$ the loops and the surjectivity of $\pi$ show that the vertex-cover condition forces $\mathcal C_{\Theta_0}\setminus w_0\mathrm{Th}\subseteq\mathrm{Th}$. Conversely, if $\sigma\notin\mathfrak S_\Theta$ then $\pi(\sigma)\in\mathrm{Th}$, and every $\sigma'\in E(\sigma)$ has $\sigma'\le\sigma$ by Proposition~\ref{prop:Eorder}, so $\pi(\sigma')\in\mathrm{Th}$. The second assertion is immediate.

(ii) For data $(\Theta,k,k',b)$ the flags are $\xi^+=kP_{\Theta_0}$ and $\xi^-=k'^{-1}w_0P_{\Theta_0}$. The bad targets are $k\bigcup_{\mathfrak T_\Theta}P_\Theta\sigma Q=\mathrm{Th}(\xi^+)$. The bad sources are $k'^{-1}P^-_\Theta\sigma Q=(k'^{-1}w_0)P_{\iota\Theta}w_0\sigma Q$, with label $w_0\pi(\sigma)$, so they form $\mathrm{Th}(\xi^-)$. Each $\mathrm{Th}(\xi)$ is a union of Schubert varieties and depends continuously on $\xi$, so the union over the compact set $\Lambda_{\Theta_0}$ is closed. Theorem~\ref{thm:balanced} gives properness. The last assertion is the same cell argument: a source outside $\mathrm{Th}_{I^\vee}(\xi^-)$ has label $\sigma\in I$, and every destination label is $\le\sigma$.
\end{proof}

The source--target dictionary above recovers the fat-thickening properness theorem~\cite[Theorems~6.8 and~6.13]{KLP}; their cocompactness theorem additionally assumes CEA and balance~\cite[Theorem~7.12]{KLP}. The quantifier in \eqref{eq:absorption} is universal: every destination of a source outside the wall must lie in it. For the cyclic action generated by $g(z)=4z$ on $\C\PP^1$, every source has a destination in $L=\{\infty\}$ along positive powers (at $0$, take $z_m=2^{-m}$). Yet the constant source $0$ also returns to $0$, so $\C\PP^1\setminus L$ is not a domain.

\begin{example}\label{ex:sl3-table}
We make the cell calculation explicit for the full flag manifold
of $\PSL(3,\mathbb K)$, where $\mathbb K=\R$ or $\C$.
Let $X=\mathrm{Fl}(\K^3)=\{(q,\ell):q\in\ell\}$ and let $(g_m)$ be regular, with attracting flag $F^+=(p,L)$ and repelling flag $F^-=(p^-,L^-)$. The cells are indexed by $W=S_3$, with the Bruhat order, graded by length, and Theorem~\ref{thm:S} reads as follows; here $v$ is the relative position of $k'x$ in $B^-vB$, and the destination is $kX_v$. All dimensions in this example are taken over $\mathbb K$.
\begin{center}\footnotesize
\begin{tabular}{@{}ccp{.30\linewidth}p{.44\linewidth}@{}}
\toprule
$\ell(v)$ & $\dim_{\K}$ source & source condition (relative to $F^-$)
 & destination (relative to $F^+$)\\\midrule
$0$ & $3$ & $x$ opposite to $F^-$ & the point $F^+$\\
$1$ & $2$ & $q\in L^-$, $p^-\notin\ell$
 & $\{(q',L)\in X\}$, a line\\
$1$ & $2$ & $p^-\in\ell$, $q\notin L^-$
 & $\{(p,\ell')\in X\}$, a line\\
$2$ & $1$ & $\ell=L^-$, $q\ne p^-$
 & $\{(q',\ell')\in X:q'\in L\}$, a surface\\
$2$ & $1$ & $q=p^-$, $\ell\ne L^-$
 & $\{(q',\ell')\in X:p\in\ell'\}$, a surface\\
$3$ & $0$ & $x=F^-$ & all of $X$\\\bottomrule
\end{tabular}
\end{center}
Membership in $X$ includes the incidence condition $q\in\ell$. The alternatives in the length-one and length-two rows are exchanged by the duality of $\K\PP^2$, and in each row the dimensions add up to $3$.

The cell graph is the Bruhat order of $S_3$. Its down-sets $\mathrm{Th}$ with $\mathrm{Th}\cup w_0\mathrm{Th}=W$ are the fat thickenings: $\{e,s_1,s_2\}$, the two sets obtained by adding $s_1s_2$ or $s_2s_1$, the set $W\setminus\{w_0\}$, and $W$. Only $\{e,s_1,s_2\}$ is balanced, since $w_0\{e,s_1,s_2\}=\{s_1s_2,s_2s_1,w_0\}$.

For a sequence of type $\{\alpha_1\}$ (blocks $(1,2)$, attracting point $p$, repelling line $L^-$), the cells form a chain of dimensions $3,2,1$. With the source and target frames fixed, the bounded Levi datum induces a projectivity $\beta$ from $L^-$ to the pencil of lines through $p$. If $q\notin L^-$, the flag $(q,\ell)$ has the single destination $(p,\beta(\ell\cap L^-))$. If $q\in L^-$ and $\ell\ne L^-$, put $M=\beta(q)$; the destination is $\{(q',M):q'\in M\}$. If $\ell=L^-$, put again $M=\beta(q)$; the destination is the surface $\{(q',\ell')\in X:q'\in M\}$. The Levi datum moves these destinations inside their cells and is not determined by the attracting and repelling flags.
\end{example}

\begin{example}\label{ex:hitchin-sl3}
A Hitchin image $\Gamma<\mathrm{SL}(3,\R)$ of a closed surface group preserves a strictly convex domain $\Omega_c\subset\R\PP^2$, with boundary $C$ and flag limit set $\Lambda=\{(b,\ell_b):b\in C\}$, where $\ell_b$ is the tangent at $b$.
The thickening $W\setminus\{w_0\}$ leaves only flags $(p,\ell)$ with $\ell\cap C=\varnothing$ and $p$ on no $\ell_b$; this domain is empty, since every point of a line missing $\overline{\Omega_c}$ lies on a tangent.
The balanced thickening $\{e,s_1,s_2\}$ instead gives the dense domain
$\{(p,\ell):p\notin C,\ \ell\text{ not tangent to }C\}$, on which the action is cocompact~\cite[Theorem~7.12]{KLP}.
\end{example}

\section{Coalescence and recurrence criteria}\label{sec:recurrence}

\subsection{The diagonal of a coalescing datum}

To determine which wall points are recurrent, we separate the relative-position condition from the Levi condition in Theorem~\ref{thm:dest}. Use its notation for an extracted datum $\eta=(\Theta',k,k',b)$, with retraction $\rho_{\Theta'}$ along $N^-_{\Theta'}$. This retraction is $N^-_{\Theta'}$-invariant and $L_{\Theta'}$-equivariant. If $\Theta\subseteq\Theta'$, then
$P_{\Theta'}\subseteq P_\Theta$, $L_{\Theta'}\subseteq L_\Theta$, and
$N_{\Theta'}=N_\Theta\rtimes N^L_{\Theta'}$, where $N^L_{\Theta'}=N_{\Theta'}\cap L_\Theta$.

\medskip
From now on $\Theta$ is self-opposed, $\iota\Theta=\Theta$. Then $\dot w_0$ normalises $L_\Theta$ and
$\dot w_0 N_\Theta\dot w_0^{-1}=N_\Theta^-$. The relative positions of $y\in X=G/Q$ with respect to the base flag $P_\Theta$ are
the $P_\Theta$-orbits $C_\tau=P_\Theta\dot\tau Q/Q$, $\tau\in\mathcal C_\Theta=W_\Theta\backslash W/W_\theta$. Each
cell fibres over its Levi flag manifold
\[
Y_\tau=L_\Theta\dot\tau Q/Q ,
\]
with fibres the $N_\Theta$-orbits; $\rho^+:C_\tau\to Y_\tau$ denotes this retraction. In the notation of Theorem~\ref{thm:S}, $C_\tau=O^+_\tau$ and $\rho^+$ is the target-side retraction defined there. Left multiplication by $w_0$
acts on $\mathcal C_\Theta$ and reverses the closure order. Every $\Theta$-flag $\xi=kP_\Theta$, $k\in K$, has a
$K$-opposite flag $k\dot w_0 P_\Theta$. It does not depend on the choice of $k$, since two choices differ by
$K_\Theta=K\cap L_\Theta$, which $\dot w_0$ normalises. So the cells $C_\tau(\xi)=kC_\tau$, their Levi manifolds and their
retractions are canonically attached to $\xi$.

\begin{definition}
The datum $\eta=(\Theta',k,k',b)$, with $\Theta'\supseteq\Theta$, is \emph{coalescing at $\Theta$} if its attracting and
repelling flags have the same projection $\xi=kP_\Theta=k'^{-1}\dot w_0 P_\Theta$ to $G/P_\Theta$. Its \emph{effective Levi
element} is $e_\eta=bk'k$, and its \emph{Levi destination} on $Y_\sigma$ is
\[
\D^L_\eta(q)=\overline{N^L_{\Theta'}\,\rho_{\Theta'}(q)}\subseteq Y_\sigma ,\qquad q\in Y_\sigma .
\]
\end{definition}

If $\eta$ is coalescing at $\Theta$, then $k'k\in\dot w_0 P_\Theta\cap K=\dot w_0 K_\Theta$, so
$e_\eta\in L_{\Theta'}\dot w_0 K_\Theta\subseteq\dot w_0 L_\Theta$. Thus $e_\eta$ maps $Y_\tau$ onto $Y_{w_0\tau}$. The
Levi destination is well defined: $Y_\sigma$ is closed and invariant under $A_{\Theta'}\subseteq L_\Theta$, so
$\rho_{\Theta'}(Y_\sigma)\subseteq Y_\sigma$, and $N^L_{\Theta'}\subseteq L_\Theta$. It has the shape of the destination formula
for the reductive group $L_\Theta$ acting on $Y_\sigma$, with data of type $\Theta'\setminus\Theta$. If $\Theta'=\Theta$, then
$N^L_\Theta$ is trivial, $\rho_\Theta$ is the identity on $Y_\sigma$, and $\D^L_\eta(q)=\{q\}$.

\begin{definition}
For $\tau\in\mathcal C_\Theta$, the \emph{Levi incidence} of $\tau$ is
$\mathcal I_\tau=\{(p,q)\in Y_\tau\times Y_{w_0\tau}:p\in\overline{N_\Theta\,q}\}$. The \emph{canonical thickening} is
$\mathcal T_\star=\{\tau\in\mathcal C_\Theta:\tau\le w_0\tau\}$, and $\mathcal T_\star(\xi)=\bigcup_{\tau\in\mathcal T_\star}C_\tau(\xi)$.
\end{definition}

Since $\overline{N_\Theta q}\subseteq\overline{C_{w_0\tau}}$ for $q\in Y_{w_0\tau}$, we have $\mathcal I_\tau=\varnothing$ unless $\tau\le w_0\tau$.
If $\tau\le\sigma\le w_0\sigma$ then $w_0\sigma\le w_0\tau$, so $\tau\le w_0\tau$. Hence $\mathcal T_\star$ is a down-set and
$\mathcal T_\star(\xi)$ is closed.

\begin{theorem}\label{thm:diagonal}
Let $\eta=(\Theta',k,k',b)$ be coalescing at the self-opposed type $\Theta$, with coalescence flag $\xi=kP_\Theta$ and
effective Levi element $e=e_\eta$. Let $x=ky$ with $y\in C_\tau$ and $p=\rho^+(y)\in Y_\tau$. Then
\[
\D_\eta(x)=k\,\overline{N_\Theta\,\D^L_\eta(ep)},
\qquad\text{and}\qquad
x\in\D_\eta(x)\iff p\in\overline{N_\Theta\,\D^L_\eta(ep)} .
\]
Consequently:
\begin{enumerate}[label=\textup{(\alph*)},nosep]
\item if $x\in\D_\eta(x)$, then $\tau\in\mathcal T_\star$, so the diagonal of $\D_\eta$ is contained in $\mathcal T_\star(\xi)$;
\item if $(p,q)\in\mathcal I_\tau$ for some $q\in\D^L_\eta(ep)$, then $x\in\D_\eta(x)$;
\item if $\Theta'=\Theta$, then $x\in\D_\eta(x)$ if and only if $(p,ep)\in\mathcal I_\tau$.
\end{enumerate}
\end{theorem}

\begin{proof}
Write $y=np$ with $n\in N_\Theta$. We first show that the destination depends only on the Levi coordinate $p$ of $y$.
Since $k'=b^{-1}ek^{-1}$, we have $bk'x=ey$; Theorem~\ref{thm:dest} and the $L_{\Theta'}$-equivariance of $\rho_{\Theta'}$ give
$\D_\eta(x)=k\,\overline{N_{\Theta'}\rho_{\Theta'}(ey)}$. Since
$e\in\dot w_0 L_\Theta$,
\[
ey=(ene^{-1})\,(ep),\qquad ene^{-1}\in\dot w_0 N_\Theta\dot w_0^{-1}=N^-_\Theta\subseteq N^-_{\Theta'} ,
\]
and the left $N^-_{\Theta'}$-invariance of the retraction gives $\rho_{\Theta'}(ey)=\rho_{\Theta'}(ep)=:q'\in Y_{w_0\tau}$. Now
$N_{\Theta'}q'=N_\Theta N^L_{\Theta'}q'$. Because each element of $N_\Theta$ acts as a homeomorphism,
$\overline{N_\Theta N^L_{\Theta'}q'}=\overline{N_\Theta\,\overline{N^L_{\Theta'}q'}}=\overline{N_\Theta\D^L_\eta(ep)}$, which is the first
formula. This closed set is $N_\Theta$-invariant ($L_\Theta$ normalises $N_\Theta$), and $y=np$, which gives the
criterion. For (a): $\overline{N_\Theta\D^L_\eta(ep)}\subseteq\overline{C_{w_0\tau}}$ meets $C_\tau$ only if $\tau\le w_0\tau$.
Part (b) is immediate, and (c) is the case $\D^L_\eta(ep)=\{ep\}$.
\end{proof}

The order obstruction is independent of the Levi data. The remaining condition combines the Schubert incidence $\mathcal I_\tau$ with the Levi destination produced by the additional gaps. Full incidence makes (b) automatic; for self-opposite $\tau$, the incidence contains the diagonal and in the cases below equals it.

\subsection{Regular groups: the Bruhat-pair formula}

If $\Theta=\Delta$, then $L_\Delta=MA$ fixes every Weyl point, each $Y_\tau$ is a point, and the Levi condition disappears.
For $\xi=gB$ and $x\in X$ write $\operatorname{pos}(\xi,x)\in W/W_\theta$ for the coset with $g^{-1}x\in B\dot\tau Q$.

The necessary Bruhat inequality is due to Kapovich--Leeb--Porti~\cite[Proposition~6.2]{KLP}; Stecker proves the realisation of complementary positions for regular sequences~\cite[Lemmas~3.6--3.8]{SteckerBalanced}. We use the conventions of the destination formula.

Call $\Gamma$ \emph{regular} if every sequence of distinct elements has a subsequence of type $\Delta$. Let $\Lambda\subseteq G/B$ be the set of its attracting flags and $\mathcal P_\Gamma\subseteq\Lambda\times\Lambda$ the set of attracting--repelling pairs. Inversion puts the second endpoint in $\Lambda$; both sets are closed by diagonal extraction.

\begin{lemma}\label{lem:pairs}
$\mathcal P_\Gamma$ is invariant under $\Gamma\times\Gamma$, acting by $(\gamma,\gamma')\cdot(\xi^+,\xi^-)=(\gamma\xi^+,\gamma'\xi^-)$.
If $\Lambda$ is minimal, then $\mathcal P_\Gamma=\Lambda\times\Lambda$.
\end{lemma}
\begin{proof}
On $G/B$, Theorem~\ref{thm:dest} identifies $\xi^+$ as the singleton destination on the open cell opposite to $\xi^-$, and recovers $\xi^-$ from the transposed relation. The identity
$\D_{(\gamma g_m\gamma'^{-1})}=(\gamma'\times\gamma)\D_{(g_m)}$
therefore proves independent invariance of the two flags.
Every product orbit $\Gamma u\times\Gamma v$ is dense when $\Lambda$ is minimal, so the closed, nonempty invariant set $\mathcal P_\Gamma$ is the whole product. For Zariski-dense groups, minimality is Benoist's theorem recalled after Theorem~\ref{thm:Quint}.
\end{proof}

\begin{theorem}\label{thm:pair}
On every flag manifold $X=G/Q$:
\begin{enumerate}[label=\textup{(\roman*)},nosep]
\item A datum $\eta$ of type $\Delta$ with attracting and repelling flags $\xi^\pm$ has
$\D_\eta(x)=\{y:\operatorname{pos}(\xi^+,y)\le w_0\operatorname{pos}(\xi^-,x)\}$.
On $G/B$, the relation recovers $\xi^+$ from its singleton fibres and $\xi^-$ from its transpose.
\item If $\Gamma<G$ is discrete and regular, then
\[
\D_\Gamma=\{(x,y):\exists(\xi^+,\xi^-)\in\mathcal P_\Gamma,\ \operatorname{pos}(\xi^+,y)\le w_0\operatorname{pos}(\xi^-,x)\}.
\]
If $\Lambda$ is minimal and $P(x)=\{\operatorname{pos}(\xi,x):\xi\in\Lambda\}$, this is equivalent to $u\le w_0v$ for some $u\in P(y)$ and $v\in P(x)$.
\end{enumerate}
No distinctness or transversality of $\xi^+$ and $\xi^-$ is required.
\end{theorem}
\begin{proof}
If $k'x\in B^-\dot\sigma Q$, then $MA$ fixes $\dot\sigma Q$, so Theorem~\ref{thm:dest} gives
$\D_\eta(x)=k\overline{B\dot\sigma Q}$.
Since $\operatorname{pos}(\xi^-,x)=w_0\sigma$, this proves (i), including recovery of the flags on $G/B$.
Taking the union over extracted sequences as in Proposition~\ref{prop:sequence-union} proves (ii), and Lemma~\ref{lem:pairs} gives its last assertion.
\end{proof}

\begin{proposition}\label{prop:unique-fat}
The canonical down-set $\mathcal T_\star$ is contained in every fat down-set of $W/W_\theta$.
It is fat if and only if every position is comparable with its opposite, equivalently if and only if there is a unique minimal fat down-set.
In this case it is balanced exactly when no position is self-opposite.
If $\tau$ and $w_0\tau$ are incomparable, there are distinct minimal fat down-sets $I_1\not\ni w_0\tau$ and $I_2\not\ni\tau$.
\end{proposition}
\begin{proof}
Fatness puts either $\tau$ or $w_0\tau$ in a down-set; if $\tau\le w_0\tau$, either choice puts $\tau$ there.
Moreover $\mathcal T_\star\cup w_0\mathcal T_\star$ consists of the positions comparable with their opposite, and the intersection consists of the self-opposite positions.

A down-set is fat exactly when its complement $J$ is an up-set with $J\cap w_0J=\varnothing$.
For incomparable $\tau,w_0\tau$, the up-set $J_0=\{\sigma:\sigma\ge w_0\tau\}$ satisfies
$J_0\cap w_0J_0=\{\sigma:w_0\tau\le\sigma\le\tau\}=\varnothing$.
Extending it maximally gives $I_1$; exchanging the two positions gives $I_2$, distinct because fatness forces $\tau\in I_1$.
Conversely, when $\mathcal T_\star$ is fat it is the smallest fat down-set by the first assertion.
\end{proof}

\begin{theorem}\label{thm:regular-recurrence}
Let $\Gamma<G$ be discrete and regular, with minimal limit set $\Lambda\subseteq G/B$ (for instance, Zariski dense), and
let $X=G/Q$. If $\mathcal T_\star\subseteq W/W_\theta$ is fat, then
\[
\Rec_\Gamma=\mathcal T_\star(\Lambda)=\{x:P(x)\cap\mathcal T_\star\ne\varnothing\},
\]
and $X\setminus\mathcal T_\star(\Lambda)$ is the largest domain of proper discontinuity.
\end{theorem}

\begin{proof}
Let $(x,y)\in\D_\Gamma$, witnessed by $u=\operatorname{pos}(\xi^+,y)\le w_0v$ with $v=\operatorname{pos}(\xi^-,x)$. If $v\notin\mathcal T_\star$, then
$w_0v\in\mathcal T_\star$ by fatness, so $u\in\mathcal T_\star$. Thus $\mathcal T_\star(\Lambda)$, which is closed and invariant, is a vertex cover, and taking
$y=x$ gives $\Rec_\Gamma\subseteq\mathcal T_\star(\Lambda)$. Conversely, if $\operatorname{pos}(\xi,x)=\tau\in\mathcal T_\star$, the pair $(\xi,\xi)$ lies in
$\mathcal P_\Gamma$ by Lemma~\ref{lem:pairs}, and $\tau\le w_0\tau$ shows $x\in\D_\Gamma(x)$. The relational criterion
concludes.
\end{proof}

Here recurrence follows directly from the coalescing pairs $(\xi,\xi)$ and the destination formula. No hyperbolicity, Anosov property, irreducibility or return construction is required.

\begin{example}
In rank one $\mathcal T_\star=\{e\}$, which is Proposition~\ref{prop:D-rank1}. On $\PP(\K^d)$ the positions
$i=\min\{j:x\in\xi^j\}$ are totally ordered, $w_0i=d+1-i$, and $\mathcal T_\star=\{i\le\lceil d/2\rceil\}$. Put $r(x)=\min P(x)$.
Theorem~\ref{thm:pair} gives $\D_\Gamma(x)=U_{d+1-r(x)}$, where $U_j=\bigcup_{\xi\in\Lambda}\PP(\xi^j)$, and Theorem~\ref{thm:regular-recurrence} gives $\Rec_\Gamma=U_{\lceil d/2\rceil}$. This is the form of Theorem~\ref{thm:cusped} for every regular group with minimal limit set.
\end{example}

\begin{example}
For $\Delta$-Anosov representations of non-elementary word-hyperbolic groups, Stecker classifies maximal domains by minimal fat ideals~\cite[Corollary~3.16]{SteckerBalanced}, covering the maximality in Example~\ref{ex:hitchin-sl3}. His balanced-ideal counts are one on $\operatorname{Fl}(\R^3)$, ten on $\operatorname{Fl}(\R^4)$, two on $\operatorname{Gr}(3,6)$, and on $\mathbb K\PP^{n-1}$ one for even $n$ and none for odd $n$, $\mathbb K=\R,\C$~\cite[\S\S4.1--4.2 and Appendix~A.1]{SteckerBalanced}. Proposition~\ref{prop:unique-fat} distinguishes uniqueness of a minimal fat ideal from balance when self-opposite positions occur.
\end{example}

\begin{example}
On $\operatorname{Gr}(2,n)$ the order on positions $\{s_1<s_2\}$ is componentwise and $w_0S=\{n+1-s_2,n+1-s_1\}$. Thus $S\le w_0S$ if and
only if $s_1+s_2\le n+1$, and the order is total in this respect: $\mathcal T_\star$ is fat for every $\operatorname{Gr}(2,n)$.

For a non-trivial full flag manifold, the parity obstruction is expressed in Coxeter length. Put $N=\ell(w_0)>0$. If $N$ is even, a prefix of length $N/2$ of a reduced word for $w_0$ gives an element $w$ with $\ell(w)=N/2$. The distinct elements $w$ and $w_0w$ have the same length, since $\ell(w_0w)=N-\ell(w)$, and hence are incomparable in the Bruhat order. Thus neither lies in $\mathcal T_\star=\{v:v\le w_0v\}$, so $\mathcal T_\star$ is not fat. Consequently, fatness on $G/B$ requires $\ell(w_0)$ to be odd. This argument is independent of restricted root multiplicities and excludes $\operatorname{Fl}(\K^4)$, $\operatorname{Fl}(\K^5)$ and the full flags of $\mathrm{Sp}(4)$ and $G_2$.
\end{example}

\begin{proposition}\label{prop:diagonal-triples}
Let $\Gamma_0$ be a uniform lattice in $G_0=\PSL(2,\C)$, and let
$\Gamma=\{(\gamma,\gamma,\gamma):\gamma\in\Gamma_0\}$
act on the full flag manifold $X=(\C\PP^1)^3$ of $G=G_0^3$.
Then $\Gamma$ is regular with minimal flag limit set. Its largest domain of
proper discontinuity consists of the ordered triples of distinct points,
although $\mathcal T_\star$ is not fat.
\end{proposition}

\begin{proof}
Cocompactness in hyperbolic three-space implies that every point of a
geodesic ray lies at uniformly bounded distance from a fixed $\Gamma_0$-orbit. Thus its classical
limit set is the whole sphere $\C\PP^1$, on which $\Gamma_0$ acts minimally.
For any sequence of distinct elements of $\Gamma_0$, its rank-one Cartan
gap tends to infinity. The three gaps in $G$ are equal, so $\Gamma$ is
regular and its flag limit set is the diagonal copy of $\C\PP^1$ in $X$.

Put $\Omega=\{(x_1,x_2,x_3):x_i\ne x_j\text{ for }i\ne j\}$.
The map $g\mapsto(g0,g1,g\infty)$ identifies $G_0$ homeomorphically with
$\Omega$, since a M\"obius transformation is determined continuously by
the images of three distinct points. Under this identification $\Gamma$
acts by left translation by the discrete subgroup $\Gamma_0$, hence properly.

We show that every point of $X\setminus\Omega$ is recurrent.
Lemma~\ref{lem:pairs} gives, for any $p\ne q$, distinct $\gamma_m\in\Gamma_0$ with attractor $p$ and repeller $q$.
Rank-one convergence gives $\gamma_m p\to p$ and
$\gamma_m^{-1}q\to q$. Thus
$(p,p,\gamma_m^{-1}q)\to(p,p,q)$ and its images tend to $(p,p,q)$.
Permuting the coordinates covers every triple with exactly two distinct
values. For $(p,p,p)$ choose any $q\ne p$ and use the constant source.
A domain of proper discontinuity contains no recurrent point, so every
such domain is contained in $\Omega$.

The Weyl order on $X$ is inclusion in the subsets of $\{1,2,3\}$.
The position relative to $(p,p,p)$ is $\{i:x_i\ne p\}$, and opposition
takes a subset to its complement. Consequently $\mathcal T_\star$
contains only the empty subset. Its union with its opposite misses all
singletons and all two-element subsets.
\end{proof}

The geometric thickenings in this example explain why counting minimal
fat ideals alone cannot decide whether a largest domain exists. The balanced
ideal $\{S:|S|\le1\}$ thickens the diagonal limit set to the triples with
a repeated coordinate. For each $i$, the distinct balanced ideal
$\{S:i\notin S\}$ thickens that same limit set to all of $X$, since one
may take $p=x_i$. Thus distinct balanced ideals can give identical
thickenings, and their domains can be contained in another thickening
domain.

\begin{remark}\label{rem:fat-question}
The diagonal lattice lies in a proper algebraic subgroup of a product.
This leaves the following question. If $G$ is simple and $\Gamma<G$ is
discrete, regular and Zariski dense, does the existence of a nonempty
largest domain on $G/Q$ imply that $\mathcal T_\star$ is fat?
\end{remark}

\subsection{Singular self-opposed types}

Let $\Theta=\iota\Theta$, retaining the definition of $\Theta$-convergence and the attracting set $\Lambda_\Theta$ from Section~\ref{sec:domains}. For singular types the Levi condition may require a different coalescing datum at each point.

\begin{definition}\label{def:recurrence-seed}
A flag $\xi\in\Lambda_\Theta$ is a \emph{recurrence seed} if for every $\tau\in\mathcal T_\star$ and every $x\in Y_\tau(\xi)$ there is a
datum $\eta=(\Theta',k,k',b)$ of an extracted sequence of distinct elements of $\Gamma$, coalescing at $\Theta$ with coalescence flag $\xi$,
such that, in the standard frame $p=k^{-1}x\in Y_\tau$,
\[
p\in\overline{N_\Theta\D^L_\eta(e_\eta p)}.
\]
Equivalently, $x\in\D_\eta(x)$; this makes the condition independent of the chosen compact frame.
\end{definition}

\begin{theorem}\label{thm:singular-recurrence}
Let $\Gamma<G$ be discrete and $\Theta$-convergent, with $\Theta$ self-opposed, and let $X=G/Q$. Assume that
\begin{enumerate}[label=\textup{(\roman*)},nosep]
\item $\Lambda_\Theta$ is minimal;
\item $\mathcal T_\star\subseteq\mathcal C_\Theta$ is fat;
\item $\Lambda_\Theta$ contains a recurrence seed.
\end{enumerate}
Then $\Rec_\Gamma=\mathcal T_\star(\Lambda_\Theta)$, and its complement is the largest domain of proper discontinuity.
\end{theorem}

\begin{proof}
By (ii) and Theorem~\ref{thm:KLPfull}, which is stated for $\Theta$-convergent groups, $\mathcal T_\star(\Lambda_\Theta)$ is a closed
invariant vertex cover, so $\Rec_\Gamma\subseteq\mathcal T_\star(\Lambda_\Theta)$.

If $\xi$ is a seed, Theorem~\ref{thm:diagonal} gives
$\mathcal T_\star(\xi)\subseteq\Rec_\Gamma$, because every point of $C_\tau(\xi)$ lies over some $p\in Y_\tau(\xi)$.

For
$\zeta\in\Lambda_\Theta=\overline{\Gamma\xi}$, choose $\gamma_m\xi\to\zeta$. Put $\xi_0=P_\Theta$. A local section of $G\to G/P_\Theta$ gives $h_m\to h$ in $G$ with $\gamma_m\xi=h_m\xi_0$ and
$\zeta=h\xi_0$. Every point $hz$ of $\mathcal T_\star(\zeta)=h\mathcal T_\star(\xi_0)$ is the limit of $h_mz\in\gamma_m\mathcal T_\star(\xi)\subseteq\Rec_\Gamma$.
The recurrent set is closed, and the relational criterion concludes.
\end{proof}

\section{Largest domains on Grassmannians}\label{grlevi:sec:grassmann-levi}
Let $U=\C^2$, let $V$ be Hermitian of dimension $r\ge2$, and let
$\Gamma_0<\operatorname{SL}(U)$ be discrete and non-elementary, with
$-I\notin\Gamma_0$ and classical limit set $\Lambda\subseteq\PP(U)$.
For a unitary representation $\vartheta:\Gamma_0\to U(V)$ put
\[
 H=\overline{\vartheta(\Gamma_0)},\qquad W=U\otimes V,\qquad
 \Gamma_\vartheta=\{[\gamma\otimes\vartheta(\gamma)]:\gamma\in\Gamma_0\}.
\]
For $p\in\PP(U)$ put $E_p=p\otimes V$.
The matrix action of Section~\ref{sec:levi-recovery} is recovered by
$V=(\C^r)^*$ and $\vartheta=\chi^\vee$, the contragredient of the
right twist $\chi$. We work intrinsically on this compact module.
Theorem~\ref{thm:marked-product} gives faithfulness and discreteness:
the singular values are $s,s^{-1}$, each repeated $r$ times.
Thus every escaping sequence, after extraction, has type $(r,r)$,
with attracting and repelling planes $E_p,E_q$, $p,q\in\Lambda$,
independently of the twist. In root notation
$\Theta=\{\alpha_r\}$ and $b=1$; the carrier family
$\Lambda_\Theta=\{E_p:p\in\Lambda\}$ is minimal.
Lemma~\ref{lem:saturate} realises every joint limit $(p,q,h)$.
The sufficient criteria below allow $\Lambda=\PP(U)$;
$\Lambda\subsetneq\PP(U)$ is imposed explicitly for converses and nonemptiness.
All Grassmannians are the full complex Grassmannians.

\subsection{The incidence relation in every degree}
For $1\le k\le2r-1$ let $Y_k=\operatorname{Gr}(k,W)$.
For $q\in\PP(U)$ choose $\ell_q\in U^*$ with kernel $q$ and set
\begin{equation}\label{grlevi:eq:KP}
 K_q(S)=\{v\in V:q\otimes v\in S\},\qquad
 P_q(S)=(\ell_q\otimes\mathrm{id}_V)(S),\qquad
 \pi_q=\ell_q\otimes\mathrm{id}_V.
\end{equation}
These subspaces are independent of the scale of $\ell_q$ and satisfy
$\dim K_q(S)+\dim P_q(S)=k$.

\begin{theorem}\label{grlevi:thm:allk}
For every $\Lambda\subseteq\PP(U)$, the group relation on $Y_k$ is
\begin{equation}\label{grlevi:eq:allk}
\D_{\Gamma_\vartheta}=\bigcup_{p,q\in\Lambda,\ h\in H}R^{(k)}_{p,q,h},\qquad
R^{(k)}_{p,q,h}(S)=\{T\in Y_k:
p\otimes hP_q(S)\subseteq T,\ \pi_p(T)\subseteq hK_q(S)\}.
\end{equation}
For $d=\dim K_q(S)$, necessarily
$\max\{0,k-r\}\le d\le\min\{k,r\}$;
the fibre is isomorphic to
$\operatorname{Gr}(d,r-k+2d)$, of complex dimension $d(r-k+d)$.
Every point of each fibre is realised along the full extracted sequence.
For $\Lambda\subsetneq\PP(U)$ and $k\le r$, the complement of
$\Sigma^{(k)}_\Lambda=\{S:S\cap E_q\ne0\text{ for some }q\in\Lambda\}$
is a nonempty properness domain for every unitary twist.
\end{theorem}
\begin{proof}
We first compute the relation for two fixed blocks; joint saturation will then give the group relation. Let $a_n$ act on
$V^+\oplus V^-$ as $a_n(u,v)=(s_nu,s_n^{-1}v)$ with
$s_n\to\infty$. For a $k$-plane $S$ write
$K=S\cap V^-$, $P=\operatorname{pr}_+S$, and $d=\dim K$.
We claim
\begin{equation}\label{grlevi:eq:twoblock-proof}
\D_{(a_n)}(S)=\{T\in\operatorname{Gr}(k,V^+\oplus V^-):
P\subseteq T\subseteq V^+\oplus K\}.
\end{equation}
To prove the first inclusion, take $S_n\to S$, $a_nS_n\to T$.
Lift a basis of $P$ to $k-d$ vectors $x_i\in S$ and choose
$x_{i,n}\in S_n$ with $x_{i,n}\to x_i$. The vectors
$s_n^{-1}a_nx_{i,n}\in a_nS_n$ converge to the fast projections
of $x_i$, which form a basis of $P$; hence $P\subseteq T$.
Now let $z_n=(u_n,v_n)\in a_nS_n$ converge to $z=(u,v)\in T$.
If $v\ne0$, the normalised inverse vectors
\[
\frac{a_n^{-1}z_n}{\|a_n^{-1}z_n\|}
=\frac{(s_n^{-1}u_n,s_nv_n)}{\|(s_n^{-1}u_n,s_nv_n)\|}
\longrightarrow (0,v/\|v\|)
\]
belong to $S_n$, so their limit belongs to $S$ and $v\in K$.
The assertion is automatic if $v=0$. Consequently
$T\subseteq V^+\oplus K$.

Conversely, among the planes $T$ with $P\subseteq T\subseteq V^+\oplus K$,
those for which projection onto $K$ is surjective are open and dense.
They have the form
$T_F=P\oplus\{(Fv,v):v\in K\}$, with a linear map
$F:K\to V^+$ understood modulo maps into $P$.
Choose a linear section of $S\to P$, written $x(u)=(u,f(u))$.
Then
\[
S_n=\operatorname{span}\{(u,f(u)):u\in P\}
 +\operatorname{span}\{(s_n^{-2}Fv,v):v\in K\}
\longrightarrow S.
\]
Normalise the images of the first $k-d$ generators by $s_n^{-1}$
and the remaining $d$ by $s_n$. This yields
$a_nS_n\to T_F$. If $T$ belongs to the boundary of this open set, choose
linear maps $F_j:K\to V^+$ with $T_{F_j}\to T$.
Let $S_n^{(j)}$ be the above sources for $F_j$. Choose $N_j$ increasing
so that, for all $n\ge N_j$,
\[
 d_{\mathrm{Gr}}(S_n^{(j)},S)+d_{\mathrm{Gr}}(a_nS_n^{(j)},T_{F_j})<1/j
\]
for a compatible metric. Taking $j(n)\to\infty$ with $n\ge N_{j(n)}$
gives $S_n=S_n^{(j(n))}\to S$ and $a_nS_n\to T$ along the
\emph{whole} sequence, including every boundary stratum.
The cases $d=0$ (a singleton fibre) and $d=k$ ($P=0$) are included.
This proves~\eqref{grlevi:eq:twoblock-proof}.
Quotienting by $P$ identifies the fibre with
$\operatorname{Gr}(d,(V^+/P)\oplus K)$.

For a sequence of distinct elements of $\Gamma_\vartheta$, extract a
singular-value decomposition in the first factor and a convergent
subsequence in the compact unitary factor. The source slow plane and
target fast plane tend respectively to $E_q$ and $E_p$; the residual
map on both copies of $V$ tends to $h$. Transporting the two-block
statement gives the two incidence inclusions in~\eqref{grlevi:eq:allk}.
Conversely, Lemma~\ref{lem:saturate} realises every $(p,q,h)$ and the
whole-sequence construction gives all its incidence pairs. The
quotient has dimension $(r-(k-d))+d=r-k+2d$,
so the fibre has complex dimension $d(r-k+d)$.

For the final statement, the wall is closed by projection of the
compact incidence set over $\Lambda$ and is invariant. A source
outside it has $K_q(S)=0$ for every $q\in\Lambda$, so every
destination is a $k$-plane contained in some $E_p$, hence belongs to
the wall. The relational criterion (Theorem~\ref{thm:relational})
proves properness. For nonemptiness, choose $s\notin\Lambda$ and a $k$-plane
$D\subseteq V$; then $s\otimes D$ meets every $E_q$, $q\in\Lambda$,
trivially.
\end{proof}

\subsection{Compact incidence criteria and their converses}
Fix $p\in\Lambda$ and write $E=E_p$, $E'=E_{p^\perp}$.
The relative positions are
$\tau_j=\{S:\dim(S\cap E)=j\}$,
$\max(0,k-r)\le j\le\min(k,r)$, with $\tau_j\le\tau_i$ exactly
when $j\ge i$, and $w_0\tau_j=\tau_{k-j}$. Thus
\[
 \mathcal T_\star=\{\tau_j:j\ge k/2\},\qquad
 \Sigma^\star_k=\{S:\dim(S\cap E_p)\ge k/2
                         \text{ for some }p\in\Lambda\}.
\]
This ideal is fat, and balanced exactly when $k$ is odd.
The Levi manifolds are
$Y_{\tau_j}=\operatorname{Gr}(j,E)\times\operatorname{Gr}(k-j,E')$,
represented by $A\oplus B$, and $N_\Theta=\operatorname{Hom}(E',E)$.
Theorem~\ref{grlevi:thm:allk} gives
$\overline{N_\Theta(A'\oplus B')}=\{T:A'\subseteq T\subseteq E\oplus B'\}$;
hence
\[
 \mathcal I_{\tau_j}=\{((A,B),(A',B')):A'\subseteq A,\ B\subseteq B'\},
\]
which is the diagonal when $j=k/2$.
For a coalescing triple $(p,p,h)$ the tensor Cartan limits are
$k=(\lim k_m)\otimes h$ and $k'=(\lim k'_m)\otimes I$.
Since $b=1$, the effective element is $e_\eta=k'k$; in the splitting
$E\oplus E'$ it becomes $ke_\eta k^{-1}=w\otimes h$, where $w$
exchanges $p,p^\perp$. It sends the residual pair $(K,P)$ to
$(hP,hK)$, and every $h\in H$ occurs by Lemma~\ref{lem:saturate}.
Equivalently, setting $q=p$ and $T=S$ in~\eqref{grlevi:eq:allk}
gives the following diagonal criterion.

\begin{proposition}\label{prop:grk}
A plane $S$ with $\dim(S\cap E_p)=j$ lies in the diagonal of the datum $(p,p,h)$ if and only if
\[
hP_p(S)\subseteq K_p(S)\qquad\text{and}\qquad h^{-1}P_p(S)\subseteq K_p(S).
\]
Consequently $E_p$ is a recurrence seed if and only if $H$ satisfies, for every $j$ with $k/2\le j\le\min(k,r)$:
\[
(\mathrm R_{j,k-j})\qquad \forall K\in\operatorname{Gr}(j,V),\ P\in\operatorname{Gr}(k-j,V)\ \ \exists h\in H:\ hP+h^{-1}P\subseteq K .
\]
\end{proposition}

Here $(\mathrm R_{b,b})$ means pair reversal on
$\operatorname{Gr}(b,V)$; $(\mathrm R_{a,0})$ and $(\mathrm R_{r,b})$
are automatic. For lines we write \emph{projective pair reversibility} as
\begin{equation*}\label{grlevi:eq:swap}
 \forall L,M\in\PP(V)\quad\exists h\in H:
 \qquad hL=M\quad\text{and}\quad hM=L.\tag{*}
\end{equation*}

\begin{lemma}\label{grlevi:lem:swap-transitivity}
Let $H\leq U(r)$ be compact. Then $H$ has property~\eqref{grlevi:eq:swap}
if and only if it acts transitively on $\PP^{r-1}(\C)$.
\end{lemma}
\begin{proof}
The pair-interchange condition plainly implies transitivity. Conversely,
assume that $H$ is transitive. Its identity component $H^0$ is also
transitive: the quotient $\PP^{r-1}(\C)/H^0$ is connected, whereas the
finite group $H/H^0$ acts transitively on it; hence this quotient is a
single point. By Onishchik's classification of connected compact groups transitive on
complex projective space, the effective projective action of $H^0$ is
conjugate to the standard action of $PU(r)$ or, when $r=2n$, of $PSp(n)$
\cite[Theorem~6(b)]{Onishchik1963}. Pair interchange is invariant under
conjugacy, so it is enough to verify it for $SU(r)$ and $Sp(n)$.

For $SU(r)$, multiply unit representatives $x,y$ by unit complex phases
so that $\langle x,y\rangle$ is real. If the two lines coincide there is nothing to prove. Otherwise, the unitary Householder reflection
$R(z)=z-2\langle z,x-y\rangle(x-y)/\|x-y\|^2$
(up to the convention for the Hermitian product) exchanges $x$ and $y$; the reality of $\langle x,y\rangle$ gives $\|x-y\|^2=2(1-\langle x,y\rangle)$. Since $R$ has determinant $-1$, choose a unit complex scalar $\zeta$ with $\zeta^r=-1$; then $\zeta R\in SU(r)$ and induces the same transformation of projective lines.

For $Sp(n)$, identify $\C^{2n}$ with $\mathbb H^n$ and write the
quaternionic Hermitian product of unit representatives (with $jz=\bar z j$) as
$\langle x,y\rangle_{\mathbb H}=a+bj$, with $a,b\in\C$.
Choose $\alpha,\beta\in U(1)$ so that
$\overline\alpha\,\overline{(a+bj)}\,\beta=a+bj$.
If $a,b\ne0$ this is equivalent to
$\beta/\alpha=a/\bar a$ and $\alpha\beta=-1$.
Choose any $\alpha\in U(1)$ with
$\alpha^2=-\bar a/a$ and set $\beta=-\bar\alpha$;
then both equations hold. If $a=0$ and $b\ne0$, take
$\alpha=1$, $\beta=-1$; if $b=0$ and $a\ne0$, take
$\alpha=1$, $\beta=a/\bar a$; if $a=b=0$,
take $\alpha=\beta=1$.

Thus the ordered quaternionic Gram
matrices of $(x,y)$ and $(y\alpha,x\beta)$ agree. When the pair is quaternionically dependent, the resulting isometry is defined on their common quaternionic line; otherwise it is defined on their two-dimensional quaternionic span. In either case, quaternionic orthogonal extension gives an element of $Sp(n)$ exchanging the two complex projective lines.
\end{proof}

\begin{lemma}\label{lem:Rab}
\begin{enumerate}[label=\textup{(\alph*)},nosep]
\item If $H$ reverses pairs on $\operatorname{Gr}(b,V)$, then $(\mathrm R_{a,c})$ holds for $c\le b\le a$.
\item $SU(r)$ reverses pairs on every $\operatorname{Gr}(b,V)$.
\item For $r=2n\ge4$, $Sp(n)$ reverses pairs on $\PP(V)$ but not on $\operatorname{Gr}(2,V)$.
\end{enumerate}
\end{lemma}

\begin{proof}
(a) Let $K\in\operatorname{Gr}(a,V)$ and $P\in\operatorname{Gr}(c,V)$. Choose $P'\supseteq P$ and $K_0\subseteq K$ of dimension $b$, and $h\in H$ with $hP'=K_0$ and $hK_0=P'$. Then $hP\subseteq hP'=K_0\subseteq K$ and $h^{-1}P\subseteq h^{-1}P'=K_0\subseteq K$.

(b) $\operatorname{Gr}(b,V)$
is a compact symmetric space of $U(r)$. Two points are joined by a minimising geodesic with midpoint $M$, and the
geodesic symmetry $2\pi_M-I\in U(r)$ exchanges them. A unit scalar puts it in $SU(r)$ without changing its action on
subspaces.

(c) The first part is Lemma~\ref{grlevi:lem:swap-transitivity}. For the second, let $\omega$ be the complex
symplectic form preserved by $Sp(n)$. The number $|\omega(u,v)|$, for an orthonormal basis $(u,v)$ of a $2$-plane, is an
$Sp(n)$-invariant with values $0$ on isotropic planes and $1$ on planes of the form $\langle u,Ju\rangle$, and both kinds
exist for $n\ge2$.
\end{proof}

\begin{theorem}\label{thm:grk}
Let $1\le k\le2r-1$ and suppose $H$ satisfies $(\mathrm R_{j,k-j})$ for $k/2\le j\le\min(k,r)$. Then the recurrent set of
$\Gamma_\vartheta$ on $\operatorname{Gr}(k,W)$ is $\Sigma^\star_k$, and
\[
\Omega^\star_k=\{S\in\operatorname{Gr}(k,W):\ \dim(S\cap E_p)<k/2\ \text{for all }p\in\Lambda\}
\]
is the largest domain of proper discontinuity. This holds for every $k$ if $H\supseteq SU(r)$. For $H=Sp(n)$, $r=2n\ge4$,
it holds for $k\in\{1,2,3,2r-3,2r-2,2r-1\}$, and the seed condition fails for $k=4$.
\end{theorem}

\begin{theorem}\label{thm:grk-converse}
Assume $\Lambda\subsetneq\PP(U)$. The following also holds for
$r=1$, with the other common hypotheses unchanged.
For $1\le k\le2r-1$ put
\[
 \ell=\min(k,2r-k),\qquad a=\lceil\ell/2\rceil,\qquad b=\lfloor\ell/2\rfloor .
\]
The action of $\Gamma_\vartheta$ on $\operatorname{Gr}(k,U\otimes V)$ has a
largest domain of proper discontinuity if and only if $(\mathrm R_{a,b})$ holds,
that is, if for every $K\in\operatorname{Gr}(a,V)$ and $P\in\operatorname{Gr}(b,V)$
some $h\in H$ satisfies $hP+h^{-1}P\subseteq K$.
In this case the recurrent set is $\Sigma^\star_k$ and the largest domain is
the nonempty set $\Omega^\star_k$. If the condition fails, two non-recurrent
$k$-planes are related by $\D_{\Gamma_\vartheta}$.
For $\ell=2b$ the condition says precisely that $H$ exchanges every pair of
$b$-planes in $V$.
\end{theorem}

\begin{proof}[Proof of Theorems~\ref{thm:grk} and~\ref{thm:grk-converse}]
For $r=1$ the twist acts projectively by scalars, so
Proposition~\ref{prop:D-rank1} gives the assertion on $\PP(U)$;
$(\mathrm R_{1,0})$ is automatic. Assume $r\ge2$.
For any $k$ and any compact twist, if $S\notin\Sigma^\star_k$ then
$d_q=\dim K_q(S)<k/2$ for all $q\in\Lambda$.
Every destination in $R^{(k)}_{p,q,h}(S)$ contains
$p\otimes hP_q(S)$ of dimension $k-d_q>k/2$, and therefore belongs
to $\Sigma^\star_k$. Theorem~\ref{thm:relational} proves properness
of $\Omega^\star_k$, without requiring $\Lambda$ to be proper.
Under the conditions of Theorem~\ref{thm:grk}, every point of the wall
is recurrent by Proposition~\ref{prop:grk}; its complement is therefore
the largest domain. This also follows from
Theorem~\ref{thm:singular-recurrence}.

For the sharper criterion first take $k\le r$, so $k=a+b$.
If $(\mathrm R_{a,b})$ holds, given a $j$-plane $K$ and a
$(k-j)$-plane $P$ with $j\ge a$, choose an $a$-plane $K_0\subseteq K$
and a $b$-plane $P_0\supseteq P$. An $h$ with
$hP_0+h^{-1}P_0\subseteq K_0$ gives
$hP+h^{-1}P\subseteq K$, proving all the sufficient conditions.
If $\Lambda$ is proper, $s\otimes D\in\Omega^\star_k$ for any
$s\notin\Lambda$ and any $k$-plane $D\subseteq V$.

Conversely, suppose $(\mathrm R_{a,b})$ fails; then $b>0$.
The pairs $(K,P)$ satisfying its inclusions for some $h\in H$ form
a closed subset of $\operatorname{Gr}(a,V)\times\operatorname{Gr}(b,V)$,
by compactness of $H$ and closedness of incidence. Since $a+b\le r$,
the pairs with $K\cap P=0$ form a dense open subset: their complement
is a proper determinantal subvariety of the irreducible product.
Choose a transverse violating pair, $q\in\Lambda$, $s\notin\Lambda$,
and set
\[
 S=q\otimes K\oplus s\otimes P,\qquad
 T=q\otimes P\oplus s\otimes K.
\]
A tensor with nonzero components in both summands has rank two,
because both its left and right factors are independent. Thus $E_q$
is the only $E_{q'}$, $q'\in\Lambda$, meeting either plane nontrivially,
and
$K_q(S)=K$, $P_q(S)=P$, $K_q(T)=P$, $P_q(T)=K$.
If $S\in R^{(k)}_{p',q',h}(S)$ with $q'\ne q$, the second inclusion
of~\eqref{grlevi:eq:allk} forces $S\subseteq E_{p'}$, impossible.
Hence $q'=q$; since $b>0$, the first inclusion forces $p'=q$.
The two inclusions now give $hP\subseteq K$ and $P\subseteq hK$,
contradicting the choice of the pair. Thus $S$ is non-recurrent.
If $a>b$, then $T\notin\Sigma^\star_k$, so properness proved above
makes $T$ non-recurrent. If $a=b$, recurrence of $T$ would force an
element exchanging $P$ and $K$, by the same two inclusions in reversed
order. Finally $T\in R^{(k)}_{q,q,I}(S)$; the related non-recurrent
pair excludes a largest domain by Theorem~\ref{thm:relational}.

For $k>r$, annihilation gives an equivariant isomorphism
$\operatorname{Gr}(k,W)\to\operatorname{Gr}(2r-k,W^*)$ for the
contragredient action on $W^*=U^*\otimes V^*$.
An alternating form identifies $U^*$ complex linearly and
$\mathrm{SL}(U)$-equivariantly with $U$. The compact group becomes
$H^\vee=\{(h^{-1})^t:h\in H\}$, where $^t$ is the algebraic dual.
The conjugate-linear Riesz isomorphism $V\to V^*$ intertwines
$H,H^\vee$ and preserves complex subspaces and incidence, so
$(\mathrm R_{a,b})$ is equivalent for both groups.
Under these identifications $E_q^\circ=q\otimes V^*$, and
\[
 \dim(S^\circ\cap E_q^\circ)
   =2r-\dim(S+E_q)=r-k+\dim(S\cap E_q).
\]
Subtracting $(2r-k)/2$ gives $\dim(S\cap E_q)-k/2$.
Thus annihilation carries the canonical walls onto each other and
preserves the relation by carrying convergent source and image
sequences to their dual counterparts. The case $2r-k\le r$ proves
the remaining assertions of the converse.

For $SU(r)$ all sufficient conditions follow from
Lemma~\ref{lem:Rab}(a,b). For $Sp(n)$, degrees $k\le3$ require only
$(\mathrm R_{1,1})$ and $(\mathrm R_{2,1})$, which follow from
Lemma~\ref{lem:Rab}(a,c); duality gives $k\ge2r-3$.
For $k=4$, the seed condition $(\mathrm R_{2,2})$ fails by
Lemma~\ref{lem:Rab}(c).
\end{proof}

For $k\ge3$, the canonical wall $\Sigma^\star_k$ is strictly smaller
than the incidence wall $\Sigma^{(k)}_\Lambda$ of
Theorem~\ref{grlevi:thm:allk}, so under the sufficient hypotheses its
complement gives a larger domain. On $\operatorname{Gr}(3,6)$ the
singular type has a chain of relative positions, whereas a
$\Delta$-Anosov action has two balanced ideals in Stecker's classification.
For even $\ell$, the criterion is reversal of equal-dimensional
subspaces; for odd $\ell$, it requires
$hP+h^{-1}P\subseteq K$ for every $b$-plane $P$ and $(b+1)$-plane $K$.

\subsection{Two-planes and strongly irreducible examples}
Put $Y=Y_2$ and
\begin{equation}\label{grlevi:eq:wall}
 \Sigma_\Lambda=\{S\in Y:S\cap E_p\ne0\text{ for some }p\in\Lambda\},
 \qquad \Omega_0=Y\setminus\Sigma_\Lambda.
\end{equation}
These are $\Sigma^\star_2$ and $\Omega^\star_2$.

\begin{proposition}\label{grlevi:prop:incidence}
For $Y=Y_2$ formula~\eqref{grlevi:eq:allk} reads
\begin{equation}\label{grlevi:eq:fibre}
R_{p,q,h}(S)=\{T\in Y:p\otimes hP_q(S)\subseteq T,
                  \ \pi_p(T)\subseteq hK_q(S)\}.
\end{equation}
The fibres for $d=\dim K_q(S)=0,1,2$ are respectively a point,
$\PP^{r-1}$ and $\operatorname{Gr}(2,r+2)$.
\end{proposition}
\begin{proof}
This is Theorem~\ref{grlevi:thm:allk} with $k=2$.
\end{proof}

\begin{lemma}\label{grlevi:lem:proper}
The wall $\Sigma_\Lambda$ is closed and invariant, and $\Omega_0$ is a
properness domain for every $H$. If $\Lambda\subsetneq\PP(U)$,
then $\Omega_0$ is nonempty.
\end{lemma}
\begin{proof}
Apply the final assertion of Theorem~\ref{grlevi:thm:allk} with $k=2$;
its properness argument does not require $\Lambda$ to be proper.
\end{proof}

\begin{lemma}\label{grlevi:lem:recurrence}
Let $q\in\Lambda$, let $s\notin\Lambda$, and let $K,P\in\PP(V)$ be distinct lines. Put
\[
 S(q,s;K,P)=\operatorname{span}(q\otimes K,\ s\otimes P).
\]
Then $S(q,s;K,P)$ is recurrent if and only if some $h\in H$ interchanges $K$ and $P$.
\end{lemma}
\begin{proof}
The only decomposable projective lines of $S=S(q,s;K,P)$ have left
factors $q,s$, so $E_q$ is the only limiting carrier meeting $S$.
As in the proof of Theorem~\ref{thm:grk-converse}, recurrence in
\eqref{grlevi:eq:fibre} forces the source and target labels both to be $q$.
Its inclusions then say $hP=K$ and $hK=P$. Conversely those equalities
give $S\in R_{q,q,h}(S)$ by the same formula.
\end{proof}

\begin{theorem}\label{grlevi:thm:main}
Assume $\Lambda\subsetneq\PP(U)$. Then every escaping sequence in $\Gamma_\vartheta$ has, after extraction, the same singular Cartan block type $(r,r)$. The attracting and repelling $r$-planes are respectively $E_p$ and $E_q$ with $p,q\in\Lambda$, independently of $\vartheta$. The open set $\Omega_0$ is nonempty, invariant, and a domain of proper discontinuity for every unitary twist. Moreover, the following conditions are equivalent:
\begin{enumerate}[label=\textup{(\roman*)}]
\item $H$ satisfies projective pair reversibility~\eqref{grlevi:eq:swap};
\item the recurrent set of $\Gamma_\vartheta$ on $Y$ is $\Sigma_\Lambda$;
\item the action admits a \emph{largest} domain of proper discontinuity.
\end{enumerate}
When these conditions hold, that largest domain is $\Omega_0$. When they fail, there are two nonrecurrent planes related by $\D_{\Gamma_\vartheta}$, and hence at least two incomparable maximal properness domains. Equivalently, all three conditions hold exactly when $H$ is transitive on $\PP(V)$, by Lemma~\ref{grlevi:lem:swap-transitivity}.
\end{theorem}
\begin{proof}
The singular type and carriers were established in the common setup,
and Lemma~\ref{grlevi:lem:proper} gives the nonempty domain.
For $k=2$, the criterion of Theorem~\ref{thm:grk-converse} is exactly
$(\mathrm R_{1,1})$, or~\eqref{grlevi:eq:swap}; it gives (i)$\Leftrightarrow$(iii)
and identifies the recurrent set when these hold.
If reversibility fails for $K,P$, choose $p,q\in\Lambda$ and
$s,t\notin\Lambda$. The geometric witnesses are
\[
 S=\operatorname{span}(q\otimes K,s\otimes P),\qquad
 T=\operatorname{span}(p\otimes P,t\otimes K).
\]
Both are non-recurrent by Lemma~\ref{grlevi:lem:recurrence}, while
$T\in R_{p,q,I}(S)$ by~\eqref{grlevi:eq:fibre}. In particular the wall
contains a non-recurrent point, so (ii) fails as well.
Theorem~\ref{thm:relational} supplies incomparable maximal extensions
of properness domains containing these two points.
\end{proof}

The coalescing interpretation has
$\mathcal T_\star=\{\tau_1,\tau_2\}$ and wall $\Sigma_\Lambda$.
Its recurrence-seed condition is pair reversibility, by
Proposition~\ref{prop:grk}; thus Theorem~\ref{thm:singular-recurrence}
gives sufficiency, while the full relation supplies the obstruction above.

\begin{remark}
For a direct coalescing model, take independent $u,v\in V$ and
$n_t=\left(\begin{smallmatrix}1&t\\0&1\end{smallmatrix}\right)\otimes I_V$.
Then
\[
\begin{aligned}
 S_t&=\operatorname{span}(e_1\otimes u-t^{-1}e_2\otimes u,e_2\otimes v)
       \longrightarrow\operatorname{span}(e_1\otimes u,e_2\otimes v),\\
 n_tS_t&=\operatorname{span}(-t^{-1}e_2\otimes u,t e_1\otimes v+e_2\otimes v)
       \longrightarrow\operatorname{span}(e_2\otimes u,e_1\otimes v).
\end{aligned}
\]
The two limiting planes are distinct, related and non-recurrent for
$\langle n_1\rangle$. Joint saturation makes this model available to our
discrete actions without requiring an actual parabolic in $\Gamma_0$.
\end{remark}

\begin{corollary}\label{grlevi:cor:two-actions}
For every integer $r\ge3$ there are two faithful, discrete, strongly
irreducible representations of the same rank-two free Kleinian group into
$\PGL(2r,\C)$ such that every escaping sequence has singular Cartan
type $(r,r)$ and the attracting and repelling rank-$r$ carrier families
coincide, but their actions on $\operatorname{Gr}(2,\C^{2r})$ have
opposite largest-domain behavior. The compact residual closures can
be chosen as
\[
H_+=\operatorname{SU}(r),\qquad
H_-=\operatorname{Sym}^{r-1}(\operatorname{SU}(2))\subset\operatorname{SU}(r).
\]
The action with $H_+$ has the nonempty largest domain
$Y\setminus\Sigma_\Lambda$, while the action with $H_-$ has no largest
domain. On $\PP^{2r-1}$ both actions have the same largest domain,
the complement of $\bigcup_{p\in\Lambda}\PP(E_p)$.
\end{corollary}
\begin{proof}[Proof of Corollary~\ref{grlevi:cor:two-actions}]
Take a classical rank-two Schottky subgroup $\Gamma_0\cong F_2<\operatorname{SL}_2(\C)$, and identify it with the free group on two generators. By Kuranishi's theorem~\cite{Kuranishi}, choose homomorphisms
$\xi_1:F_2\to\operatorname{SU}(2)$ and
$\xi_r:F_2\to\operatorname{SU}(r)$ with dense images. Define
$\vartheta_-=\operatorname{Sym}^{r-1}\circ\xi_1$ and
$\vartheta_+=\xi_r$. The group $H_+=\operatorname{SU}(r)$ satisfies pair reversibility
by Lemma~\ref{grlevi:lem:swap-transitivity}.
The smaller group $H_- =\operatorname{Sym}^{r-1}(\operatorname{SU}(2))$
preserves the rational normal curve $\nu_{r-1}(\PP^1)$ of pure powers,
a proper subvariety of $\PP^{r-1}$ for $r\ge3$.
A line on this curve and one outside it cannot be interchanged.
Theorem~\ref{grlevi:thm:main} proves the opposite largest-domain behaviour.

Faithfulness and discreteness follow from Theorem~\ref{thm:marked-product}.

For strong irreducibility, the paired complex Zariski closures
project onto $\operatorname{SL}_2$ on the left and respectively onto
$\operatorname{SL}_r$ and $\operatorname{SL}_2$ on the right (before
taking the symmetric power). When $r\ge3$, $\operatorname{SL}_2$
and $\operatorname{SL}_r$ are non-isogenous simple groups, so algebraic
Goursat yields the full product. In the smaller case a proper
connected subdirect algebraic subgroup of
$\operatorname{SL}_2\times\operatorname{SL}_2$ would, up to finite
centres, identify the factors through an algebraic automorphism.
On a finite-index subgroup it would preserve squared traces.
But for powers of any loxodromic $a\in\Gamma_0$,
$|\operatorname{tr}(a^n)|\to\infty$, while
$|\operatorname{tr}(\xi_1(a)^n)|\le2$.
Thus the connected closure is the full product in both cases.
The outer tensor-product modules
$\C^2\otimes\C^r$ and
$\C^2\otimes\operatorname{Sym}^{r-1}\C^2$
are irreducible for these connected product groups.
Every finite-index subgroup remains Zariski dense: indeed, the
connected product group cannot be a finite union of proper Zariski-closed
cosets of the closure of such a subgroup. This proves strong
irreducibility. The carriers coincide by the common setup. On $\PP(W)$ the
condition $(\mathrm R_{1,0})$ is automatic, so
Theorem~\ref{thm:grk} gives the same largest domain for both twists.
\end{proof}

\begin{remark}\label{grlevi:rem:exterior-sip}
The strong irreducibility in Corollary~\ref{grlevi:cor:two-actions}
concerns the defining representation on $W=U\otimes V$, where
$U=\C^2$ and $\dim V=r$. It does not extend to all exterior
powers. Indeed, the natural $\operatorname{GL}(U)\times
\operatorname{GL}(V)$--equivariant decomposition
\[
 \bigwedge\nolimits^2(U\otimes V)
 \cong
 \bigl(\bigwedge\nolimits^2U\otimes\operatorname{Sym}^2V\bigr)
 \oplus
 \bigl(\operatorname{Sym}^2U\otimes\bigwedge\nolimits^2V\bigr)
\]
has two nonzero proper invariant summands for every $r\ge3$.
Thus neither satisfies exterior SIP, which requires strong
irreducibility and proximality in every exterior degree. Our contrast
concerns strong irreducibility of $W$, not an all-exterior-powers hypothesis.
\end{remark}

\begin{lemma}\label{lem:so3-r21}
On $V=\C^3$ the standard real rotation group $\mathrm{SO}(3)$ satisfies
$(\mathrm R_{2,1})$. It does not act transitively on $\PP(V)$.
\end{lemma}

\begin{proof}
Write the line as $\C x$ and the plane as $\ker y^t$, for nonzero
$x,y\in\C^3$ and the complex bilinear pairing $y^t x$. For $u$ on the
real unit sphere $S^2$, the half-turn $h_u=2uu^t-I$ belongs to
$\mathrm{SO}(3)$ and satisfies $h_u^{-1}=h_u$. It suffices to find a zero of
the continuous function $F(u)=y^t h_u x$.
Choose real orthonormal vectors $e,f$ whose span contains $\Re x,\Im x$.
On the great circle $u=e\cos t+f\sin t$, put $A=(y^te)(e^tx)-(y^tf)(f^tx)$
and $B=(y^te)(f^tx)+(y^tf)(e^tx)$. Expansion gives
\[
 F(u)=A\cos(2t)+B\sin(2t).
\]
If $A,B\in\C$ are linearly dependent over $\R$, this expression vanishes
for some $t$. Otherwise the real linear map $(v,w)\mapsto Av+Bw$
is invertible, so this loop winds twice about the origin, with sign given
by its orientation. A zero-free $F$ on a closed hemisphere bounded by
the great circle would extend this loop to a disk in $\C\setminus\{0\}$,
forcing winding number zero. This proves the assertion about
$(\mathrm R_{2,1})$. The projective action preserves the nonempty proper
subset of complex lines spanned by real vectors, which proves
nontransitivity.
\end{proof}

\begin{corollary}\label{cor:so3-odd}
For a twist with $r=3$, $\Lambda\subsetneq\PP(U)$ and
$H=\mathrm{SO}(3)$, the action has largest domain $\Omega^\star_3$
on $\operatorname{Gr}(3,\C^6)$ and has no largest domain on
$\operatorname{Gr}(2,\C^6)$. In particular, pair reversal on
$\operatorname{Gr}(1,V)$ is not necessary for a largest domain on
$\operatorname{Gr}(3,U\otimes V)$.
\end{corollary}

\begin{proof}
The cases $k=3$ and $k=2$ of Theorem~\ref{thm:grk-converse} require
$(\mathrm R_{2,1})$ and $(\mathrm R_{1,1})$, respectively. The lemma gives
the first. The second implies projective transitivity and hence fails.
The example is realised by the $r=3$ symmetric-square twist in
Corollary~\ref{grlevi:cor:two-actions}. Indeed, the alternating form
identifies the defining module of $\mathrm{SU}(2)$ with its complex dual.
Under the resulting decomposition of its endomorphism module, the scalar
summand corresponds to the alternating square, and the symmetric square
is the traceless endomorphism module. Conjugation preserves its real form
$\mathfrak{su}(2)$ and the positive inner product $-\operatorname{tr}(XY)$
there. The connected image lies in $\mathrm{SO}(3)$. Its kernel is
$\{\pm I\}$, since a matrix commuting with every traceless endomorphism is
scalar, so the image has dimension three and equals $\mathrm{SO}(3)$.
Complexifying this real action identifies the two compact images up to a
change of basis.
\end{proof}

\begin{corollary}\label{cor:sp-su}
Let $r=2n\ge4$ and let $\Gamma_0\cong F_2$ be a Schottky group as in Corollary~\ref{grlevi:cor:two-actions}. There are
two faithful, discrete, strongly irreducible twists $\Gamma_{\mathrm{Sp}}$ and $\Gamma_{\mathrm{SU}}$ of $\Gamma_0$, with compact closures
$\mathrm{Sp}(n)$ and $\mathrm{SU}(2n)$, such that:
\begin{enumerate}[label=\textup{(\alph*)},nosep]
\item every divergent sequence of either group has Cartan type $(r,r)$, with attracting and repelling $r$-planes $E_p$, $p\in\Lambda$;
\item their complete unmarked destination relations on $\PP(W)$ coincide, and so do their largest domains there;
\item on $\operatorname{Gr}(k,W)$, for $k\in\{1,2,3,2r-3,2r-2,2r-1\}$, both have the same largest domain $\Omega^\star_k$;
\item on $\operatorname{Gr}(4,W)$ the action of $\Gamma_{\mathrm{SU}}$ has the largest domain $\Omega^\star_4$, whereas
$\Gamma_{\mathrm{Sp}}$ has two non-recurrent $4$-planes related by $\D_{\Gamma_{\mathrm{Sp}}}$, and hence no largest domain.
\end{enumerate}
\end{corollary}

\begin{proof}
Choose dense homomorphisms $F_2\to\mathrm{Sp}(n)$ and $F_2\to\mathrm{SU}(2n)$, as in Proposition~\ref{prop:sharpness}. The
Zariski closures of the two compact factors are $\mathrm{Sp}(2n,\C)$ and $\mathrm{SL}(2n,\C)$, neither isogenous to
$\mathrm{SL}_2(\C)$, so the Goursat argument of Corollary~\ref{grlevi:cor:two-actions} gives strong irreducibility.
The common setup gives faithfulness, discreteness, (a), and identical
Cartan projections under the natural marking. Both groups are transitive on $\PP(V)$, so (b) follows from
Theorem~\ref{thm:unmarked}, applied in the dual module (the contragredient of a transitive group is transitive), and (c) from Theorem~\ref{thm:grk}. In (d), the twist by $\mathrm{SU}(2n)$
is covered by Theorem~\ref{thm:grk}. For $\mathrm{Sp}(n)$, Lemma~\ref{lem:Rab}(c) shows that pair reversal fails on
$\operatorname{Gr}(2,V)$, and Theorem~\ref{thm:grk-converse} with $b=2\le n$ applies.
\end{proof}

For $r=4$ the obstruction can be made explicit. Choose a symplectic basis $e_1,e_2,f_1,f_2$ of $V$, with $\omega(e_i,f_j)=\delta_{ij}$, and put
\[
K=\langle e_1,e_2\rangle,\qquad P=\langle f_1,e_1+f_2\rangle.
\]
Then $K$ is isotropic, $P$ is symplectic, and $K\cap P=0$. No element of $\mathrm{Sp}(2)$ can exchange them.
For $p,q\in\Lambda$ and $s,t\in\PP(U)\setminus\Lambda$, the planes
\[
S=q\otimes K\oplus s\otimes P,\qquad T=p\otimes P\oplus t\otimes K
\]
are the non-recurrent related pair furnished by Theorem~\ref{thm:grk-converse}.

The construction also applies to any non-elementary free Schottky group: for each compact factor, send two free generators to a dense generating pair and all remaining free generators to the identity. The Schottky factor is Zariski dense in $\mathrm{SL}_2(\C)$, so the same Goursat argument proves strong irreducibility. For $n=2$ this gives Theorem~\ref{thmC}, with the explicit witness
above; Proposition~\ref{prop:sharpness} gives its rank-two realisation.

\paragraph{Relation with earlier work.}
Left--right twisted actions were developed by Ghys~\cite{Ghys1995}.
Brion's compactification treats graph degenerations of complex flag
varieties~\cite{r2grBrion2003}; Kapovich--Leeb--Porti and Stecker treat
dynamic relations and properness on flag manifolds~\cite{KLP,SteckerBalanced},
and Carvajales--Stecker study homogeneous spaces with finitely many
relative-position orbits, including non-flag spaces~\cite{CarvajalesStecker2024}.
A close precedent is Stecker's reducible,
$\{\alpha_2\}$--Anosov (not full-flag Anosov) representation into
$\operatorname{Sp}(4,\R)$ with uncountably many incomparable maximal
domains on the real Lagrangian Grassmannian~\cite[Section~3.1.6]{SteckerThesis}.
Our criteria concern residual complex incidence and the full
Levi-twisted relation on complex Grassmannians; they give strongly
irreducible examples with identical carriers. We claim no priority
for incomparable maximal domains, projective transitivity
classifications or properness of balanced thickenings. The distinction
between an individual graph, its marked family and the unmarked union
remains essential: unmarked orbit recovery uses $s\notin\Lambda$;
no recovery of Levi subgroups is claimed for arbitrary discrete groups,
and equality of orbit partitions need not mean equality of compact groups.

\subsection{Lagrangian Grassmannians} Now let $V=\C^n$ carry the standard complex symmetric bilinear form $\beta(v,w)=v^tw$, and let
$\vartheta:\Gamma_0\to O(n)$ be real orthogonal, with closure $H\subseteq O(n)$. Then $\omega=\omega_U\otimes\beta$ is
symplectic on $W$, every $\gamma\otimes\vartheta(\gamma)$ lies in
$\mathrm{Sp}(W)\cong\mathrm{Sp}(2n,\C)$, and
$k_m\otimes\vartheta(\gamma_m)$ lies in its maximal compact subgroup. The type is the Siegel type $\Theta=\{\alpha_n\}$, and the planes $E_p=p\otimes V$ are
Lagrangian. Let $X=\mathrm{LGr}(W)$. Here isotropy ties the two
residual subspaces together, so recurrence will be governed by one
inclusion involving $\beta$-orthogonal complements. We write $A^{\perp\beta}$ for the bilinear orthogonal complement of a complex subspace $A\subseteq V$, and $A^\perp$ for its Hermitian orthogonal complement.

\begin{proposition}\label{prop:lag-relation}
For a Lagrangian $S$ and $q\in\PP(U)$ one has $P_q(S)=K_q(S)^{\perp\beta}$. On $\mathrm{LGr}(W)$,
\[
\D_{\Gamma_\vartheta}=\bigcup_{p,q\in\Lambda,\ h\in H}R^L_{p,q,h},\qquad
R^L_{p,q,h}(S)=\{T\in\mathrm{LGr}(W):\ T\supseteq p\otimes hP_q(S)\}.
\]
\end{proposition}

\begin{proof}
Write $s=q\otimes a+q'\otimes b\in S$ with $q'\notin q$. For $k\in K_q(S)$, $\omega(q\otimes k,s)=\omega_U(q,q')\beta(k,b)$, so
$\beta(K_q(S),P_q(S))=0$. Writing $d=\dim K_q(S)$, the dimensions are $d$ and $n-d$, which gives the first assertion. The necessary target inclusion follows from~\eqref{grlevi:eq:allk};
its second condition is automatic for a Lagrangian containing
$p\otimes hP_q(S)$, since
\[
 T\subseteq(p\otimes hP_q(S))^{\perp\omega}
   =E_p\oplus p^\perp\otimes hK_q(S).
\]
For sufficiency take the two-block model with $V^\pm$ Lagrangian,
$K=S\cap V^-$ and $P=\operatorname{pr}_+S$.
The upper symplectic unipotent group $N_\Theta$ consists of
$\omega$-symmetric maps $V^-\to V^+$. Its orbit of $P\oplus K$
consists of the Lagrangians containing $P$ with intersection $P$
with $V^+$. In the reduced space $P^{\perp\omega}/P$ these are
graphs of symmetric maps, open and dense among the Lagrangians
containing $P$. Thus the orbit closure is exactly that closed set.

We now realise every destination by moving sources that remain Lagrangian. Put $S_0=P\oplus K$, which is Lagrangian.
There is a lower symplectic unipotent $u^-$ with $S=u^-S_0$:
a section of $S\to P$ defines a symmetric bilinear form on $P$,
which extends to $V^+$.
For an orbit target write $T=u^+S_0$, with $u^+$ upper symplectic unipotent.
With $a_m(u,v)=(s_mu,s_m^{-1}v)$ and $s_m\to\infty$, set
\[
 S_m=u^-a_m^{-1}u^+a_mS_0.
\]
All factors are symplectic, so every $S_m$ is Lagrangian.
The off-diagonal blocks of $a_m^{-1}u^+a_m$ and
$a_mu^-a_m^{-1}$ are multiplied by $s_m^{-2}$, and
$a_mS_0=S_0$. Hence
\[
 S_m\longrightarrow S,\qquad
 a_mS_m=(a_mu^-a_m^{-1})u^+S_0\longrightarrow T.
\]
For a boundary target choose orbit targets $T_j\to T$.
For each fixed $j$, the associated sources $S_m^{(j)}$ satisfy both limits above.
Choose increasing $N_j$ so that, for all $m\ge N_j$,
\[
 d_{\mathrm L}(S_m^{(j)},S)+d_{\mathrm L}(a_mS_m^{(j)},T_j)<1/j
\]
where $d_{\mathrm L}$ is any compatible metric.
Taking $j(m)\to\infty$ with $m\ge N_{j(m)}$
gives the required Lagrangian sources along the whole sequence.

Transporting to the twisted frames and applying
Lemma~\ref{lem:saturate}, as in Theorem~\ref{grlevi:thm:allk}, gives
all triples $(p,q,h)$.
\end{proof}

Relative to a Lagrangian $E$, the positions are
$\tau_j=\{S:\dim(S\cap E)=j\}$, $0\le j\le n$, ordered as before,
with $w_0\tau_j=\tau_{n-j}$. Hence
\[
 \mathcal T_\star=\{\tau_j:j\ge n/2\},\qquad
 \Sigma^L=\{S\in\mathrm{LGr}(W):\dim(S\cap E_p)\ge n/2
                                  \text{ for some }p\in\Lambda\}.
\]
The ideal is fat, and balanced exactly when $n$ is odd.
The Levi group is $GL(E)$, with $Y_{\tau_j}\cong\operatorname{Gr}(j,E)$
represented by $A\oplus\operatorname{ann}(A)$.
The orbit closure just computed gives
$\mathcal I_{\tau_j}=\{(A,A'):A'\subseteq A\}$.
For $A=p\otimes A_V$, its annihilator is
$p^\perp\otimes A_V^{\perp\beta}$, and $w\otimes h$ sends the
residual subspace to $h(A_V^{\perp\beta})$.
Proposition~\ref{prop:lag-relation}, or Theorem~\ref{thm:diagonal}(c),
therefore gives the following.

\begin{proposition}\label{prop:lag}
A Lagrangian $S$ with $S\cap E_p=p\otimes A$ is in the diagonal of $(p,p,h)$ if and only if $h(A^{\perp\beta})\subseteq A$.
Hence $E_p$ is a recurrence seed if and only if
\[
(\mathrm L)\qquad \forall j\ge n/2,\ \forall A\in\operatorname{Gr}(j,\C^n)\ \ \exists h\in H:\ h(A^{\perp\beta})\subseteq A .
\]
\end{proposition}

\begin{lemma}\label{lem:normal}
For every complex subspace $A\subseteq\C^n$ there are a real orthonormal family $e_1,\dots,e_t,e'_1,f'_1,\dots,e'_s,f'_s$ and
angles $0<\theta_i\le\pi/4$ with
\[
A=\langle e_1,\dots,e_t\rangle_\C\oplus\langle\cos\theta_i\,e'_i+\mathrm i\sin\theta_i\,f'_i\rangle_\C .
\]
Consequently $SO(n)$ satisfies $(\mathrm L)$.
\end{lemma}

\begin{proof}
By Takagi's factorisation of the symmetric form $\beta|_A$, there is a Hermitian orthonormal basis $a_i=x_i+\mathrm iy_i$
of $A$ with $\beta(a_i,a_l)=\sigma_i\delta_{il}$, $0\le\sigma_i\le1$. Comparing real and imaginary parts of
$\langle a_i,a_l\rangle=\delta_{il}$ and $\beta(a_i,a_l)=\sigma_i\delta_{il}$ gives $x_i\cdot y_l=0$ for all $i,l$, together with
$2x_i\cdot x_l=(1+\sigma_i)\delta_{il}$ and $2y_i\cdot y_l=(1-\sigma_i)\delta_{il}$. Normalising gives the family, with
$\cos^2\theta_i=(1+\sigma_i)/2$; the vectors with $\sigma_i=1$ are real.

We now use this normal form to prove $(\mathrm L)$. Let $R$ be the real orthogonal complement of the
family, of dimension $n-t-2s$. Then $A^{\perp\beta}=\overline{A^{\perp}}$ is spanned by $R_\C$ and the vectors
$-\mathrm i\sin\theta_i\,e'_i+\cos\theta_i\,f'_i$. The rotation $e'_i\mapsto-f'_i$, $f'_i\mapsto e'_i$ in each plane maps
these vectors to $\cos\theta_i\,e'_i+\mathrm i\sin\theta_i\,f'_i$. If $j=\dim A\ge n/2$, then $n-t-2s\le t$, so an isometry maps
$R$ into $\langle e_1,\dots,e_t\rangle$. The sum of these maps extends to an element of $O(n)$ carrying
$A^{\perp\beta}$ into $A$.

If $t>0$ and this extension has determinant $-1$, precompose it with the real reflection that changes the sign of $e_1$ and fixes its real orthogonal complement. This reflection fixes $A^{\perp\beta}$ pointwise, so the inclusion is preserved and the determinant becomes $1$. If $t=0$, the inequalities $j=s\ge n/2$ and $2s\le n$ force $n=2s$; the product of the rotations on the displayed real two-planes already has determinant $1$. This gives an element of $SO(n)$ in all cases.
\end{proof}

\begin{theorem}\label{thm:lag}
If $H$ satisfies $(\mathrm L)$, in particular if $H\supseteq SO(n)$, the recurrent set of $\Gamma_\vartheta$ on
$\mathrm{LGr}(W)$ is $\Sigma^L$, and
\[
\Omega^L=\{S\in\mathrm{LGr}(W):\ \dim(S\cap E_p)<n/2\ \text{for all }p\in\Lambda\}
\]
is the largest domain of proper discontinuity. If $H\subseteq\{\pm I\}$ and $n\ge2$, condition $(\mathrm L)$ fails.
\end{theorem}

\begin{theorem}\label{thm:lag-converse}
Assume $\Lambda\subsetneq\PP(U)$ in the orthogonal-twist setup
above, with $n\ge1$, and put $j=\lceil n/2\rceil$.
The action of $\Gamma_\vartheta$ on $\mathrm{LGr}(W)$ has a largest
domain of proper discontinuity if and only if
\[
 (\mathrm L_j)\qquad
 \forall A\in\operatorname{Gr}(j,V)\quad
 \exists h\in H\quad h(A^{\perp\beta})\subseteq A .
\]
In this case the recurrent set is $\Sigma^L$ and the largest domain is the
nonempty set $\Omega^L$. Otherwise two non-recurrent Lagrangians are related
by $\D_{\Gamma_\vartheta}$. For $n=2m$ the inclusion in $(\mathrm L_m)$
is equality.
\end{theorem}

\begin{proof}[Proof of Theorems~\ref{thm:lag} and~\ref{thm:lag-converse}]
For $n=1$ the only $j$-plane is $V$, and $V^{\perp\beta}=0$.
The Lagrangian Grassmannian is $\PP(U)$ and the twist acts by projective
scalars on the second factor. Proposition~\ref{prop:D-rank1} gives
$\Rec_{\Gamma_\vartheta}=\Lambda$ and largest domain $\PP(U)\setminus\Lambda$.

For $n\ge2$, Proposition~\ref{prop:lag-relation} sends every source
outside $\Sigma^L$ to a point inside it: if $d_q<n/2$, its destination
contains $p\otimes hP_q(S)$ of dimension $n-d_q>n/2$.
Thus $\Omega^L$ is a properness domain for any $H$ and any $\Lambda$.
Under $(\mathrm L)$ every point of the wall is recurrent by
Proposition~\ref{prop:lag}, so Theorem~\ref{thm:relational} gives the
largest domain. Lemma~\ref{lem:normal} supplies $(\mathrm L)$ for
$H\supseteq SO(n)$; equivalently one may apply
Theorem~\ref{thm:singular-recurrence}. If $H\subseteq\{\pm I\}$,
a nondegenerate $\lceil n/2\rceil$-plane $A$ violates
$A^{\perp\beta}\subseteq A$, so $(\mathrm L)$ fails.

If $(\mathrm L_j)$ holds and $\dim B\ge j$, choose a $j$-plane
$A\subseteq B$. An $h$ working for $A$ satisfies
$h(B^{\perp\beta})\subseteq h(A^{\perp\beta})\subseteq A\subseteq B$.
Thus $(\mathrm L)$ holds. When $\Lambda$ is proper,
$E_s\in\Omega^L$ for every $s\notin\Lambda$.

Suppose instead that $(\mathrm L_j)$ fails. The $j$-planes satisfying its
inclusion for some $h$ form a closed subset, by compactness of $H$,
continuity of $A\mapsto A^{\perp\beta}$ and closedness of incidence.
The locus where $\beta|_A$ is degenerate is a proper algebraic subset
of the Grassmannian: the Gram determinant vanishes there, whereas it is
nonzero on a coordinate $j$-plane. We can therefore choose a violating
$A$ on which $\beta$ is nondegenerate. Write $B=A^{\perp\beta}$, so
$V=A\oplus B$ and $\dim B=n-j>0$. Pick $q\in\Lambda$ and
$s\in\PP(U)\setminus\Lambda$ and put
\[
 S=q\otimes A\oplus s\otimes B,\qquad
 T=q\otimes B\oplus s\otimes A .
\]
Both have dimension $n$ and are isotropic, since $\beta(A,B)=0$, and
hence are Lagrangian. A tensor with nonzero parts in both summands
has rank two, since $q\ne s$ and $A\cap B=0$. Thus $E_q$ is the only
$E_{q'}$ with $q'\in\Lambda$ that meets either plane nontrivially.
Here $K_q(S)=A$, $P_q(S)=B$, $K_q(T)=B$ and $P_q(T)=A$.

If $S\in R^L_{p',q',h}(S)$, Proposition~\ref{prop:lag-relation}
requires $p'\otimes hP_{q'}(S)\subseteq S$. For $q'\ne q$ one has
$P_{q'}(S)=V$, which would force $S=E_{p'}$. For $q'=q$ the inclusion
uses a nonzero subspace $hB$, and hence forces $p'=q$ and $hB\subseteq A$.
Both alternatives are impossible. Thus $S$ is non-recurrent.
If $n$ is odd, $\dim B<n/2$, so $T\notin\Sigma^L$ and the
properness already proved makes it non-recurrent. If $n$ is even,
the same recurrence test used for $S$ would require $hA\subseteq B$.
Dimension gives $hA=B$ and then $h^{-1}B=A$, contradicting the
choice of $A$. In both cases $T$ contains $q\otimes B=q\otimes P_q(S)$.
Hence $T\in R^L_{q,q,I}(S)$, and Theorem~\ref{thm:relational} excludes a
largest domain.
\end{proof}

\begin{corollary}\label{cor:lag-n2}
Let $n=2$, so that $\mathrm{LGr}(W)$ is the three-dimensional quadric, and assume $\Lambda\subsetneq\PP(U)$. A largest domain of
proper discontinuity exists if and only if $R_{\pi/2}\in H$. It is then $\{S:S\cap E_p=0\text{ for all }p\in\Lambda\}$. If $R_{\pi/2}\notin H$, there are two non-recurrent Lagrangians related by $\D_{\Gamma_\vartheta}$.
\end{corollary}

\begin{proof}
On lines of $\C^2$, the polarity is
$A^{\perp\beta}=R_{\pi/2}A$, with precisely the two isotropic fixed
lines $[1:\mathrm i]$ and $[1:-\mathrm i]$. If $R_{\pi/2}\in H$, then $h=R_{\pi/2}^{-1}$ satisfies $hA^{\perp\beta}=A$ for every
line $A$, and the case $j=2$ of $(\mathrm L)$ is empty; Theorem~\ref{thm:lag} applies. If $R_{\pi/2}\notin H$, then $H$ does not
contain $SO(2)$, so it is finite, and $-R_{\pi/2}=R_{\pi/2}^{-1}\notin H$ either. For $h\in H$ the matrix $hR_{\pi/2}$ is therefore not
scalar, and it fixes at most two lines. A line $A$ outside these finitely many lines violates $(\mathrm L)$, and
Theorem~\ref{thm:lag-converse} applies.
\end{proof}

For the trivial twist, $A\mapsto A\otimes I_2$ is the complexification
of the real embedding in Stecker's Fuchsian representation
$\iota\circ\rho_0<\mathrm{Sp}(4,\R)$~\cite[\S3.6]{SteckerBalanced},
whose action on the real Lagrangian Grassmannian has infinitely many
maximal domains. The corollary concerns the complex Grassmannian.
If $\Gamma_0$ is free and $\Lambda\subsetneq\PP(U)$, sending one
generator to $R_{\pi/2}$ and the others to $I$ gives closure $\langle R_{\pi/2}\rangle\cong\Z/4$ and the same type,
Lagrangian limit planes and wall as a twist trivial up to sign.
The finite quarter-turn makes that wall recurrent and its complement
the largest domain; the twist with $H\subseteq\{\pm I\}$ has none.

For general $n$, the rotation in Lemma~\ref{lem:normal} depends on $A$.
For an $m$-plane $A\subset\C^{2m}$, a fixed complex structure
$J\in O(2m)$ satisfies $J(A^{\perp\beta})=A$ precisely when $A$ is
Lagrangian for $\beta(J\cdot,\cdot)$. For $n\ge4$ a single such map
does not classify the compact subgroups satisfying $(\mathrm L)$.

\subsection{Further questions}
Theorems~\ref{thm:grk-converse} and~\ref{thm:lag-converse} give exact criteria for compact twists with proper classical limit set. The remaining classification problem concerns compact group actions. Which closed subgroups of $U(r)$ satisfy $(\mathrm R_{a,b})$ for specified $a,b$, particularly when $a=b+1$? Lemma~\ref{lem:so3-r21} shows why a classification by pair reversal alone would miss the odd-dimensional case. Likewise, which proper closed subgroups of $\mathrm{SO}(n)$ satisfy $(\mathrm L_{\lceil n/2\rceil})$? For $n=2m$ this is $(\mathrm L_m)$. Beyond these twist families, Remark~\ref{rem:fat-question} asks whether simplicity and Zariski density restore the necessity of fatness in the regular case.

\section{Projective applications}\label{sec:applications}
\subsection{The critical level and the largest domain}\label{sec:critical}

Let $G=\PSL(V)$, $V=\C^d$, $d=n+1\ge2$. We study when the middle Cartan wall consists of recurrent points. The definitions through Theorem~\ref{thm:crit-recurrence} also apply over $\R$. Assume $\Gamma$ infinite; for finite groups the wall is empty and the whole space is the largest domain.

\medskip
\emph{Block notation.} For an extracted divergent sequence
$g_m=k_ma_mk'_m$, use Proposition~\ref{prop:proj}: write
$V=\bigoplus_{j=1}^rV_j$, $m_j=\dim V_j$, $r_j=\sum_{i\le j}m_i$,
and $a_m|_{V_j}=z_{j,m}D_{j,m}$ with
$z_{j+1,m}/z_{j,m}\to0$ and $D_{j,m}\to D_j$ invertible.
Set $r_0=0$, $G_0=F_{r+1}=0$, and
\[
G_j=k(V_1\oplus\cdots\oplus V_j),\qquad
F_j=k'^{-1}(V_j\oplus\cdots\oplus V_r).
\]
Thus $\dim G_j=r_j$, $\dim F_j=d-r_{j-1}$, and the $F_j$
are the inverse target spaces. For $x=[\xi]$ with
$\xi\in F_j\setminus F_{j+1}$, Proposition~\ref{prop:proj} gives
\begin{equation}\label{eq:blockform}
\D_{(g_m)}(x)=\PP(G_{j-1}\oplus\C v_x)\subseteq\PP(G_j),
\qquad v_x=kD_j\operatorname{pr}_{V_j}(k'\xi).
\end{equation}
In particular $E_{(g_m)}=\PP(F_2)$. The first pseudo-projective
limit is $\tau=[kD_1\pi_1k']$, where $\pi_1$ projects onto
$V_1$ along the other blocks; its kernel and image are $F_2$ and $G_1$.

\begin{proposition}\label{prop:kernel}
Let $\mathcal L_\Gamma\subset\PP(\mathrm{End}\,V)$ be the set of limits of sequences of distinct elements of $\Gamma$, the pseudo-projective boundary of $\Gamma$ in the sense of Cano--Seade~\cite{CanoSeade2010}. Then
\[
\Lambda_{\Myr}(\Gamma)=\Lambda_{\mathrm{exc}}(\Gamma)=\bigcup_{\tau\in\mathcal L_\Gamma}\PP(\ker\tau),\qquad\Omega_{\Myr}(\Gamma)=\Eq(\Gamma),
\]
the set of points at which $\Gamma$ is equicontinuous.
\end{proposition}
\begin{proof}
For an extracted sequence, formula~\eqref{eq:blockform} gives
$E_{(g_m)}=\PP(\ker\tau)$, where $\tau$ is its first
pseudo-projective limit. Conversely, every
$\tau\in\mathcal L_\Gamma$ arises in this way. Thus
$\Lambda_{\mathrm{exc}}$ is the displayed union of kernels.
The set $\mathcal L_\Gamma$ is compact by diagonal extraction,
so the incidence set
$\{(\tau,x):\tau\in\mathcal L_\Gamma,\ x\in\PP(\ker\tau)\}$
is compact. Its projection is closed, giving
$\Lambda_{\Myr}=\Lambda_{\mathrm{exc}}$.

Equicontinuity at $x$ forces any two moving-source limits along
one sequence to coincide. Conversely, if equicontinuity fails at
$x$, there are $x_m\to x$ and $\gamma_m\in\Gamma$ with
$d(\gamma_mx_m,\gamma_mx)$ bounded below. The elements may be
chosen distinct after extraction; extract again so that
$\gamma_m\to\tau\in\mathcal L_\Gamma$. If
$x\notin\PP(\ker\tau)$, both image sequences converge to
$\tau x$, a contradiction. This proves the asserted equality
with the equicontinuity region.
\end{proof}

This identifies $\Omega_{\Myr}$ with the usual equicontinuity region; for complex hyperbolic groups and $\PSL(3,\C)$ see~\cite{CanoSeade2010,BCNLimit2011,CLL} and~\cite[Sections~7.4--7.5]{CNS}.

The \emph{Conze--Guivarc'h set} $\Lambda_{\CoG}$ is the closure of
rank-one images $\PP(\mathrm{Im}\,\tau)$, $\tau\in\mathcal L_\Gamma$.
For strongly irreducible $\Gamma$ containing a proximal element it is
the unique minimal closed invariant set~\cite[author version, Lem.~3.2(c)]{BenoistQuint}; see~\cite{ConzeGuivarch,AbelsMargulisSoifer}
for the earlier real theory. It equals $\Lambda^{(0)}_{\mathrm Z}$;
more generally $\Lambda^{(2k)}_{\mathrm Z}$ is the closure of the
union of target spaces $\PP(G_j)$ of complex dimension at most $k$.
Each rank-one image is the limit of every fixed source outside
$\PP(\ker\tau)$, a dense set, so it belongs to $B_\Gamma$.
Theorem~\ref{thm:L} therefore gives, without irreducibility assumptions,
\begin{equation}\label{eq:proj-chain}
\Lambda_{\CoG}\subseteq B_\Gamma\subseteq\mathscr R_\Gamma\subseteq\Lambda_{\Kul}\subseteq\Lambda_{\Myr}=\C\PP^n\setminus\Eq.
\end{equation}

\begin{example}\label{ex:cp5}
A cyclic action on $\C\PP^5$ separates the Cartan levels from the
Kulkarni and Myrberg sets. Let
$V=E\oplus F$, where $E=\langle e_1,e_2,h\rangle$ and
$F=\langle f_1,f_2,f_3\rangle$, and let $\Gamma$ be generated by
\[
 g=\diag(4e^{N_2},4,\tfrac14e^{N_3})
 \quad\text{on }\langle e_1,e_2\rangle\oplus\langle h\rangle\oplus F,
\]
where $N_2e_2=e_1$, $N_3f_3=f_2$, $N_3f_2=f_1$, and the
remaining images are zero. For $m\to+\infty$, the singular values
of $g^m$ are, up to bounded positive factors,
\[
 4^m m,\quad4^m,\quad4^m/m,\quad
 4^{-m}m^2,\quad4^{-m},\quad4^{-m}m^{-2}.
\]
Thus the type is the full flag. Put
$A=E+\langle f_1,f_2\rangle$ and $C=F+\langle e_1\rangle$.
The target flags of $g^m$ and $g^{-m}$ are respectively
\[
 \langle e_1\rangle\subset\langle e_1,h\rangle\subset E
 \subset E+\langle f_1\rangle\subset A\subset V,
\]
\[
 \langle f_1\rangle\subset\langle f_1,f_2\rangle\subset F
 \subset C\subset F+\langle e_1,h\rangle\subset V.
\]
The fixed-source set is
$B_\Gamma=\PP\langle e_1,h\rangle\cup\{[f_1]\}$.
For sources outside $B_\Gamma$, formula~\eqref{eq:blockform}
reaches at most the fifth forward target space $A$ and the fourth
backward target space $C$; these entire spaces occur from the
sources $[f_2]$ and $[e_2]$, respectively. The subsequent source
strata lie in $B_\Gamma$. Proposition~\ref{prop:L2} therefore gives
\[
 \Lambda_{\Kul}=\PP(A)\cup\PP(C),\qquad
 \Lambda_{\Myr}=\PP(A)\cup\PP(F+\langle e_1,h\rangle).
\]
Moreover,
\[
 \Lambda^{(6)}_{\mathrm Z}
 =\PP(E+\langle f_1\rangle)\cup\PP(C)
 \subsetneq\Lambda_{\Kul}
 \subsetneq\Lambda^{(8)}_{\mathrm Z}=\Lambda_{\Myr}.
\]
Indeed, $[e_2+f_2]\in\Lambda_{\Kul}\setminus\Lambda^{(6)}_{\mathrm Z}$
and $[h+f_3]\in\Lambda_{\Myr}\setminus\Lambda_{\Kul}$.
\end{example}

\medskip
\emph{The critical wall.} Put $c=\lceil d/2\rceil$. An extracted divergent sequence $h=(h_m)$ is \emph{admissibly oriented} if $r_{r-1}\ge c$, that is, if its last block has dimension at most $d-c$; at least one of $h$, $h^{-1}$ is admissible, since $m_1+m_r\le d$.

For admissible $h$ let $t_h=\min\{j:r_j\ge c\}<r$ and $B_h=G_{t_h}$, the first target space reaching the middle. Let
\[
p_\Gamma=\min\{\dim B_h:\ h\text{ admissible}\}\in[c,d-1],\qquad q_\Gamma=p_\Gamma-1,
\]
and let $B^*_h$ be the smallest target space of $h^{-1}$ of dimension $\ge d-p_\Gamma$. Define
\[
\Lambda_{\mathrm{crit}}(\Gamma)=\overline{\textstyle\bigcup_{h\text{ admissible}}\bigl(\PP(B_h)\cup\PP(B^*_h)\bigr)},\qquad\Omega_{\mathrm{crit}}=\C\PP^n\setminus\Lambda_{\mathrm{crit}}.
\]
The paired space $B_h^*$ may exceed dimension $d-p_\Gamma$ if this threshold lies inside an inverse block; it is retained even when $h^{-1}$ is not admissible. For $n=1$ the wall is the classical limit set.

\begin{proposition}\label{prop:crit-cover}
$\Lambda_{\mathrm{crit}}$ is a closed invariant vertex cover. Hence $\Gamma$ acts properly discontinuously on $\Omega_{\mathrm{crit}}$, $\mathscr R_\Gamma\subseteq\Lambda_{\mathrm{crit}}$, and the trichotomy of Corollary~\ref{cor:trichotomy} holds for $\Lambda_{\mathrm{crit}}$.
\end{proposition}
\begin{proof}
Multiplication by fixed elements changes singular values by bounded factors.
Conjugating $h_m$ by $g\in\Gamma$ therefore preserves admissibility
and, after extraction, sends the target and source flags to $(gG_j)$
and $(gF_j)$. Its selected spaces are $gB_h,gB_h^*$, proving invariance.

To prove the vertex-cover property, let $(x,y)\in\D_\Gamma$ be witnessed by $h$; inverting $h$ and exchanging $x,y$ if necessary, $h$ is admissible. Put $t=t_h$. Since $p_\Gamma\le r_t$, the inverse target space $F_{t+1}$ has dimension $d-r_t\le d-p_\Gamma$, so $F_{t+1}\subseteq B^*_h$. If $x\notin\PP(B^*_h)$, the source stratum of $x$ is $F_j\setminus F_{j+1}$ with $j\le t$, and~\eqref{eq:blockform} gives $y\in\PP(G_t)=\PP(B_h)$. Theorem~\ref{thm:relational} concludes.
\end{proof}

For $\Gamma=\langle[\diag(2,1,1)]\rangle$, only negative powers are admissible and $\Lambda_{\mathrm{crit}}=\PP\langle e_2,e_3\rangle\cup\{[e_1]\}=\mathscr R_\Gamma$; allowing both orientations would fill $\C\PP^2$. In general, largestness requires recurrence of every point of the critical wall.

\medskip
\emph{Uniform critical type.} Fix integers $1\le a<p<d$ with $a+p=d$. The group has \emph{self-opposed uniform critical type} $(a,p)$ if along every divergent sequence, after extraction,
\begin{equation}\label{eq:uniform-type}
\frac{\sigma_a}{\sigma_{a+1}}\to\infty,\qquad\frac{\sigma_p}{\sigma_{p+1}}\to\infty,\qquad\sup_m\frac{\sigma_{a+1}}{\sigma_p}<\infty,
\end{equation}
so the middle crossing is one block in positions $a+1,\dots,p$.
Since $(d-p,d-a)=(a,p)$, inverses have the same crossing, both
orientations are admissible, and $p_\Gamma=p$. The paired space
$B_h^*$ is the inverse critical $a$-space, contained in the already
selected $p$-space $B_{h^{-1}}$. Thus, if
$\mathscr F\subset\mathrm{Fl}(a,p;V)$ is the closure of the
target critical flags $(A\subset B)=(G_{t-1}\subset G_t)$,
\[
\Lambda_{\mathrm{crit}}=\textstyle\bigcup_{(A\subset B)\in\mathscr F}\PP(B).
\]
Let $W=W_{a,p}$ be the Cartan component of $\bigwedge^aV\otimes\bigwedge^pV$, with the Pl\"ucker--Segre embedding $\jmath:\mathrm{Fl}(a,p;V)\to\PP(W)$ and the induced representation $\tau$. The \emph{critical module} is the span of the lines $\jmath(F)$, $F\in\mathscr F$; it is \emph{full} if it equals $W$. The action is \emph{critically strongly irreducible} if every finite-index subgroup acts irreducibly on $V$ and on $W$.

For a singular value decomposition $h=U\diag(\sigma_1,\dots,\sigma_d)V^*$ write $\widehat A^+_h=U\langle e_1,\dots,e_a\rangle$, $\widehat B^+_h=U\langle e_1,\dots,e_p\rangle$, and $\widehat A^-_h=V\langle e_{p+1},\dots,e_d\rangle$, $\widehat B^-_h=V\langle e_{a+1},\dots,e_d\rangle$ (the slow source flag, which is the target flag of $h^{-1}$).

\medskip\noindent\emph{The diagonal criterion.} Put $\Theta=\{\alpha_a,\alpha_p\}$. For a flag
$\xi=(A\subset B)$ with $K$-opposite frame $V=A\oplus C\oplus H$, $B=A\oplus C$, there are three cells:
\[
c_0=\PP(A),\qquad c_1=\PP(B)\setminus\PP(A),\qquad c_2=\PP(V)\setminus\PP(B),
\]
with $Y_{c_0}=\PP(A)$, $Y_{c_1}=\PP(C)$ and $Y_{c_2}=\PP(H)$. Here $w_0$ exchanges $c_0$ and $c_2$ and fixes $c_1$. So
$\mathcal T_\star=\{c_0,c_1\}$ is fat and $\mathcal T_\star(\xi)=\PP(B)$. The group $N_\Theta$ consists of the block unipotents mapping
$C\to A$ and $H\to A\oplus C$. Hence $\overline{N_\Theta[h]}=\PP(A\oplus C\oplus\langle h\rangle)\supseteq\PP(A)$, so
$\mathcal I_{c_0}=Y_{c_0}\times Y_{c_2}$, and $\overline{N_\Theta[c]}\cap\PP(C)=\{[c]\}$, so $\mathcal I_{c_1}$ is the diagonal.

The bounded middle block in~\eqref{eq:uniform-type} prevents further Levi degeneration on $\PP(C)$.

\begin{proposition}\label{prop:proj-diagonal}
Let $\Gamma<\operatorname{PSL}(V)$ have self-opposed uniform critical type $(a,p)$, and let $\eta$ be any datum coalescing at
$\Theta$ with flag $\xi=(A\subset B)$ and effective Levi element $e_\eta$. In the $K$-opposite frame attached to $\xi$, put
$e_{\xi,\eta}=ke_\eta k^{-1}$ on $\PP(C)$ and $\rho^+_\xi(x)=k\rho^+(k^{-1}x)$. The diagonal of $\D_\eta$ is
\[
\PP(A)\ \cup\ \{x\in\PP(B)\setminus\PP(A):\ e_{\xi,\eta}\,\rho^+_\xi(x)=\rho^+_\xi(x)\}.
\]
In particular $\PP(B)$ lies in the diagonal if and only if $e_{\xi,\eta}$ acts trivially on $\PP(C)$. This is automatic when $p=a+1$.
\end{proposition}

\begin{proof}
Work in the standard frame, writing $C_0=k^{-1}C$ and $H_0=k^{-1}H$.
Uniform critical type means that the middle block is bounded, so $\Theta'\setminus\Theta$ contains only roots $\alpha_i$ with $i<a$ or $i>p$.
Thus $N^L_{\Theta'}$ acts trivially on $C_0$, and $A_{\Theta'}$ acts on $C_0$ by scalars. Hence $\D^L_\eta(q)=\{q\}$ on $\PP(C_0)$.
For $x\in\PP(B)\setminus\PP(A)$, put $p=\rho^+(k^{-1}x)$. Theorem~\ref{thm:diagonal} gives the criterion $e_\eta p=p$, which, after applying $k$, is exactly $e_{\xi,\eta}\rho^+_\xi(x)=\rho^+_\xi(x)$.
On $c_0$, the relevant Levi destination is a nonempty subset of $\PP(H_0)$ and $\mathcal I_{c_0}$ is full, so part~(b) of that theorem applies. Points of $c_2$ are excluded by part~(a).
\end{proof}

\begin{definition}\label{def:seed}
A flag $F=(A\subset B)\in\mathscr F$ is a
\emph{regular critical seed} if there are $s_j\in\Gamma$ with
the following properties. The elements $s_j$ have strict spectral
gaps $|\mu_a|>|\mu_{a+1}|$ and $|\mu_p|>|\mu_{p+1}|$,
and their attracting and repelling spectral flags
$F^\pm_j=(A^\pm_j\subset B^\pm_j)$ both tend to $F$.

Their middle spectral blocks $C_j=B^-_j\cap B^+_j$, with
$V=A^-_j\oplus C_j\oplus A^+_j$, carry relatively compact
cyclic projective actions
$\overline{\langle[s_j|_{C_j}]\rangle}\subset\PGL(C_j)$.
After extraction, $C_j\to C_\infty$ with $B=A\oplus C_\infty$.
\end{definition}

Theorem~\ref{thm:singular-recurrence} propagates a recurrence seed
through a minimal $\mathscr F$. We first prove minimality and,
when $p=a+1$, construct coalescing spectral flags.

\begin{lemma}\label{lem:rec-contractions}
If the critical module is full and the action is critically strongly irreducible, then $\tau(\Gamma)$ contains a proximal element and $\mathscr F=\overline{\Gamma F_0}$ for every $F_0\in\mathscr F$. If moreover $p=a+1$, then for every $F\in\mathscr F$ there are $s_j\in\Gamma$ with strict spectral gaps at $a$ and $p$ whose attracting and repelling spectral flags both tend to $F$.
\end{lemma}
\begin{proof}
The highest singular line of
$\bigwedge^aV\otimes\bigwedge^pV$ lies in $W$, has exponent
$\sum_{i\le a}\log\sigma_i+\sum_{i\le p}\log\sigma_i$,
and its second-to-first singular ratio is
$\max\{\sigma_{a+1}/\sigma_a,\sigma_{p+1}/\sigma_p\}$.
Since $W$ is a compact-invariant orthogonal summand, restriction
cannot increase this ratio. Thus a sequence $h_m$ realising $F$
has normalised limit $T=v\otimes\phi$ of rank one with
$[v]=\jmath(F)$. Irreducibility supplies $g\in\Gamma$ with
$\phi(\tau(g)v)\ne0$; then $T\tau(g)$ has a simple nonzero
eigenvalue, and $\tau(h_mg)$ is proximal for large $m$.

Left multiplication moves the rank-one exterior images by
$\bigwedge^ag,\bigwedge^pg$, so $\mathscr F$ is invariant.
For any $F_0\in\mathscr F$, irreducibility also supplies $g$
with $\jmath(gF_0)\notin\PP(\ker T)$, whence $h_mgF_0\to F$.
This proves $\mathscr F=\overline{\Gamma F_0}$.

Suppose $p=a+1$ and $\tau(s)$ is proximal. In the commuting
Jordan decomposition $s=s_{\mathrm s}s_{\mathrm u}$, the
semisimple factor determines the eigenvalue moduli. Order them
as $|\mu_1|\ge\cdots\ge|\mu_d|$. The highest weight
$\lambda=\omega_a+\omega_p$ has modulus
$\prod_{i\le a}|\mu_i|\prod_{i\le p}|\mu_i|$.
For $i=a,p$, lowering its vector by $\alpha_i$ is nonzero
because $\langle\lambda,\alpha_i^\vee\rangle=1$.
Equality $|\mu_i|=|\mu_{i+1}|$ would therefore give two independent
weight spaces of maximal modulus, contradicting proximality.
Both required spectral gaps follow.

Let $L$ be the Zariski closure of $\tau(\Gamma)$. The finite-index
subgroup $\tau(\Gamma)\cap L^0$ is Zariski dense in $L^0$, and its
line orbits span $W$ by strong irreducibility. Thus avoiding
$\PP(\ker T)$ for each of $\jmath(gF_s^+),\jmath(gF_s^-)$
gives two nonempty Zariski-open conditions on $L^0$; some
$g\in\Gamma$ satisfies both. The conjugates
$s_j=h_jgsg^{-1}h_j^{-1}$ have spectral flags
$h_jgF_s^\pm\to F$.
\end{proof}

\begin{lemma}
\label{lem:spectral-critical-powers}
Let $s\in\mathrm{GL}(V)$ have strict gaps between the moduli of its
$a$th and $(a+1)$st eigenvalues and between those of its $p$th and
$(p+1)$st eigenvalues. For $q\in\{a,p\}$ let $E_q^+(s)$ be
the sum of the generalised eigenspaces corresponding to the $q$
leading eigenvalues, counted with algebraic multiplicity. Then the
top left singular $q$-plane of $s^M$ converges to $E_q^+(s)$ as
$M\to\infty$, and the top left singular $q$-plane of $s^{-M}$
converges to the corresponding spectral plane of the inverse.
These conclusions allow nontrivial Jordan blocks within either
spectral group.
\end{lemma}
\begin{proof}
For $q\in\{a,p\}$, the strict cut makes
$\bigwedge^qE_q^+(s)$ the simple dominant eigenline of
$\bigwedge^qs$. In Jordan form all other terms decay exponentially
relative to this eigenvalue, despite their polynomial nilpotent
factors. Thus every normalised subsequential limit of
$\bigwedge^qs^M$ has rank one with this image, and its
second-to-first singular ratio tends to zero. The dominant left
singular line converges to $\bigwedge^qE_q^+(s)$; the Pl\"ucker
embedding gives the asserted plane. Apply the same argument to
$s^{-1}$ for the inverse flags.
\end{proof}

The regular-seed argument combines compact return on the middle
block with a small correction through the expanding block.
The following two steps allow the extreme spectral spaces to become
nontransverse in the limit. They are unnecessary when $p=a+1$.

\begin{lemma}
\label{lem:compact-middle-return}\label{lem:quant-return}
\begin{enumerate}[label=\textup{(\alph*)},nosep]
\item Let $T\in\mathrm{GL}(E\oplus C\oplus H)$ preserve the complex summands, where
$\overline{\langle[S]\rangle}\subset\PGL(C)$ is compact for
$S=T|_C$. Suppose the eigenvalue moduli on $E$ and $H$ are
strictly below and above the common modulus $q$ on $C$.
For every $\varepsilon>0$ and $M_0$ there are $M>M_0$ and
$\lambda\in\C^*$ with $|\lambda|q^M=1$ and
\[
 \|\lambda T^M|_E\|<\varepsilon,\qquad
 \|(\lambda T^M|_H)^{-1}\|<\varepsilon,\qquad
 \|\lambda S^M-I_C\|<\varepsilon.
\]
\item Let $T_j\in\mathrm{GL}(V)$ carry the summands of
$V=E_j^-\oplus C_j^-\oplus H_j^-$ onto those of
$V=E_j^+\oplus C_j^+\oplus H_j^+$, respectively, and let $\lambda_j\in\C^*$. Suppose
$E_j^-\to A$, $H_j^+\to A$, $C_j^\pm\to C$ and
$A\cap C=0$. For identifications $J_j^\pm:C\to C_j^\pm$
converging to inclusion, assume
\[
 \|\lambda_jT_j|_{E_j^-}\|\to0,\qquad
 \|(\lambda_jT_j|_{H_j^-})^{-1}\|\to0,\qquad
 (J_j^+)^{-1}\lambda_jT_jJ_j^-\to I_C.
\]
Then every $[v]\in\PP(A\oplus C)$ has representatives
$v_j\to v$ with $[T_jv_j]\to[v]$. If the $[T_j]$ are
pairwise distinct elements of a discrete group, this space is recurrent.
\end{enumerate}
\end{lemma}
\begin{proof}
(a) Compactness conjugates the projective cyclic closure into a
projective unitary group. Hence $S$ is diagonalisable with common
eigenvalue modulus $q$, and $q^{-1}S$ is similar to a unitary matrix.
Compact return gives unbounded positive integers $M$ and phases
$\omega_M$ with $\omega_Mq^{-M}S^M\to I_C$.
Take $\lambda_M=\omega_Mq^{-M}$. Jordan form on the extreme
summands gives, for some $K$ and $0<\vartheta<1$,
\[
 \|q^{-M}T^M|_E\|+\|(q^{-M}T^M|_H)^{-1}\|
 \le K M^{\dim V}\vartheta^M
\]
for all large $M$, hence also along these returns. This proves all
three estimates simultaneously, with the stated scalar normalisation.

(b) Write $v=u+c$, and choose $u_j^-\in E_j^-$ and
$u_j^+\in H_j^+$ tending to $u\in A$. Set
\[
 v_j=u_j^-+J_j^-c+(\lambda_jT_j|_{H_j^-})^{-1}u_j^+.
\]
The inverse estimate gives $v_j\to u+c$, and the other two limits
give $\lambda_jT_jv_j\to u+c$. This also covers $u=0$ or $c=0$.
The inverse norm controls the correction even when the angles
between the extreme summands degenerate; no uniform transversality
is required.
\end{proof}

Over $\R$, strict modulus cuts select conjugation-invariant sums
of generalised eigenspaces, so the flags are real. In part~(a),
$q^{-1}S$ is similar to an orthogonal matrix; its powers return to
the identity with real scalars $\lambda_M=q^{-M}$ and the same estimates.

\begin{theorem}\label{thm:crit-recurrence}
Let $\Gamma<\PSL(V)$ be infinite and discrete, of self-opposed uniform critical type $(a,p)$. Assume that its critical module is full and that the action is critically strongly irreducible. Suppose that $p-a=1$, or that $\mathscr F$ contains a regular critical seed. Then
\[
\Lambda_{\mathrm{crit}}=\mathscr R_\Gamma,
\]
and $\Omega_{\mathrm{crit}}$ is the largest domain of proper discontinuity \textup{(}possibly empty\textup{)}.
\end{theorem}
\begin{proof}
Lemma~\ref{lem:rec-contractions} makes $\mathscr F$ minimal.
The three-cell calculation preceding Proposition~\ref{prop:proj-diagonal}
identifies the fat wall $\mathcal T_\star(\mathscr F)$ with
$\Lambda_{\mathrm{crit}}$. We construct a recurrence seed.

\emph{The case $p=a+1$.} For any $F\in\mathscr F$, choose
$s_j$ with both spectral flags tending to $F$ by
Lemma~\ref{lem:rec-contractions}. Lemma~\ref{lem:spectral-critical-powers}
also gives divergent singular gaps at $a,p$ for powers of each
fixed $s_j$, since the corresponding exterior powers have
rank-one normalised limits. Choose $M_j>j$ so that
$\gamma_j=s_j^{M_j}$ are pairwise distinct, both gap ratios
exceed $j$, and the forward and inverse singular critical flags
are within $1/j$ of their spectral flags. Such distinct elements
are available because strict spectral gaps make the projective
powers of each $s_j$ escape compact sets. A Cartan extraction
therefore coalesces at $\Theta$ with flag $F$. Its middle block
is a line, so Proposition~\ref{prop:proj-diagonal} puts all of
$\PP(B)$ on the diagonal, making $F$ a recurrence seed.

\emph{The regular-seed case.} Take $s_j$, $C_j\to C_\infty$
and $F=(A\subset B)$ from Definition~\ref{def:seed}, so
$B=A\oplus C_\infty$. Fix identifications
$J_j:C_\infty\to C_j$ tending to inclusion, and then choose
\[
 \varepsilon_j=\bigl[j(1+\|J_j\|\|J_j^{-1}\|)\bigr]^{-1}.
\]
Apply Lemma~\ref{lem:compact-middle-return}(a) to
$s_j|_{A_j^-\oplus C_j\oplus A_j^+}$ with this tolerance.
Take its lower bound on the power above $j$, large enough for
the forward and inverse singular flags to be within $1/j$ of
the spectral flags (Lemma~\ref{lem:spectral-critical-powers}),
and to avoid all previously chosen projective elements.
Again strict spectral gaps justify this last, separate requirement.
With $q_j$ the common middle eigenvalue modulus, the resulting
$T_j=s_j^{N_j}$ and $\lambda_j$ satisfy $N_j\to\infty$,
$|\lambda_j|q_j^{N_j}=1$, and simultaneously
\[
 J_j^{-1}\lambda_jT_jJ_j\to I_{C_\infty},\qquad
 \|\lambda_jT_j|_{A_j^-}\|\to0,\qquad
 \|(\lambda_jT_j|_{A_j^+})^{-1}\|\to0.
\]
The tolerance was fixed after the norms of $J_j,J_j^{-1}$,
so the middle estimate holds on the fixed space $C_\infty$.
Part~(b), with $E_j^\pm=A_j^-$, $H_j^\pm=A_j^+$ and
$C_j^\pm=C_j$, now puts $\PP(B)$ on the diagonal of this
sequence and of any Cartan extraction. Its singular flags
coalesce at $F$; Proposition~\ref{prop:proj-diagonal} identifies
the resulting middle Levi action as projectively trivial.
Thus $F$ is a recurrence seed.

Uniform critical type makes $\Gamma$ $\Theta$-convergent.
Theorem~\ref{thm:singular-recurrence} now gives both
$\mathscr R_\Gamma=\Lambda_{\mathrm{crit}}$ and the largest-domain
assertion.
\end{proof}

\begin{remark}
For $p-a=1$ and lattices, fullness and strong irreducibility follow from Borel density. If $\Gamma$ is the image of a lattice in a connected semisimple group $S$ without compact factors, and $V$ and $W_{a,p}$ are irreducible $S$-modules, then the nonzero critical module is $S$-invariant by Furstenberg's form of Borel density~\cite{r3bor:Borel1960,r3bor:Furstenberg1976}, hence equals $W_{a,p}$. Applying the same argument to finite-index subgroups gives critical strong irreducibility. The uniform critical-type hypothesis remains separate: it does not follow from Borel density or irreducibility.
\end{remark}

\begin{proposition}\label{prop:wide-seed}
For every $d\ge4$ there is an infinite discrete, complex-Zariski-dense subgroup $\Gamma<\PSL_d(\C)$ contained in the projective image of $\mathrm{SU}(d-1,1)$, of uniform critical type $(1,d-1)$ with a regular critical seed. Its largest domain is nonempty.
\end{proposition}
\begin{proof}
Write $V=\C e_1\oplus C\oplus\C e_d$, where $C=\langle e_2,\ldots,e_{d-1}\rangle$, and equip $V$ with the Hermitian form
\[
J=\begin{pmatrix}0&0&1\\0&I_{d-2}&0\\1&0&0\end{pmatrix}.
\]
On $\langle e_1,e_d\rangle$, in that order, let
\[
u=\begin{pmatrix}1&i\\0&1\end{pmatrix},\qquad
s=\begin{pmatrix}2&i\\-3i&2\end{pmatrix},\qquad
w=\begin{pmatrix}0&i\\i&0\end{pmatrix},
\]
extending each by the identity on $C$. For $z\in\{2e_k,2ie_k:2\le k\le d-1\}$ also use
\[
n(z)=\begin{pmatrix}1&-z^*&-\|z\|^2/2\\0&I_{d-2}&z\\0&0&1\end{pmatrix}.
\]
These matrices preserve $J$, have determinant one, and have Gaussian-integer entries. Let $\Gamma$ be their projective group. It is infinite since $u$ has infinite order, and discrete since bounded Gaussian-integer matrices form a finite set and the projective scalar ambiguity is finite.

The logarithms of the two unipotents $n(2e_k)$ and $n(2ie_k)$
are respectively $-2E_{1k}+2E_{kd}$ and
$2iE_{1k}+2iE_{kd}$. Their complex span therefore contains
$E_{1k}$ and $E_{kd}$ in the Lie algebra of the Zariski closure.
Conjugation by $w$ gives $E_{dk}$ and $E_{k1}$, and brackets
then give every off-diagonal matrix unit and all traceless
diagonal matrices. Thus the closure is $\mathrm{SL}_d(\C)$. Since $\mathrm{SL}_d(\C)$ is connected, every finite-index subgroup has the same Zariski closure; it acts irreducibly on $V$ and on $W_{1,d-1}\simeq\mathfrak{sl}_d(\C)$.

Put $v=(e_1-e_d)/\sqrt2$, so $J=I-2vv^*$. For $g\in\Gamma$, the identity $g^*Jg=J$ gives
\[
g^*g-I=2(g^*vv^*g-vv^*),\qquad \operatorname{rank}(g^*g-I)\le2.
\]
Since $\det(g^*g)=1$, the singular values are $e^t,1,\ldots,1,e^{-t}$. Along every escaping sequence $t\to\infty$, so the only diverging gaps are at $1$ and $d-1$: the type is uniformly $(1,d-1)$.

Finally, $s$ has eigenvalues $2+\sqrt3$, $1$ on $C$, and $2-\sqrt3$. For $s_j=u^jsu^{-j}$, the two extremal eigenlines are
\[
\C\bigl((1\pm j\sqrt3)e_1\mp i\sqrt3e_d\bigr),
\]
and both tend to $A=\C e_1$. The middle block is the fixed space $C$, on which every $s_j$ acts identically. Hence the attracting and repelling critical flags tend to $F=(A\subset B)$, where $B=A\oplus C$. Lemma~\ref{lem:spectral-critical-powers} and a diagonal choice of powers put $F$ in the critical limit-flag set, so it is a regular critical seed. The critical module is then nonzero and invariant, hence full by the irreducibility already proved.

The $J$-negative projective lines form a nonempty complex hyperbolic ball on which $\mathrm{SU}(d-1,1)$, and hence $\Gamma$, acts properly. Theorem~\ref{thm:crit-recurrence} shows that the largest domain of $\Gamma$ contains this ball.
\end{proof}

\begin{remark}
For non-elementary Borel--Anosov groups, Stecker~\cite{SteckerBalanced} classifies maximal domains by minimal fat ideals. Kapovich--Leeb--Porti~\cite{KLP} prove properness for fat-thickening complements and cocompactness for balanced ones; Guichard--Wienhard~\cite{GuichardWienhard} construct proper cocompact domains for Anosov representations. Here largestness follows from moving-source recurrence, through critical seeds, the cusped relation below, or divisors. For $\C\PP^2$, see also~\cite{BCNSFlag2026,CCL2026}.
\end{remark}

\subsection{Cusped Hitchin groups}\label{sec:cusped}
Canary--Zhang--Zimmer supply diverging singular gaps and limiting
flags even along parabolic powers. These determine the complete
relation below without traditional Borel--Anosov or Zariski-density
assumptions.

\begin{lemma}\label{lem:regular-rectangles}
Let $(g_m)$ be a sequence of complex invertible $d\times d$ matrices
for which every adjacent singular ratio diverges. Suppose the
limiting left singular full flag is $(G_j)$ and the limiting
slow-right full flag is $(F_j)$, with $\dim G_j=j$ and
$\dim F_j=d-j+1$. Then
\[
\D_{(g_m)}=\bigcup_{j=1}^d\PP(F_j)\times\PP(G_j).
\]
In particular, the formula remains valid when the two full flags
have coinciding distinguished lines.
\end{lemma}
\begin{proof}
Every subsequence admits an extraction with converging singular
frames. All blocks are lines, so Proposition~\ref{prop:proj} gives
$\D(x)=\PP(G_i)$ for
$x\in\PP(F_i)\setminus\PP(F_{i+1})$; the union of these fibres
is the displayed rectangle relation. Theorem~\ref{thm:dest}(d)
gives convergence of the extracted graphs to it. This limit depends
only on the two prescribed flags, regardless of any coincidences.
Compactness of the hyperspace and uniqueness of every extracted
limit therefore give convergence of the original graphs. Each
pair is consequently approximated on those graphs along the whole
original sequence, as required for moving sources.
\end{proof}

\begin{theorem}\label{thm:cusped}
Let $\Pi<\mathrm{PSL}(2,\R)$ be torsion-free, non-elementary and
geometrically finite. Let $d\ge2$ and let
$\rho:\Pi\to\mathrm{SL}_d(\R)$ be a Hitchin representation in the
sense of Canary--Zhang--Zimmer~\cite{r11CZZCusped}.
Write $\xi:\Lambda(\Pi)\to\mathrm{Fl}(\R^d)$ for its positive
flag map, $\xi^q_\C(x)=\xi^q(x)\otimes_\R\C$ and
$\Gamma=[\rho(\Pi)]<\mathrm{PSL}_d(\C)$. Put
\[
 U_q=\bigcup_{z\in\Lambda(\Pi)}\PP(\xi^q_\C(z))
 \qquad(1\le q\le d),\qquad c=\lceil d/2\rceil.
\]
Then
\begin{equation}\label{eq:cusped-relation}
 \D_\Gamma=\bigcup_{j=1}^d U_{d-j+1}\times U_j.
\end{equation}
For $x\in\C\PP^{d-1}$ with least index $r$ such that $x\in U_r$,
its destination is $\D_\Gamma(x)=U_{d-r+1}$. The sets
$U_1,\ldots,U_d$ form a strictly increasing chain, and each occurs
as a destination fibre. Moreover, $q_\Gamma=c-1$ and
$\mathscr R_\Gamma=\Lambda_{\mathrm{crit}}(\Gamma)=U_c$.
The complement of $U_c$ is nonempty and is the largest domain of
proper discontinuity. If $\Pi$ contains a parabolic element,
$\rho$ is not Borel--Anosov in the traditional word-hyperbolic sense.
\end{theorem}

\begin{proof}
Canary--Zhang--Zimmer prove that $\rho$ is $P_k$--Anosov in their
sense for every $1\le k<d$, with boundary maps given by the
$k$-planes of the positive flag map. Hyperbolic elements have
images with distinct eigenvalue moduli, and parabolic elements
have regular projectively unipotent images
\cite[Theorem~1.4]{r11CZZCusped}. For a sequence $\delta_j$
escaping in $\Pi$, one has
$d_{\mathbb H^2}(o,\delta_j o)\to\infty$. Their gap estimate and
Cartan-flag convergence give
\[
 \frac{\sigma_k(\rho(\delta_j))}
 {\sigma_{k+1}(\rho(\delta_j))}\longrightarrow\infty
 \quad (1\le k<d),
\]
and if $\delta_j o\to x\in\Lambda(\Pi)$, the target
$k$-planes converge to $\xi^k(x)$
\cite[Theorem~1.1(1) and Corollary~4.2(1)]{r11CZZCusped}.
The singular ratios are unchanged by complexification.
A compact set of projective matrices has bounded singular ratios,
so its preimage in $\Pi$ is finite. Thus $\Gamma$ is discrete and
the projective map has finite kernel, which is trivial since
$\Pi$ is torsion-free.

For any escaping $\delta_m$, extract so that
$\delta_mo\to z_+$ and $\delta_m^{-1}o\to z_-$.
The resulting full-type sequence has target and slow-source flags
\[
 G_j=\xi^j_\C(z_+),\qquad
 F_j=\xi^{d-j+1}_\C(z_-)\qquad(1\le j\le d).
\]
Every boundary point is approached by an escaping orbit sequence,
so the full flag limit set is exactly
$\xi_\C(\Lambda(\Pi))$. It is minimal by equivariance and
minimality of the classical limit-set action of a non-elementary
Fuchsian group. Thus all hypotheses of Lemma~\ref{lem:pairs}
and Theorem~\ref{thm:pair} hold: the limiting flag pairs fill
the square of this set. Lemma~\ref{lem:regular-rectangles}
expresses each pair's relation as
$\bigcup_j\PP(\xi^{d-j+1}_\C(z_-))\times
\PP(\xi^j_\C(z_+))$, including $z_-=z_+$.
Their union is~\eqref{eq:cusped-relation}.

The sets $U_q$ are nested and closed by compact incidence
over $\Lambda(\Pi)$. If $r$ is the least index with $x\in U_r$,
the products containing $x$ in the first factor are precisely those
with $j\le d-r+1$, which gives the stated fibre.

Fix $z\in\Lambda(\Pi)$ and $1\le j<d$. Choose real vectors
$0\ne u\in\xi^1(z)$ and $v\in\xi^{j+1}(z)\setminus\xi^j(z)$.
Then $[u+iv]\in U_{j+1}$. If it lay in $\PP(\xi^j_\C(w))$,
reality would give $u,v\in\xi^j(w)\subseteq\xi^{d-1}(w)$.
Transversality of $\xi^1(z)$ to $\xi^{d-1}(w)$ for $w\ne z$
would force $w=z$, contradicting the choice of $v$. Thus
$[u+iv]\in U_{j+1}\setminus U_j$. Since $U_1$ is nonempty, every
least index occurs, so the fibre formula yields exactly $d$
distinct fibres.

The diagonal of the $j$th product is $U_{\min(j,d-j+1)}$,
so $\mathscr R_\Gamma=U_c$. Each product has a factor in $U_c$,
making $U_c$ a vertex cover. Theorem~\ref{thm:relational} gives
the largest-domain assertion.

Both orientations are admissible: all blocks are lines and
$d-1\ge c$. Their first crossing is $r_c=c$, so $q_\Gamma=c-1$;
the selected dimensions are $c,d-c$. The boundary realisation
already proved gives $\Lambda_{\mathrm{crit}}=U_c\cup U_{d-c}=U_c$.

To see that this domain is nonempty, choose distinct
$x,y\in\Lambda(\Pi)$ and nonzero real vectors $u\in\xi^1(x)$,
$v\in\xi^1(y)$. If $u+iv\in\xi^c_\C(z)$, then reality gives
$u,v\in\xi^c(z)\subseteq\xi^{d-1}(z)$. Transversality forces
$x=z=y$, which is impossible. Hence $[u+iv]$ is in the domain.

Finally, a parabolic element has infinite order and its image has
all eigenvalues of the same modulus, whereas the image of every
infinite-order element of a traditional Borel--Anosov representation
has distinct eigenvalue moduli.
\end{proof}

\begin{proposition}\label{prop:cusped-existence}
Arbitrarily close to the principal Fuchsian representation of the
three-punctured sphere group into $\mathrm{SL}_5(\R)$, there are
real Zariski-dense cusped Hitchin representations with all three
peripheral images regular unipotent. In particular,
Theorem~\ref{thm:cusped} applies beyond the principal and orthogonal
loci.
\end{proposition}
\begin{proof}
Put $G=\mathrm{SL}_5(\R)$ and let $\mathcal U$ be its regular
unipotent conjugacy class. The principal representation on
$\operatorname{Sym}^4\R^2$ of a complete finite-area structure on
$S_{0,3}$ gives a triple $(u_1,u_2,u_3)\in\mathcal U^3$ with
$u_1u_2u_3=I$. A regular unipotent has centraliser dimension $4$
in $\mathfrak{sl}_5$, so $\dim G=24$ and $\dim\mathcal U=20$.
For the multiplication map $\mu:\mathcal U^3\to G$, the
annihilator of $d\mu$ under the trace pairing is the common
Lie-algebra centraliser of the principal triple, after translating
the tangent spaces. It is zero by absolute irreducibility of
$\operatorname{Sym}^4\R^2$. Thus $\mu$ is a submersion there and
$M=\mu^{-1}(I)$ is locally smooth of dimension $3\cdot20-24=36$.
Projection to the first two factors identifies this relative family
with a local $36$-dimensional submanifold of $G^2$, since the third
factor is the inverse of their product.

For positivity, the unique invariant complete flag of $u\in\mathcal U$ has
$j$-plane $\ker((u-I)^j)$, $1\le j<5$, and varies continuously
by constancy of rank. Thus the triples carry continuous peripheral
framings. Fix an ideal triangulation and the orientation data for
its Fock--Goncharov chart, applied to the projectivised framed
representation in $\mathrm{PGL}_5(\R)$. Its coordinates are locally
rational functions of finitely many flags and monodromy words.
Their denominators are nonzero at the principal representation,
where all coordinates are strictly positive
\cite[Theorem~6.2 and Corollary~7.2]{r11FGPositive}.
Finiteness gives one relative neighbourhood in $M$ on which all
these denominators stay nonzero and all coordinates stay positive.
By \cite[Theorem~1.6]{r11FGPositive}, they give a positive
equivariant framing on the Farey set, which here is the set of
peripheral fixed points. Every parabolic in the surface group is
conjugate to a nonzero power of a peripheral generator; its image
is therefore still regular unipotent. These representations are
positive type preserving, hence cusped Hitchin by
\cite[Definition~9.1 and Theorem~9.2]{r11CZZCusped}.

The source is a torsion-free Fuchsian lattice, so the Zariski-closure
classification \cite[Theorem~11.10]{r11CZZEntropy} applies:
in degree five a proper projective closure is principal
$\mathrm{PSL}_2$ or conjugate to $\mathrm{PSO}(3,2)$.
Both preserve a nondegenerate symmetric form of signature $(3,2)$.
Let $\mathcal Q$ be the open set of their projective classes,
of dimension $14$. If $A\in G$ preserves $[Q]$, then
$A^tQA=cQ$ and determinants give $c^5=1$, hence $c=1$.
The stabiliser of $[Q]$ therefore has dimension $10$.
Consequently the semialgebraic incidence
\[
 \mathcal I=\{(A,B,[Q])\in G^2\times\mathcal Q:
 A^tQA=Q,\ B^tQB=Q\}
\]
has dimension $14+2\cdot10=34$, and its projection to $G^2$
has dimension at most $34$. It cannot contain an open subset
of the local $36$-dimensional relative family $M$.
Every sufficiently small neighbourhood of the principal triple
in $M$ therefore contains a positive relative representation
preserving no such form. The classification forces full projective
Zariski closure, and hence full real Zariski closure $G$.
\end{proof}

\subsection{Symmetric powers}\label{sec:symmetric} Let $\rho_N:\PSL(2,\C)\to\PGL(\Sym^N\C^2)$ be the irreducible representation of degree $N\ge1$, so that $\PP(\Sym^N\C^2)=\Sym^N\C\PP^1$ is the space of effective divisors of degree $N$ on $\C\PP^1$.

Let $\Gamma<\PSL(2,\C)$ be non-elementary and discrete, with limit set $\Lambda$, and put $G_N=\rho_N(\Gamma)$ and $t=\lceil N/2\rceil$. Use the Hermitian metric on $\Sym^N\C^2$ induced from the tensor power. If a determinant-one lift of $g$ has singular values $s,s^{-1}$ with $s\ge1$, its symmetric power has singular values $s^{N-2i}$, $0\le i\le N$, as seen in the orthonormal basis of symmetrised monomials. Thus its extreme singular ratio is $s^{2N}$ and every adjacent singular ratio is $s^2$. Along an escaping sequence $s\to\infty$, so $G_N$ is discrete and every divergent sequence has full type.

\begin{lemma}\label{lem:coalesce}
For every $\xi\in\Lambda$ there are loxodromic $s_j\in\Gamma$ whose attracting and repelling fixed points both tend to $\xi$.
\end{lemma}
\begin{proof}
Take distinct $h_j$ converging to $\xi$ uniformly on compact subsets of $\C\PP^1\setminus\{\eta\}$. Since the orbit of $\eta$ is infinite, some conjugate $\ell$ of a loxodromic element has both fixed points different from $\eta$, and $s_j=h_j\ell h_j^{-1}$ has fixed points $h_j\ell^\pm\to\xi$.
\end{proof}

\begin{theorem}\label{thm:C}
\begin{enumerate}[label=\textup{(\alph*)},nosep]
\item $\Lambda_{\mathrm{crit}}(G_N)=\mathscr R_{G_N}=C_\Lambda:=\{D:\ \mathrm{mult}_\xi D\ge t\text{ for some }\xi\in\Lambda\}$. Hence $\Omega_{\max}=\{D:\mathrm{mult}_\xi D<t\ \forall\xi\in\Lambda\}$ is the largest domain of proper discontinuity. No finite generation, convex cocompactness or absence of parabolics is assumed.
\item If $\Lambda\neq\C\PP^1$, then $\Lambda_{\Kul}(G_N)=\Lambda_{\Myr}(G_N)=\{D:\ \mathrm{supp}\,D\cap\Lambda\neq\varnothing\}$, and $\Omega_{\Kul}=\Eq=\Sym^N(\C\PP^1\setminus\Lambda)$.
\item If $\Lambda\neq\C\PP^1$, then $\Eq=\Omega_{\Kul}=\Omega_{\max}$ for $N\le2$, while $\Eq=\Omega_{\Kul}\subsetneq\Omega_{\max}$ for $N\ge3$.
\end{enumerate}
\end{theorem}
\begin{proof}
(a) \emph{Properness.} If $D_j\to D$ and $\gamma_jD_j\to D'$
for distinct $\gamma_j$, extract a rank-one limit with attractor
$a$ and repeller $b$ in $\Lambda$. For $k=\mathrm{mult}_bD$,
at least $N-k$ roots stay away from $b$, so
$\mathrm{mult}_aD'\ge N-k$. Hence $D\notin C_\Lambda$ implies
$k<t$ and $D'\in C_\Lambda$. The set $C_\Lambda$ is compact
as the image of $(\xi,E)\mapsto t\xi+E$ on
$\Lambda\times\Sym^{N-t}\C\PP^1$, and invariant;
Theorem~\ref{thm:relational}(i) applies.

\emph{Recurrence.} Let $D=k\xi+E$ with $\xi\in\Lambda$,
$k\ge t$ and $\xi\notin\mathrm{supp}\,E$. Take $s_j$ from
Lemma~\ref{lem:coalesce}, with fixed points $a_j,b_j\to\xi$
off $\mathrm{supp}\,E$, and choose distinct powers
$\gamma_j=s_j^{M_j}$ such that $\gamma_jE$ is within $1/j$
of $(N-k)a_j$ and $\gamma_j^{-1}E$ is within $1/j$ of
$(N-k)b_j$. Place the remaining $2k-N\ge0$ roots at the fixed point $a_j$:
\[
D_j=E+\gamma_j^{-1}E+(2k-N)a_j\longrightarrow D,\qquad\gamma_jD_j=\gamma_jE+E+(2k-N)a_j\longrightarrow D.
\]
So $C_\Lambda\subseteq\mathscr R_{G_N}$, and with properness on the complement, $C_\Lambda=\mathscr R_{G_N}$.

\emph{The wall.} The target flag of a sequence with attractor $a$ is the osculating flag of the rational normal curve, $\PP(\mathcal O_k(a))=\{D:\mathrm{mult}_aD\ge N-k\}$. Here $c=\lceil(N+1)/2\rceil$, every sequence is admissible, $B_h=\mathcal O_{c-1}(a)$ gives $\{\mathrm{mult}_a\ge t\}$, and $B^*_h$ is the space $\{\mathrm{mult}_b\ge N-t+1\}$ at the repeller, contained in $C_\Lambda$. Attractors are dense in $\Lambda$, so $\Lambda_{\mathrm{crit}}=C_\Lambda$.

(b) The first limit of a sequence with attractor $a$ and repeller $b$ has kernel $H_b=\{D:b\in\mathrm{supp}\,D\}$, and repellers fill $\Lambda$; Proposition~\ref{prop:kernel} gives $\Lambda_{\Myr}$ and $\Eq$.

For $\Lambda_{\Kul}$, first $B_{G_N}\subseteq\Sym^N\Lambda$. An infinite stabiliser has an infinite finite-index subgroup fixing each support point. Every such point belongs to $\Lambda$, by applying rank-one convergence to distinct elements of this subgroup. For an escaping sequence with attractor $a$ and repeller $b$, every root of a fixed divisor other than $b$ tends to $a$, and all accumulation points of the images of $b$ belong to the closed invariant set $\Lambda$. This also puts the limits defining $L_1$ in $\Sym^N\Lambda$. One inclusion is Proposition~\ref{prop:kulkarni-myrberg}.

Conversely let $D=k\xi+E$ with $\xi\in\Lambda$, $k\ge1$, $\xi\notin\mathrm{supp}\,E$, and fix $z_0\notin\Lambda$. With $s_j$ from Lemma~\ref{lem:coalesce}, choose pairwise distinct powers $\gamma_j$ with $\gamma_jz_0\to\xi$ and $\gamma_j^{-1}E\to(N-k)\xi$. The sources $D_j=kz_0+\gamma_j^{-1}E$ form, with their limit $kz_0+(N-k)\xi$, a compact set disjoint from $B_{G_N}$ (each has a root at $z_0$), and $\gamma_jD_j\to D$. So $D\in L_2$.

(c) For $N\le2$, $t=1$. For $N\ge3$ the divisor $\xi+(N-1)z_0$, $\xi\in\Lambda$, $z_0\notin\Lambda$, lies in $\Omega_{\max}\setminus\Omega_{\Kul}$.
\end{proof}

For $\Lambda=\C\PP^1$ the complement of $C_\Lambda$ is the set of GIT-stable divisors~\cite[Theorem~3.10]{ThomasGIT}; for proper $\Lambda$ stability is imposed only over $\Lambda$. The regions already differ dynamically: the largest domain is governed by a multiplicity threshold, while equicontinuity excludes any root on the classical limit set.

\appendix
\section{The cell lemma}\label{sec:proofcells}
We prove Lemma~\ref{lem:cells} with the conventions of
Section~\ref{sec:setting-sum}, using real Bruhat and parabolic structure
\cite[4.8, 4.10(d), 5.15--5.18 and~5.20]{BorelTits}
\cite[Chapter~VII]{KnappLG}. The root-coordinate argument includes
non-reduced restricted systems.

\begin{description}[leftmargin=2em,nosep]
\item[(B1) Root groups] Call $\Psi\subset\Phi$ \emph{closed} if
$\lambda,\mu\in\Psi$ and $\lambda+\mu\in\Phi$ imply
$\lambda+\mu\in\Psi$ (take $\lambda=\mu$ to see that
$\lambda\in\Psi,2\lambda\in\Phi$ implies $2\lambda\in\Psi$).
If moreover $\Psi$ is contained in a positive system, then
$\bigoplus_{\Psi}\g_\lambda$ is a nilpotent subalgebra and
$U_\Psi=\exp\bigoplus_\Psi\g_\lambda$ is a closed unipotent
subgroup. For \emph{any} ordering $\lambda_1,\dots,\lambda_r$
of $\Psi$, the product of exponential slices
$\exp\g_{\lambda_1}\times\cdots\times\exp\g_{\lambda_r}\to U_\Psi$
is a diffeomorphism.

To prove the last statement, choose a positive root system
containing $\Psi$ and filter by root height. The BCH formula is
triangular: the coordinate in each root space is its input
coordinate plus a polynomial in coordinates of strictly smaller
height. Solve successively by height; nilpotence makes the inverse
polynomial and global. This argument works for any ordering of
the slices.

In a non-reduced system $\exp\g_\alpha$ need not itself be a
subgroup when $2\alpha$ is a root; its brackets are included at
the $2\alpha$ level of the filtration.
\item[(B2) Big cell] For a standard parabolic $Q=P_\theta$, the product $N_\theta^-\times Q\to G$ is a diffeomorphism onto an open subset. Equivalently $n\mapsto nQ$ is a diffeomorphism of $N_\theta^-$ onto an open neighbourhood of the base point of $G/Q$.
\end{description}

Write $W_\Theta\subset W$ for the Weyl group of the Levi factor $L_\Theta$, generated by the reflections in the roots of $\Delta\setminus\Theta$; similarly $W_\theta$. Representatives of $W$ in $N_K(\ak)$ will be denoted by the same letters, and we choose each double-coset representative $\sigma$ to be of minimal length in $W_\Theta\sigma W_\theta$. In particular $\sigma^{-1}\lambda>0$ for every positive root $\lambda$ in the span of $\Delta\setminus\Theta$; the analogous condition on the right holds for the Levi roots of $Q$.

\paragraph{Step 1: part (i).}
The parabolic Bruhat decomposition $G=\bigsqcup_{w\in W_I\backslash W/W_J}P_IwP_J$ for standard parabolics~\cite[5.15--5.20]{BorelTits} applies to $P_{\Theta^*}$ and $Q$, where $\Theta^*=-w_0\Theta$. Since $P^-_\Theta=w_0P_{\Theta^*}w_0^{-1}$ and $w_0^{-1}W_\Theta w_0=W_{\Theta^*}$, translating by $w_0$ gives $G=\bigsqcup_\sigma P^-_\Theta\sigma Q$ over $\sigma\in W_\Theta\backslash W/W_\theta$, and dividing by $Q$ gives (i).

\paragraph{Step 2: part (ii).}
Since $\sigma$ normalises $\ak$, $A=\sigma A\sigma^{-1}\subset\sigma Q\sigma^{-1}$, so $\sigma Q$ is fixed by $A_\Theta\subset A$; and $L_\Theta$ centralises $A_\Theta$, so every point of $L_\Theta\sigma Q$ is fixed. Conversely, let $x\in X^{A_\Theta}$. By the parabolic Bruhat decomposition $G=\bigsqcup_\sigma P_\Theta\sigma Q$ (Step~1 without the translation by $w_0$), $x=nl\sigma Q$ with $n\in N_\Theta$, $l\in L_\Theta$. Choose $z_t\in A_\Theta$ with $\beta(z_t)\to\infty$ for $\beta\in\Theta$ (as $t\to\infty$). Then $z_t^{-1}x=x$, and also $z_t^{-1}x=(z_t^{-1}nz_t)\,l\sigma Q\to l\sigma Q$, because $\mathrm{Ad}(z_t^{-1})$ multiplies $\g_\lambda$ by $e^{-\lambda(\log z_t)}\to0$ for $\lambda\in\Phi(N_\Theta)$. Hence $x=l\sigma Q$, and the fixed set is $\bigcup_\sigma L_\Theta\sigma Q$; the pieces are disjoint because they lie in distinct $P_\Theta$-double cosets.

To prove compactness of $L_\Theta\sigma Q/Q$, use the Iwasawa decomposition of the reductive group $L_\Theta$: $L_\Theta=K_\Theta\,A\,N_L$ with $K_\Theta=K\cap L_\Theta$ and $N_L=N\cap L_\Theta$. We have $A\subset\sigma Q\sigma^{-1}$, and $N_L\subset\sigma Q\sigma^{-1}$ because for a positive root $\lambda$ of $L_\Theta$ our choice of $\sigma$ gives $\sigma^{-1}\lambda>0$, a root of $Q$. Hence $L_\Theta\sigma Q=K_\Theta\sigma Q$ is the image of the compact group $K_\Theta$.

\paragraph{Step 3: part (iii).}
By Step~1 and $P_\Theta^-=N_\Theta^-L_\Theta$, a point of $O_\sigma$ is $x=nl\sigma Q$ with $n\in N_\Theta^-$, $l\in L_\Theta$. For $z\in A_\Theta$ with $\beta(z)\to\infty$ ($\beta\in\Theta$), every root $\lambda$ of $N_\Theta^-$ has a non-positive expansion $\lambda=-\sum n_\alpha\alpha$ with $n_\beta>0$ for some $\beta\in\Theta$, so $\lambda(\log z)=-\sum_{\beta\in\Theta}n_\beta\beta(\log z)\to-\infty$, and $znz^{-1}\to1$ (uniformly for $n$ in compact sets). Since $z$ commutes with $l$ and fixes $\sigma Q$,
\[
zx=(znz^{-1})\,l\,\sigma Q\longrightarrow l\sigma Q.
\]
The limit is intrinsic to $x$, hence $\rho(x)=l\sigma Q$ is well defined and independent of how $z\to\infty$.

\paragraph{Step 4: part (iv).}
Put $n_\sigma=\mathrm{Ad}(\sigma)\,\mathfrak n_\theta^-
=\bigoplus_{\lambda\in\Sigma}\g_\lambda$ with
$\Sigma=\sigma\Phi(N_\theta^-)$. The set
$\Phi(N_\theta^-)$ is closed (coefficients add), hence so is
$\Sigma$. Moreover, $\Sigma\subset\sigma\Phi^-$, a positive
system, so (B1) applies to $\Sigma$ and to its closed subsets.
By (B2) translated by $\sigma$, the map $n\mapsto n\sigma Q$
is a diffeomorphism of
$N_\sigma=\sigma N_\theta^-\sigma^{-1}$ onto an open
neighbourhood of $\sigma Q$. It is $A$-equivariant for the
conjugation action on $N_\sigma$: for $a\in A$,
$a\,n\sigma Q=(ana^{-1})\,a\sigma Q=(ana^{-1})\sigma Q$.

Define the $\Theta$-height of a root $\lambda=\sum_\alpha n_\alpha\alpha$ by $h(\lambda)=\sum_{\beta\in\Theta}n_\beta$. It is additive, $h>0$ exactly on $\Phi(N_\Theta)$, $h=0$ on $\Phi(L_\Theta)$, $h<0$ on $\Phi(N_\Theta^-)$. Then $\Sigma_+=\{\lambda\in\Sigma:h>0\}$ and $\Sigma_{\le0}=\{h\le0\}\cap\Sigma$ are closed (the height of a sum is the sum of heights), so (B1) applies to both, and
\[
N_\sigma=\prod_{\Sigma_{\le0}}\exp\g_\lambda\cdot\prod_{\Sigma_+}\exp\g_\lambda=M_-M_+
\]
with unique decomposition, taking in (B1) the ordering that lists $\Sigma_{\le0}$ first.

It remains to show $N_\Theta\sigma Q=M_+\sigma Q$. The set $\Phi(N_\Theta)$ is closed with $\Phi(N_\Theta)\cap-\Phi(N_\Theta)=\varnothing$. Split it as $\Sigma_+\sqcup\Psi$ with $\Psi=\{\lambda\in\Phi(N_\Theta):\sigma^{-1}\lambda\in\Phi(Q)\}$. Indeed, if $\lambda\in\Phi(N_\Theta)$ and $\sigma^{-1}\lambda\notin\Phi(Q)$ then $\sigma^{-1}\lambda\in\Phi(N_\theta^-)$, so $\lambda\in\Sigma\cap\Phi(N_\Theta)=\Sigma_+$. The exponential slices of the root spaces in $\Psi$ lie in $\sigma Q\sigma^{-1}$. By (B1), listing $\Sigma_+$ first, $N_\Theta=M_+\cdot\prod_\Psi\exp\g_\mu$, and the second factor fixes $\sigma Q$. Therefore $N_\Theta\sigma Q=M_+\sigma Q$.

\paragraph{Dimensions.}
We now verify the dimension counts used in Corollary~\ref{cor:conservation}
and Theorem~\ref{thm:dest}(c), taking all dimensions over $\R$.
The orbit map of $M_+$ at
$\sigma Q$ is injective (since $M_+\subset N_\sigma$ and
$N_\sigma\cap\sigma Q\sigma^{-1}=1$) and an immersion, so
$\dim M_+\sigma Q=\sum_{\Sigma_+}\dim\g_\lambda$.

For the cell, $O_\sigma$ is a homogeneous space of $P_\Theta^-$,
with tangent space at $\sigma Q$
$\mathfrak p_\Theta^-/(\mathfrak p_\Theta^-\cap
\mathrm{Ad}(\sigma)\mathfrak q)$.
The centraliser of $\ak$ lies in both Lie algebras and contributes
nothing to this quotient. The remaining summands are the root
spaces $\g_\lambda$ with $h(\lambda)\le0$ and
$\sigma^{-1}\lambda\notin\Phi(Q)$, that is, with
$\lambda\in\Sigma_0\cup\Sigma_-$. So
\[
\dim O_\sigma=\sum_{\Sigma_0\cup\Sigma_-}\dim\g_\lambda,\qquad \dim X=\sum_\Sigma\dim\g_\lambda=\dim O_\sigma+\dim M_+\sigma Q.
\]
This gives the equality in Corollary~\ref{cor:conservation}.
\qed

%--------------------------------------------------------------


\medskip


\begin{thebibliography}{99}
\bibitem{AbelsMargulisSoifer}
H.~Abels, G.~A.~Margulis and G.~A.~Soifer,
\emph{Semigroups containing proximal linear maps},
Israel J. Math. \textbf{91} (1995), no.~1--3, 1--30.
\url{https://doi.org/10.1007/BF02761637}.

\bibitem{BCNLimit2011}
W.~Barrera, A.~Cano and J.~P.~Navarrete,
\emph{The limit set of discrete subgroups of $\mathrm{PSL}(3,\mathbb C)$},
Math. Proc. Cambridge Philos. Soc. \textbf{150} (2011), no.~1, 129--146.
\url{https://doi.org/10.1017/S0305004110000423}.

\bibitem{BCNSFlag2026}
W.~Barrera, A.~Cano, J.~P.~Navarrete and J.~Seade,
\emph{Asymptotic flag geometry of complex Kleinian groups in $\mathbb P^2_{\mathbb C}$},
preprint, 2026, arXiv:2608.22641.
\url{https://arxiv.org/abs/2608.22641}.

\bibitem{Benoist97}
Y.~Benoist,
Propri\'et\'es asymptotiques des groupes lin\'eaires,
\emph{Geom. Funct. Anal.} \textbf{7} (1997), no.~1, 1--47.
\url{https://doi.org/10.1007/PL00001613}.

\bibitem{BenoistQuint}
Y.~Benoist and J.-F.~Quint,
\emph{Random Walks on Reductive Groups},
Ergebnisse der Mathematik und ihrer Grenzgebiete. 3.~Folge, vol.~62,
Springer, Cham, 2016.
\url{https://doi.org/10.1007/978-3-319-47721-3}.
Author version: \url{https://www.math.u-bordeaux.fr/~jquint/publications/walk.pdf}.

\bibitem{BB1973}
A.~Bia\l ynicki-Birula,
\emph{Some theorems on actions of algebraic groups},
Ann. of Math. (2) \textbf{98} (1973), no.~3, 480--497.
\url{https://doi.org/10.2307/1970915}.

\bibitem{BjornerBrenti}
A.~Bj\"orner and F.~Brenti, \emph{Combinatorics of Coxeter Groups},
Graduate Texts in Math. \textbf{231}, Springer, New York, 2005.
\url{https://doi.org/10.1007/3-540-27596-7}.

\bibitem{r3bor:Borel1960}
A.~Borel,
\emph{Density properties for certain subgroups of semi-simple groups
without compact components},
Ann. of Math. (2) \textbf{72} (1960), no.~1, 179--188.
\url{https://doi.org/10.2307/1970150}.

\bibitem{BorelTits}
A.~Borel and J.~Tits,
\emph{Groupes r\'eductifs},
Publ. Math. Inst. Hautes \'Etudes Sci. \textbf{27} (1965), 55--151.
\url{https://doi.org/10.1007/BF02684375}.

\bibitem{BG}
E.~Breuillard and T.~Gelander,
\emph{On dense free subgroups of Lie groups},
J. Algebra \textbf{261} (2003), no.~2, 448--467.
\url{https://doi.org/10.1016/S0021-8693(02)00675-0}.

\bibitem{BGErratum}
E.~Breuillard and T.~Gelander,
\emph{On dense free subgroups of Lie groups -- revisited},
erratum, J. Algebra \textbf{711} (2027), 1238--1243 (forthcoming).
\url{https://doi.org/10.1016/j.jalgebra.2026.10.011}.
Author version: \url{https://arxiv.org/abs/2605.20568}.

\bibitem{r2grBrion2003}
M.~Brion, \emph{Group completions via Hilbert schemes},
J. Algebraic Geom. \textbf{12} (2003), no.~4, 605--626.
\url{https://doi.org/10.1090/S1056-3911-03-00315-1}.
Preprint: \url{https://arxiv.org/abs/math/0010215}.

\bibitem{r11CZZCusped}
R.~Canary, T.~Zhang and A.~Zimmer,
\emph{Cusped Hitchin representations and Anosov representations of geometrically finite Fuchsian groups},
Adv. Math. \textbf{404} (2022), part B, article~108439.
\url{https://doi.org/10.1016/j.aim.2022.108439}.
Author version: \url{https://arxiv.org/abs/2103.06588}.

\bibitem{mlCZZ}
R.~Canary, T.~Zhang and A.~Zimmer,
\emph{Patterson--Sullivan measures for transverse subgroups},
J. Mod. Dyn. \textbf{20} (2024), 319--377.
\url{https://doi.org/10.3934/jmd.2024009}.
Author version: \url{https://websites.umich.edu/~canary/CZZ3.pdf}.

\bibitem{r11CZZEntropy}
R.~Canary, T.~Zhang and A.~Zimmer,
\emph{Entropy rigidity for cusped Hitchin representations},
J. Topol. \textbf{19} (2026), no.~1, e70064.
\url{https://doi.org/10.1112/topo.70064}.
The Zariski-closure classification is Theorem~11.10 in the published
version and Theorem~11.9 in the author manuscript:
\url{https://websites.umich.edu/~canary/CZZ2.pdf}.

\bibitem{CCL2026}
A.~Cano, L.~A.~Cano Cordero and L.~Loeza,
\emph{Regular sets for discrete projective group actions on $\mathbb C\mathbb P^2$},
Bull. Braz. Math. Soc. (N.S.) \textbf{57} (2026), no.~2, article~28.
\url{https://doi.org/10.1007/s00574-026-00516-4}.

\bibitem{CLL}
A.~Cano, B.~Liu and M.~M.~L\'opez,
\emph{The limit set for discrete complex hyperbolic groups},
Indiana Univ. Math. J. \textbf{66} (2017), no.~3, 927--948.
\url{https://doi.org/10.1512/iumj.2017.66.6048}.
Author version: \url{https://arxiv.org/abs/1506.08113v5}.

\bibitem{CNS}
A.~Cano, J.~P.~Navarrete and J.~Seade,
\emph{Complex Kleinian Groups},
Progress in Mathematics, vol.~303, Birkh\"auser, 2013.
\url{https://doi.org/10.1007/978-3-0348-0481-3}.

\bibitem{CanoSeade2010}
A.~Cano and J.~Seade,
\emph{On the equicontinuity region of discrete subgroups of $\mathrm{PU}(1,n)$},
J. Geom. Anal. \textbf{20} (2010), no.~2, 291--305.
\url{https://doi.org/10.1007/s12220-009-9107-6}.

\bibitem{CarvajalesStecker2024}
L.~Carvajales and F.~Stecker,
\emph{Anosov representations acting on homogeneous spaces: Domains of discontinuity},
Adv. Math. \textbf{459} (2024), 110022.
\url{https://doi.org/10.1016/j.aim.2024.110022}.

\bibitem{ConzeGuivarch}
J.-P.~Conze and Y.~Guivarc'h,
\emph{Limit sets of groups of linear transformations},
Sankhy\=a Ser.~A \textbf{62} (2000), no.~3, 367--385.
\url{https://perso.univ-rennes1.fr/jean-pierre.conze/publi.htm}.

\bibitem{r11FGPositive}
V.~Fock and A.~Goncharov,
\emph{Moduli spaces of local systems and higher Teichm\"uller theory},
Publ. Math. Inst. Hautes \'Etudes Sci. \textbf{103} (2006), 1--211.
\url{https://doi.org/10.1007/s10240-006-0039-4}.

\bibitem{Frances}
C.~Frances,
\emph{Lorentzian Kleinian groups},
Comment. Math. Helv. \textbf{80} (2005), no.~4, 883--910.
\url{https://doi.org/10.4171/CMH/38}.

\bibitem{r3bor:Furstenberg1976}
H.~Furstenberg,
\emph{A note on Borel's density theorem},
Proc. Amer. Math. Soc. \textbf{55} (1976), no.~1, 209--212.
\url{https://doi.org/10.1090/S0002-9939-1976-0422497-X}.

\bibitem{Ghys1995}
\'Etienne Ghys,
\emph{D\'eformations des structures complexes sur les espaces homog\`enes de
$\mathrm{SL}(2,\C)$},
J. Reine Angew. Math. \textbf{468} (1995), 113--138.
\url{https://doi.org/10.1515/crll.1995.468.113}.

\bibitem{GuichardWienhard}
O.~Guichard and A.~Wienhard,
\emph{Anosov representations: domains of discontinuity and applications},
Invent. Math. \textbf{190} (2012), no.~2, 357--438.
\url{https://doi.org/10.1007/s00222-012-0382-7}.

\bibitem{KapovichProper}
M.~Kapovich,
\emph{A note on properly discontinuous actions},
S\~ao Paulo J. Math. Sci. \textbf{18} (2024), no.~2, 807--836.
\url{https://doi.org/10.1007/s40863-023-00353-z}.
Author version: \url{https://arxiv.org/abs/2301.05325}.

\bibitem{KLP}
M.~Kapovich, B.~Leeb and J.~Porti,
\emph{Dynamics on flag manifolds: domains of proper discontinuity and cocompactness},
Geom. Topol. \textbf{22} (2018), no.~1, 157--234.
\url{https://doi.org/10.2140/gt.2018.22.157}.

\Needspace{6\baselineskip}
\bibitem{KnappLG}
A.~W.~Knapp,
\emph{Lie Groups Beyond an Introduction}, second ed.,
Progress in Mathematics, vol.~140, Birkh\"auser, Boston, 2002.
Digital second edition, 2023: \url{https://doi.org/10.3792/euclid/9798989504206}.
Author's copy: \url{https://www.math.stonybrook.edu/~aknapp/download/Beyond2.pdf}.

\bibitem{Kulkarni}
R.~S.~Kulkarni,
\emph{Groups with domains of discontinuity},
Math. Ann. \textbf{237} (1978), 253--272.
\url{https://doi.org/10.1007/BF01420180}.

\bibitem{Kuranishi}
M.~Kuranishi,
\emph{On everywhere dense imbedding of free groups in Lie groups},
Nagoya Math. J. \textbf{2} (1951), 63--71.
\url{https://doi.org/10.1017/S0027763000010059}.

\bibitem{Mendez2020}
M.~M\'endez,
\emph{Complex orthogonal geometric structures of dimension three},
Rev. Mat. Iberoam. \textbf{36} (2020), no.~6, 1721--1752.
\url{https://doi.org/10.4171/RMI/1181}.
Author version: \url{https://arxiv.org/abs/1809.06953v1}.

\bibitem{Myrberg}
P.~J.~Myrberg,
\emph{Untersuchungen \"uber die automorphen Funktionen beliebig vieler Variablen},
Acta Math. \textbf{46} (1925), 215--336.
\url{https://doi.org/10.1007/BF02564065}.

\bibitem{Onishchik1963}
A.~L. Onishchik,
\emph{Transitive compact transformation groups},
Mat. Sb. (N.S.) \textbf{60(102)} (1963), no.~4, 447--485.
\url{https://www.mathnet.ru/eng/sm4555}.

\bibitem{qfQuint2002}
J.-F.~Quint,
\emph{Mesures de Patterson--Sullivan en rang sup\'erieur},
Geom. Funct. Anal. \textbf{12} (2002), no.~4, 776--809.
\url{https://doi.org/10.1007/s00039-002-8266-4}.
Author version: \url{https://www.math.u-bordeaux.fr/~jquint/publications/patterson.pdf}.

\bibitem{RichardsonSpringer1990}
R.~W.~Richardson and T.~A.~Springer, \emph{The Bruhat order on
symmetric varieties}, Geom. Dedicata \textbf{35} (1990), 389--436.
\url{https://doi.org/10.1007/BF00147354}.

\bibitem{SteckerThesis}
F.~Stecker,
\emph{Domains of discontinuity of Anosov representations in flag
manifolds and oriented flag manifolds}, PhD thesis,
Universit\"at Heidelberg, 2019.
\url{https://doi.org/10.11588/heidok.00026825}.
Author's copy: \url{https://florianstecker.de/files/thesis-onesided.pdf}.

\bibitem{SteckerBalanced}
F.~Stecker,
\emph{Balanced ideals and domains of discontinuity of Anosov representations},
Geom. Dedicata \textbf{218} (2024), article~9.
\url{https://doi.org/10.1007/s10711-023-00850-7}.
Author version (cited numbering): \url{https://arxiv.org/abs/1810.11496v1}.

\bibitem{ThomasGIT}
R.~P.~Thomas,
\emph{Notes on GIT and symplectic reduction for bundles and varieties},
Surveys in Differential Geometry \textbf{10} (2006), 221--273.
\url{https://doi.org/10.4310/SDG.2005.v10.n1.a7}.
Author version: \url{https://arxiv.org/abs/math/0512411}.

\bibitem{Tits1971}
J.~Tits, \emph{Repr\'esentations lin\'eaires irr\'eductibles d'un groupe
r\'eductif sur un corps quelconque}, J. Reine Angew. Math. \textbf{247} (1971), 196--220.
\url{https://doi.org/10.1515/crll.1971.247.196}.
\end{thebibliography}
\end{document}